\documentclass[11pt,letterpaper]{amsart}

\usepackage[centertags]{amsmath}
\usepackage{amsthm}
\usepackage{amssymb}
\usepackage{amscd}
\usepackage{eucal}
\usepackage{epsfig}
\usepackage[matrix,arrow,curve]{xy}
\CompileMatrices

\theoremstyle{plain}

\newtheorem{thm}{Theorem}[section]
\newtheorem{theorem}[thm]{Theorem}

\newtheorem{lemma}[thm]{Lemma}

\newtheorem{prop}[thm]{Proposition}

\newtheorem{conjecture}[thm]{Conjecture}
\newtheorem*{conjecture*}{Conjecture}

\newtheorem*{question*}{Question}
\newtheorem{defn}[thm]{Definition}
\newtheorem{definition}[thm]{Definition}

\newtheorem*{definitions*}{Definitions}
\newtheorem*{warning*}{Warning}

\theoremstyle{definition}

\newtheorem*{rem*}{Remark}
\newtheorem{remark}[thm]{Remark}
\newtheorem*{remark*}{Remark}

\newtheorem*{remarks*}{Remarks}
\newtheorem*{example*}{Example}
\newtheorem{example}[thm]{Example}
\newtheorem*{examples*}{Examples}
\newtheorem{examples}[thm]{Examples}

\newtheorem*{notation*}{Notation}
\newtheorem*{convention*}{Convention}
\newtheorem*{conventions*}{Conventions}

\newcommand{\zeroindent}{\parindent0cm \parskip1ex}

\newenvironment{theoremlist}%
{\begin{list}{{\rm(\alph{enumi}) }}{\usecounter{enumi}

\leftmargin0cm \labelsep0cm \rightmargin0cm \parsep1em \listparindent0em
\itemsep0em \topsep1em
\parskip1em \setlength{\labelwidth}{\fill}}}
{\parskip0em \end{list}}

\newenvironment{normallist}%
{

\begin{enumerate}}{\end{enumerate}}

{

\begin{enumerate}}{\end{enumerate}}

{

\begin{enumerate} \parsep0cm
\itemsep0cm \parskip0cm}{\end{enumerate}}

\newenvironment{mylist}%
{\begin{list}{{\rm(\alph{enumi}) }}{\usecounter{enumi}

\itemindent0em \parsep0em \listparindent0em
\itemsep0em \topsep0.5em \labelwidth2em
\parskip0em}}
{\end{list}}

{\begin{list}{(\alph{enumi}) }{\usecounter{enumi}

\leftmargin0cm \labelsep0cm
\itemsep0em \rightmargin0cm
\setlength{\labelwidth}{\fill}}} {\end{list}}

\newcommand{\includefigure}[3]{%
\begin{figure}[#3]
\begin{center}
\epsfig{file=#2} \\ \caption{\label{fig:#1}}
\end{center}
\end{figure}}

\newcommand{\R}{\mathbb{R}}
\newcommand{\Z}{\mathbb{Z}}

\newcommand{\C}{\mathbb{C}}

\newcommand{\half}{{\textstyle\frac{1}{2}}}

\newcommand{\iso}{\cong}           %isomorphism sign
\newcommand{\htp}{\simeq}          %homotopy sign
\newcommand{\smooth}{C^\infty}

\newcommand{\RP}[1]{\R {\mathrm P}^{#1}}
\newcommand{\leftsc}{\langle}
\newcommand{\rightsc}{\rangle}

\newcommand{\suchthat}{\; | \;}
\newcommand{\id}{\mathrm{id}}
\newcommand{\Id}{\mathrm{Id}}

\newcommand{\im}{\mathrm{im}}
\renewcommand{\ker}{\mathrm{ker}}

\newcommand{\Hom}{\mathrm{Hom}}
\newcommand{\End}{\mathrm{End}}

\renewcommand{\o}{\omega}
\newcommand{\Aut}{\mathrm{Aut}}

\newcommand{\Diff}{\mathrm{Diff}}

\newcommand{\Coh}{Coh}
\renewcommand{\O}{{\mathcal O}}
\newcommand{\Pic}{\mathrm{Pic}}
\newcommand{\Qco}{Qco}
\newcommand{\SHom}{{\mathcal H}\!\text{\em om}}
\newcommand{\SExt}{{\mathcal E}\text{\em xt}}

\newcommand{\Ext}{\mathrm{Ext}}
\newcommand{\Cone}{\mathrm{Cone}}
\newcommand{\Auteq}{\mathrm{Auteq}}

\newcommand{\ab}{\mathfrak{S}}
\newcommand{\epi}{\twoheadrightarrow}

\newcommand{\Amod}{\overline{A}_m\text{-}mod}
\newcommand{\amod}{A_{m,n}\text{-}mod}

\newcommand{\dga}{\mathcal{A}}
\newcommand{\dgm}{\mathcal{M}}

\newcommand{\Symp}{\mathrm{Symp}}
\newcommand{\Symph}{\Symp^{h}}
\newcommand{\Sympgrh}{\Symp^{h,gr}}
\newcommand{\tL}{\tilde{L}}

\newcommand{\Fuk}{Fuk}

\newcommand{\Sympgr}{\Symp^{\mathrm{gr}}}
\newcommand{\ttau}{\tilde{\tau}}

\newcommand{\isotopic}{\htp}
\newcommand{\Aff}{\mathrm{Aff}}

\newcommand{\RR}{\mathbf R\hspace{1pt}}
\newcommand{\barpath}[1]{(\overline{#1})}

\zeroindent
\numberwithin{equation}{section}

\title[Braid group actions]{Braid group actions on derived categories\\of coherent sheaves}
\author{Paul Seidel and Richard Thomas}
\date{\today}

\begin{document}
\maketitle

\begin{quote}
{\small
{\bf Abstract.} This paper gives a construction of braid group actions on the
derived category of coherent sheaves on a variety $X$. The motivation for this is Kontsevich's
homological mirror conjecture, together with the occurrence of certain braid group actions in
symplectic geometry. One of the main results is that when $\dim X \geq 2$, our braid group
actions are always faithful.

We describe conjectural mirror symmetries between smoothings and resolutions
of singularities that lead us to find examples
of braid group actions arising from crepant resolutions of various
singularities. Relations with the McKay correspondence and
with exceptional sheaves on Fano manifolds are given. Moreover, the case of an elliptic
curve is worked out in some detail. }
\end{quote}

\section{Introduction}

\subsection{Derived categories of coherent sheaves\label{subsec:introduction}}

Let $X$ be a smooth complex projective variety and $D^b(X)$ the bounded derived category of
coherent sheaves. It is an interesting question how much information about $X$ is contained in
$D^b(X)$.

Certain invariants of $X$ can be shown to depend only on $D^b(X)$. This is obviously true for
$K(X)$, the Grothendieck group of both the abelian category $\Coh(X)$ of coherent sheaves
and of $D^b(X)$. A deep result of Orlov \cite{orlov96} implies that the
topological $K$-theory $K^*_{\mathrm{top}}(X)$ is also an invariant of $D^b(X)$; hence, so are
the sums of its even and odd Betti numbers. Because of the uniqueness of Serre functors
\cite{bondal-kapranov90}, the dimension of $X$ and whether it is Calabi-Yau ($\o_X \iso \O_X$)
or not, can be read off from $D^b(X)$. Using Orlov's theorem quoted above, one can prove that
the Hochschild cohomology of $X$, $HH^*(X) = \Ext^*_{X \times X}(\O_\Delta,\O_\Delta)$,
depends only on $D^b(X)$. As pointed out by Kontsevich \cite[p.\ 131]{kontsevich94}, it is
implicit in a paper of Gerstenhaber and Schack \cite{gerstenhaber-schack90} that
\[
HH^r(X) \iso \bigoplus_{p+q = r} H^p(X,\Lambda^q TX).
\]
Thus for Calabi-Yau varieties $\dim HH^r(X) = \sum_{p+q = r} h^{p,n-q}(X)$; in mirror symmetry these
are the Betti numbers of the mirror manifold. Finally, a theorem of Bondal and Orlov
\cite{bondal-orlov97} says that if the canonical sheaf $\o_X$ or its inverse is ample, $X$ can
be entirely reconstructed from $D^b(X)$. Contrary to what this list of results might suggest,
there are in fact non-isomorphic varieties with equivalent derived categories. The first
examples are due to Mukai: abelian varieties \cite{mukai81} and $K3$ surfaces \cite{mukai87}.
Examples with nontrivial $\o_X$ have been found by Bondal and Orlov \cite{bondal-orlov95}.

This paper is concerned with a closely related object, the self-equivalence group
$\Auteq(D^b(X))$.  Recall that an exact functor between two triangulated categories
$\mathfrak{C}, \mathfrak{D}$ is a pair $(F,\nu_F)$ consisting of a functor $F: \mathfrak{C}
\longrightarrow \mathfrak{D}$ and a natural isomorphism $\nu_F: F \circ [1]_{\mathfrak{C}}
\iso [1]_{\mathfrak{D}} \circ F$ (here $[1]_{\mathfrak{C}},[1]_{\mathfrak{D}}$ are the
translation functors) with the property that exact triangles in $\mathfrak{C}$ are mapped to
exact triangles in $\mathfrak{D}$. The appropriate equivalence relation between such functors
is `graded natural isomorphism' which means natural isomorphism compatible with the maps
$\nu_F$ \cite[section 1]{bondal-orlov97}. Ignoring set-theoretic difficulties, which are
irrelevant for $\mathfrak{C} = D^b(X)$, the equivalence classes of exact functors from
$\mathfrak{C}$ to itself form a monoid. $\Auteq(\mathfrak{C})$ is defined as the group of
invertible elements in this monoid. Known results about $\Auteq(D^b(X))$ parallel those for
$D^b(X)$ itself. It always contains a subgroup $A(X) \iso (\Aut(X) \ltimes \Pic(X)) \times \Z$
generated by the automorphisms of $X$, the functors of tensoring with an invertible sheaf, and
the translation. Bondal and Orlov \cite{bondal-orlov97} have shown that if $\o_X$ or
$\o_X^{-1}$ is ample then $\Auteq(D^b(X)) = A(X)$. Mukai's arguments \cite{mukai81} imply that
$\Auteq(D^b(X))$ is bigger than $A(X)$ for all abelian varieties (recent work of Orlov
\cite{orlov97} describes $\Auteq(D^b(X))$ completely in this case).

Our own interest in self-equivalence groups comes from Kontsevich's homological mirror
conjecture \cite{kontsevich94}. One consequence of this conjecture is that for Calabi-Yau
varieties to which mirror symmetry applies, the group $\Auteq(D^b(X))$ should be related to
the symplectic automorphisms of the mirror manifold. This conjectural relationship is rather
abstract, and difficult to spell out in concrete examples. Nevertheless, as a first and rather
naive check, one can look at some special symplectic automorphisms of the mirror and try to
guess the corresponding self-equivalences of $D^b(X)$. Having made this guess in a
sufficiently plausible way (which means that the two objects show similar behaviour), the next
step might be to take some unsolved questions about symplectic automorphisms and translate it
into one about $\Auteq(D^b(X))$.  Using the smoother machinery of sheaf theory one stands a
good chance of solving this analogue, and this in turn provides a conjectural answer, or
`mirror symmetry prediction', for the original problem. The present paper is an experiment in
this mode of thinking. We now state the main results independently of their motivation; the
discussion of mirror symmetry will be taken up again in the next section.

Let $X,Y$ be two (as before, smooth complex projective) varieties. The Fourier-Mukai transform
(FMT) by an object ${\mathcal P} \in D^b(X \times Y)$ is the exact functor
\[
\Phi_{\mathcal P}: D^b(X) \longrightarrow D^b(Y), \quad
\Phi_\mathcal{P}(\mathcal{G}) = {\mathbf R}\pi_{2\,*}(
\pi_1^*\mathcal{G} \stackrel{\mathbf L}{\otimes} \mathcal{P}),
\]
where $\pi_1: X \times Y \rightarrow X$, $\pi_2: X \times Y \rightarrow Y$ are the
projections. This is a very practical way of defining functors. Orlov \cite{orlov96} has
proved that any equivalence $D^b(X) \longrightarrow D^b(Y)$ can be written as a FMT. Earlier
work of Maciocia \cite{maciocia96} shows that if $\Phi_{\mathcal P}$ is an equivalence, then
${\mathcal P}$ must satisfy a partial Calabi-Yau condition: $\mathcal P \otimes \pi_1^* \o_X
\otimes \pi_2^* \o_Y^{-1} \iso \mathcal P$. Bridgeland \cite{bridgeland98} provides a partial
converse to this.

Now take an object ${\mathcal E} \in D^b(X)$ which is a complex of locally free sheaves. We
define the {\em twist functor} $T_{\mathcal E}: D^b(X) \longrightarrow D^b(X)$ as the FMT with
\begin{equation} \label{eq:p-sheaf}
{\mathcal P} = \Cone(\eta: \mathcal{E^\vee} \boxtimes \mathcal{E} \longrightarrow \O_\Delta),
\end{equation}
where ${\mathcal E}^\vee$ is the dual complex, $\boxtimes$ the exterior tensor product,
$\Delta \subset X \times X$ is the diagonal, and $\eta$ the canonical pairing. Since
quasi-isomorphic ${\mathcal E}$ give rise to isomorphic functors $T_{\mathcal E}$, one can use
locally free resolutions to extend the definition to arbitrary objects of $D^b(X)$.

\begin{definition} \label{def:spherical-sheaf}
\begin{theoremlist}
\item
$\mathcal{E} \in D^b(X)$ is called {\em spherical} if $\Hom^r_{D^b(X)}(\mathcal E,\mathcal E)$
is equal to $\C$ for $r = 0,\, \dim X$ and zero in all other degrees, and if in addition
${\mathcal E} \otimes \o_X \iso {\mathcal E}$.

\item
An $(A_m)$-configuration, $m \geq 1$, in $D^b(X)$ is a collection of $m$ spherical objects
$\mathcal{E}_1, \dots, \mathcal{E}_m$ such that
\[
\dim_\C \Hom^*_{D^b(X)}(\mathcal E_i,\mathcal E_j) =
\begin{cases}
1 & |i-j| = 1,\\ 0 & |i-j| \geq 2.
\end{cases}
\]
\end{theoremlist}
\end{definition}

Here, as elsewhere in the paper, $\Hom^r(\mathcal E,\mathcal F)$ stands for $\Hom(\mathcal E,
\mathcal F[r])$, and $\Hom^*(\mathcal E,\mathcal F)$ is the total space $\bigoplus_{r \in \Z}
\Hom^r(\mathcal E,\mathcal F)$.

\begin{theorem} \label{th:braiding}
The twist $T_{\mathcal E}$ along any spherical object $\mathcal{E}$ is an exact
self-equivalence of $D^b(X)$. Moreover, if $\mathcal{E}_1, \dots, \mathcal{E}_m$ is an
$(A_m)$-configuration, the twists $T_{\mathcal E_i}$ satisfy the braid relations up to graded
natural isomorphism:
\begin{align*}
 & T_{\mathcal E_i} T_{\mathcal E_{i+1}} T_{\mathcal E_i} \iso
 T_{\mathcal E_{i+1}} T_{\mathcal E_i} T_{\mathcal E_{i+1}} &&
 \text{for } i = 1,\dots,m-1, \\
 & T_{\mathcal E_i} T_{\mathcal E_j} \iso T_{\mathcal E_j} T_{\mathcal E_i}
 &&\text{for } |i-j| \geq 2.
\end{align*}
\end{theorem}

We should point out that the first part, the invertibility of $T_{\mathcal E}$, was also known
to Kontsevich, Bridgeland and Maciocia. Let $\rho$ be the homomorphism from the braid group
$B_{m+1}$ to $\Auteq(D^b(X))$ defined by sending the standard generators $g_1,\dots,g_m \in
B_{m+1}$ to $T_{\mathcal E_1}, \dots, T_{\mathcal E_m}$. We call this a weak braid group
action on $D^b(X)$ (there is a better notion of a group action on a category which requires
the presence of certain additional natural transformations \cite{deligne97}; we have not
checked whether these exist in our case). $\rho$ induces a representation $\rho_*$ of
$B_{m+1}$ on $K(X)$. Concretely, the twist along an arbitrary $\mathcal{E} \in D^b(X)$ acts on
$K(X)$ by
\begin{equation} \label{eq:k-theory}
(T_{\mathcal E})_*(y) = y - \leftsc [\mathcal{E}], y \rightsc
[\mathcal{E}],
\end{equation}
where $\leftsc [\mathcal F],[\mathcal G] \rightsc = \sum_i (-1)^i \dim
\Hom^i(\mathcal{F},\mathcal{G})$ is the Mukai pairing \cite{mukai87}
or `Euler form'. If $\dim X$ is even then
$\rho_*$ factors through the symmetric group $S_{m+1}$. The odd-dimensional case is slightly
more complicated, but still $\rho_*$ is far from being faithful, at least if $m$ is large.

For $\rho$ itself we have the following contrasting result:

\begin{theorem} \label{th:faithful}
Assume that $\dim X \geq 2$. Then the homomorphism $\rho$ generated by the twists in any
$(A_m)$-configuration is injective.
\end{theorem}

The assumption $\dim X \neq 1$ cannot be removed; indeed, there is a $B_4$-action on the
derived category of an elliptic curve which is not faithful (see section
\ref{subsec:elliptic-curve}).

\subsection{Homological mirror symmetry and self-equivalences\label{subsec:mirror}}

We begin by recalling Kontsevich's homological mirror conjecture \cite{kontsevich94}. On one
hand, one takes Calabi-Yau varieties $X$ and their derived categories $D^b(X)$. On the other
hand, using entirely different techniques, it is thought that one can attach to any compact
symplectic manifold $(M,\beta)$ with zero first Chern class a triangulated category, the
derived Fukaya category $D^b\Fuk(M,\beta)$ (despite the notation, this is not constructed as
the derived category of an abelian category). Kontsevich's conjecture is that {\em whenever
$X$ and $(M,\beta)$ form a mirror pair, there is a (non-canonical) exact equivalence}
\begin{equation} \label{eq:kontsevich}
D^b(X) \iso D^b\Fuk(M,\beta).
\end{equation}
A more prudent formulation would be to say that \eqref{eq:kontsevich} should hold for the
generally accepted constructions of mirror manifolds. Before discussing this conjecture
further, we need to explain what $D^b\Fuk(M,\beta)$ looks like. This is necessarily a
tentative description, since a rigorous definition does not exist yet. Moreover, for
simplicity we have omitted some of the more technical aspects.

Let $(M,\beta)$ be as before, of real dimension $2n$. To simplify things we assume that
$\pi_1(M)$ is trivial; this excludes the case of the two-torus, so that $n \geq 2$. Recall
that a submanifold $L^n \subset M$ is called Lagrangian if $\beta|L \in \Omega^2(L)$ is zero.
Following Kontsevich \cite[p.\ 134]{kontsevich94} one considers objects, denoted by $\tL$,
which are Lagrangian submanifolds with some extra structure. We will call such objects `graded
Lagrangian submanifolds' and the extra structure the `grading'. This grading amounts
approximately to an integer choice. In fact there is a free $\Z$-action, denoted by $\tL
\mapsto \tL[j]$ for $j \in \Z$, on the set of graded Lagrangian submanifolds; and if $L$ is a
connected Lagrangian submanifold, all its possible gradings (assuming that there are any) form
a single orbit of this action. For details we refer to \cite{seidel99}. For any pair
$(\tL_1,\tL_2)$ of graded Lagrangian submanifolds one expects to have a Floer cohomology group
$HF^*(\tL_1,\tL_2)$, which is a finite-dimensional graded $\R$-vector space satisfying
$HF^*(\tL_1,\tL_2[j]) = HF^*(\tL_1[-j],\tL_2) = HF^{*+j}(\tL_1,\tL_2)$. Defining this is a
difficult problem; a fairly general solution has been announced recently by Fukaya,
Kontsevich, Oh, Ohta and Ono.

The most essential property of $D^b\Fuk(M,\beta)$ is that {\em any graded Lagrangian
submanifold $\tL$ defines an object in this category}. The translation functor (which is part
of the structure of $D^b\Fuk(M,\beta)$ as a triangulated category) acts on such objects by
$\tL \mapsto \tL[1]$. The morphisms between two objects of this kind are given by the degree
zero Floer cohomology with complex coefficients:
\[
\Hom_{D^b\Fuk(M,\beta)}(\tL_1,\tL_2) = HF^0(\tL_1,\tL_2) \otimes_{\R} \C
\]
(Floer groups in other degrees can be recovered by changing $\tL_2$ to $\tL_2[j]$).
Composition of such morphisms is given by certain products on Floer cohomology, which were
first introduced by Donaldson. There is also a slight generalisation of this: any pair
$(\tL,E)$ consisting of a graded Lagrangian submanifold together with a flat unitary vector
bundle $E$ on the underlying Lagrangian submanifold, defines an object of $D^b\Fuk(M,\beta)$.
The morphisms between such objects are a twisted version of Floer cohomology. It is important
to keep in mind that {\em $D^b\Fuk(M,\beta)$ contains many objects other than those which we
have described}. This is necessarily so because it is triangulated: there must be enough
objects to complete each morphism to an exact triangle, and these objects will not usually
have a direct geometric meaning. However, it is expected that the objects of the form
$(\tL,E)$ generate the category $D^b\Fuk(M,\beta)$ in some sense.

\begin{remark} \label{remark:imaginary}
In the traditional picture of mirror symmetry, $M$ carries a $\C$-valued closed two-form
$\beta_\C$ with real part $\beta$. What we have said concerns the Fukaya category for
$\im(\beta_\C) = 0$. Apparently, the natural generalisation to $\im(\beta_\C) \neq 0$ would be
to take objects $(\tL,E,A)$ consisting of a graded Lagrangian submanifold $\tL$, a complex
vector bundle $E$ on the underlying Lagrangian submanifold $L$, and a unitary connection $A$
on $E$ with curvature $F_A = -\beta_\C|L \otimes \id_E$. The point is that to any map $w:
(D^2, \partial D^2) \longrightarrow (M,L)$ one can associate a complex number
\[
\frac{\mathrm{trace}(\text{monodromy of $A$ around $w|\partial
D^2$})}{\mathrm{rank}(E)}\,\exp(-\textstyle{\int_{D^2}} w^*\beta_{\C}),
\]
which is invariant under deformations of $w$. These numbers, as well as certain variations of
them, would be used as weights in the counting procedure which underlies the definition of
Floer cohomology. For simplicity, we will stick to the case $\im(\beta_\C) = 0$ in our
discussion.
\end{remark}

In parallel with graded Lagrangian submanifolds, there is also a notion of graded symplectic
automorphisms; in fact these are just a special kind of graded Lagrangian submanifolds on
$(M,-\beta) \times (M,\beta)$. The graded symplectic automorphisms form a topological group
$\Sympgr(M,\beta)$ which is a central extension of the usual symplectic automorphism group
$\Symp(M,\beta)$ by $\Z$. $\Sympgr(M,\beta)$ acts naturally on the set of graded Lagrangian
submanifolds. Moreover, the central subgroup $\Z$ is generated by a graded symplectic
automorphism denoted by $[1]$, which maps each graded Lagrangian submanifold $\tL$ to
$\tL[1]$; we refer again to \cite{seidel99} for details. Because $D^b\Fuk(M,\beta)$ is defined
in what are essentially symplectic terms, {\em every graded symplectic automorphism of $M$
induces an exact self-equivalence of it}. Moreover, an isotopy of graded symplectic
automorphisms will give rise to an equivalence between the induced functors. Thus one has a
canonical map
\[
\pi_0(\Sympgr(M,\beta)) \longrightarrow \Auteq(D^b\Fuk(M,\beta)).
\]
Now we return to Kontsevich's conjecture. Assume that $(M,\beta)$ has
a mirror partner $X$ such that \eqref{eq:kontsevich} holds. Then there
is an isomorphism between $\Auteq(D^b\Fuk(M,\beta))$ and
$\Auteq(D^b(X))$. Combining this with the canonical map above yields a
homomorphism
\begin{equation} \label{eq:mu}
\mu: \pi_0(\Sympgr(M,\beta)) \longrightarrow \Auteq(D^b(X)).
\end{equation}
Somewhat oversimplified, and ignoring the conjectural nature of the
whole discussion, one can say that {\em symplectic automorphisms of
$M$ induce self-equivalences of the derived category of coherent
sheaves on its mirror partner}. Note that the map $\mu$ depends on the
choice of equivalence \eqref{eq:kontsevich} and hence is not canonical.

\begin{remark}
One can see rather easily that the central element $[1] \in
\Sympgr(M,\beta)$ induces the translation functor on
$D^b\Fuk(M,\beta)$ and hence on $D^b(X)$. Passing to the quotient
yields a map
\[
\bar{\mu}: \pi_0(\Symp(M,\beta)) \longrightarrow \Auteq(D^b(X))/
\mathrm{(translations)}.
\]
This simplified version may be more convenient for those readers who are unfamiliar with the
`graded symplectic' machinery.
\end{remark}

\subsection{Dehn twists and mirror symmetry\label{subsec:dehn}}

A Lagrangian sphere in $(M,\beta)$ is a Lagrangian submanifold $S \subset M$ which is
diffeomorphic to $S^n$. One can associate to any Lagrangian sphere a symplectic automorphism
$\tau_S$ called the generalized Dehn twist along $S$, which is defined by a local construction
in a neighbourhood of $S$ (see \cite{seidel98b} or \cite{seidel99} for details; strictly
speaking, $\tau_S$ depends on various choices, but since the induced functor on
$D^b\Fuk(M,\beta)$ is expected to be independent of these choices, we will ignore them in our
discussion). These maps are symplectic versions of the classical Picard-Lefschetz
transformations. In particular, their action on $H_*(M)$ is given by
\begin{equation} \label{eq:dehn-picard-lefschetz}
(\tau_S)_*(x) =
 \begin{cases}
 x - ([S] \cdot x) [S] & \text{if $\dim(x) = n$,} \\
 x & \text{otherwise.}
 \end{cases}
\end{equation}
where $\cdot$ is the intersection pairing twisted by a dimension-dependent sign. As explained
in \cite[section 5b]{seidel99} $\tau_S$ has a preferred lift $\ttau_S \in \Sympgr(M,\beta)$ to
the graded symplectic automorphism group. Suppose that $(M,\beta)$ has a mirror partner $X$
such that Kontsevich's conjecture \eqref{eq:kontsevich} holds. Choose some lift $\tilde{S}$ of
$S$ to a graded Lagrangian submanifold, and let ${\mathcal E} \in D^b(X)$ be the object which
corresponds to $\tilde{S}$. Then
\begin{equation} \label{eq:hom-hom}
 \Hom^*_{D^b(X)}(\mathcal{E},\mathcal{E}) \iso
 \Hom_{D^b\Fuk(M,\beta)}^*(\tilde{S},\tilde{S}) =
 HF^*(\tilde{S},\tilde{S}) \otimes_{\R} \C.
\end{equation}
The Floer cohomology $HF^*(\tilde{S},\tilde{S})$ is isomorphic to the ordinary cohomology
$H^*(S;\R)$; this is not true for general Lagrangian submanifolds, but it holds for spheres.
Therefore ${\mathcal E}$ must be spherical object (this motivated our use of the word
spherical). A natural conjecture about the homomorphism $\mu$ introduced in the previous
section is that
\begin{equation} \label{eq:mu-claim}
\mu([\ttau_S]) = [T_{\mathcal E}],
\end{equation}
where $T_{\mathcal E}$ is the twist functor as defined in section \ref{subsec:introduction}.
Roughly speaking, the idea is {\em twist functors and generalized Dehn twists correspond to
each other under mirror symmetry}. At present this is merely a guess, which can be motivated
e.g.\ by comparing \eqref{eq:k-theory} with \eqref{eq:dehn-picard-lefschetz}. But supposing
that one wanted to actually prove this claim, how should one go about it? The first step would
be to observe that for any ${\mathcal F} \in D^b(X)$ there is an exact triangle
\[
\xymatrix{ {\Hom^*(\mathcal E,\mathcal F) \otimes_{\C} \mathcal E}
\ar[r] & {\mathcal F} \ar[r] & {T_{\mathcal E}(\mathcal F)}
\ar@/^1.5em/[ll]+<0ex,-1.5ex>^-{[1]}}
\]
Here $\Hom^*(\mathcal E,\mathcal F)$ is the graded group of homs in the derived category,
$\Hom^*(\mathcal E,\mathcal F) \otimes_{\C} \mathcal E$ is the corresponding direct sum of
shifted copies of $\mathcal E$, and the first arrow is the evaluation map. This exact triangle
determines $T_{\mathcal E}(\mathcal F)$ up to isomorphism; moreover, it does so in purely
abstract terms, which involve only the triangulated structure of the category $D^b(X)$. Hence
{\em if there was an analogous abstract description of the action of $\ttau_S$ on
$D^b\Fuk(M,\beta)$ one could indeed prove \eqref{eq:mu-claim}} (this is slightly imprecise,
since it ignores a technical problem about non-functoriality of cones in triangulated
categories). The first step towards such a description will be provided in \cite{seidel99b}.
Note that here, for the first time in our discussion of mirror symmetry, we have made
essential use of the triangulated structure of the categories.

Now define an $(A_m)$-configuration of Lagrangian spheres in $(M,\beta)$ to be a collection of
$m \geq 1$ pairwise transverse Lagrangian spheres $S_1,\dots,S_m \subset M$ such that
\begin{equation} \label{eq:intersect}
|S_i \cap S_j| = \begin{cases} 1 & |i-j| = 1, \\ 0 & |i-j| \geq
 2. \end{cases}
\end{equation}
Such configurations occur in K{\"a}hler manifolds that can be degenerated into
a manifold with
a singular point of type $(A_m)$ (see \cite{seidel98b} or
\cite{khovanov-seidel98}). The
generalized Dehn twists $\ttau_{S_1},\dots,\ttau_{S_m}$ along such spheres
satisfy the braid
relations up to isotopy inside $\Sympgr(M,\beta)$. For $n=2$, and ignoring
the issue of gradings, this was proved in \cite[Appendix]{seidel98b}; the argument given there
can be adapted to yield the slightly sharper and more general statement which we are using
here. Thus, by mapping the standard generators of the braid group to the classes
$[\ttau_{S_i}]$ one obtains a homomorphism from $B_{m+1}$ to $\pi_0(\Sympgr(M,\beta))$. It is
a difficult open question in symplectic geometry whether this homomorphism, which we denote by
$\rho'$, is injective; see \cite{khovanov-seidel98} for a partial result. We will now see what
mirror symmetry has to say about this.

Assume as before that Kontsevich's conjecture holds, and let ${\mathcal E}_1,\dots,{\mathcal
E}_m \in D^b(X)$ be the objects corresponding to some choice of gradings
$\tilde{S}_1,\dots,\tilde{S}_m$ for the $S_j$. We already know that each ${\mathcal E}_i$ is a
spherical object. An argument similar to \eqref{eq:hom-hom} but based on \eqref{eq:intersect}
shows that ${\mathcal E}_1,\dots,{\mathcal E}_m$ is an $(A_m)$-configuration in $D^b(X)$ in
the sense of Definition \ref{def:spherical-sheaf}.  Hence the twist functors $T_{{\mathcal
E}_i}$ satisfy the braid relations (Theorem \ref{th:braiding}) and generate a homomorphism
$\rho$ from $B_{m+1}$ to $\Auteq(D^b(X))$. Assuming that our claim \eqref{eq:mu-claim} is
true, one would have a commutative diagram
\[
\xymatrix{{B_{m+1}} \ar^-{\rho'}[r] \ar_-{\rho}[dr] &
{\pi_0(\Sympgr(M,\beta))} \ar[d]^{\mu} \\ & {\Auteq(D^b(X))}}
\]
Since $\dim_{\C} X = n \geq 2$, we have Theorem \ref{th:faithful} which says that $\rho$ is
injective. In the diagram above this would clearly imply that $\rho'$ is injective. Thus we
are led to a conjectural answer `based on mirror symmetry' to a question of symplectic
geometry:

\begin{conjecture}
Let $(M,\beta)$ be a compact symplectic manifold with $\pi_1(M)$ trivial and $c_1(M,\beta) =
0$, and $(S_1,\dots,S_m)$ an $(A_m)$-configuration of Lagrangian spheres in $M$ for some $m
\geq 1$. Then the map $\rho': B_{m+1} \rightarrow \pi_0(\Sympgr(M,\beta))$ generated by the
generalized Dehn twists $\ttau_{S_1},\dots,\ttau_{S_m}$ is injective.
\end{conjecture}

\subsection{A survey of the paper\label{subsec:contents}}

Section \ref{sec:braid-group-actions} introduces spherical objects and twists functors for
derived categories of fairly general abelian categories. The main result is the construction
of braid group actions, Theorem \ref{th:braid-group-action}.

Section \ref{subsec:sheaves} explains how the abstract framework specializes in the case of
coherent sheaves; this recovers the definitions presented in section
\ref{subsec:introduction}, and in particular Theorem \ref{th:braiding}. More generally, in section \ref{subsec:generalisations}, we
consider singular and quasi-projective varieties, as well as equivariant sheaves on varieties
with a finite group action; the latter give rise to what are probably the simplest examples of
our theory. In section \ref{subsec:exceptional} we present a more systematic way of producing spherical objects, which exploits
their relations with the (much studied) exceptional objects on Fano varieties.
Elliptic curves provide the only example where both
sides of the homological mirror conjecture are completely understood; in section \ref{subsec:elliptic-curve} the
group of symplectic automorphisms and the group of autoequivalences of the derived
category are compared in an explicit way.
Section \ref{subsec:k3} gives more explicit examples on $K3$ surfaces, then finally section \ref{subsec:singularities}
puts our results in the framework of mirror symmetry for singularities; this was the underlying motivation for much of this work.

Section \ref{sec:faithfulness} contains the proof of the faithfulness result, Theorem
\ref{th:main}. For the benefit of the reader, we provide here an outline of the argument, in
the more concrete situation stated as Theorem \ref{th:faithful} above; the general case does
not differ greatly from this. Let ${\mathcal E}_1,\dots,{\mathcal E}_m$ be a collection of
spherical objects in $D^b(X)$, and set ${\mathcal E} = {\mathcal E}_1 \oplus \dots \oplus
{\mathcal E}_m$. For a fixed $m$ and dimension $n$ of the variety, the endomorphism algebra
\[
\End^*({\mathcal E}) = \bigoplus_{i,j} \Hom^*({\mathcal E}_i,{\mathcal E}_j)
\]
is essentially the same for all $({\mathcal E}_1,\dots, {\mathcal E}_m)$. More precisely,
after possibly shifting each ${\mathcal E}_i$ by some amount, one can achieve that
$\End^*({\mathcal E})$ is equal to a specific graded algebra $A_{m,n}$ depending only on
$m,n$. Moreover, one can define a functor $\Psi^{\text{naive}}: D^b(X) \longrightarrow \amod$
into the category of graded modules over $A_{m,n}$ by mapping ${\mathcal F}$ to
$\Hom^*({\mathcal E},{\mathcal F})$. By a result of \cite{khovanov-seidel98} the derived
category $D^b(\amod)$ carries a weak action of $B_{m+1}$, and one might hope that
$\Psi^{\text{naive}}$ should be compatible with these two actions. A little thought shows that
this cannot possibly be true: $\amod$ can be embedded into $D^b(\amod)$ as the subcategory of
complexes of length one, but the braid group action on $D^b(\amod)$ does not preserve this
subcategory. Nevertheless, the basic idea can be saved, at the cost of introducing some more
homological algebra.

Take resolutions ${\mathcal E}_i'$ of ${\mathcal E}_i$ by bounded below complexes of injective
quasi-coherent sheaves. Then one can define a differential graded algebra $end({\mathcal E'})$
whose cohomology is $\End^*({\mathcal E})$. The quasi-isomorphism type of $end({\mathcal E'})$
is independent of the choice of resolutions, so it is an invariant of the
$(A_m)$-configuration ${\mathcal E}_1,\dots,{\mathcal E}_m$. As before there is an exact
functor $hom({\mathcal E'},-): D^b(X) \longrightarrow D(end({\mathcal E'}))$ to the derived
category of differential graded modules over $end({\mathcal E'})$. Now assume that
$end({\mathcal E'})$ is formal, that is to say, quasi-isomorphic to the differential graded
algebra $\dga_{m,n} = (A_{m,n},0)$ with zero differential. Quasi-isomorphic differential
graded algebras have equivalent derived categories, so what one obtains is an exact functor
\[
\Psi: D^b(X) \longrightarrow D(\dga_{m,n}),
\]
which replaces the earlier $\Psi^{\text{naive}}$. A slight modification of the arguments of
\cite{khovanov-seidel98} shows that there is a weak braid group on $D(\dga_{m,n})$; moreover,
in contrast with the situation above, the functor $\Psi$ now relates the two braid group
actions. Still borrowing from \cite{khovanov-seidel98}, one can interpret the braid group
action on $D(\dga_{m,n})$ in terms of low-dimensional topology, and more precisely geometric
intersection numbers of curves on a punctured disc. This leads to a strong faithfulness result
for it, which through the functor $\Psi$ implies the faithfulness of the original braid group
action on $D^b(X)$.

This argument by reduction to the known case of $D(\dga_{m,n})$ hinges on the formality of
$end({\mathcal E}')$. We will prove that this assumption is always satisfied when $n \geq 2$.
This has nothing to do with the geometric origin of $end({\mathcal E}')$; in fact, what we
will show is that $A_{m,n}$ is intrinsically formal for $n \geq 2$, which means that all
differential graded algebras with this cohomology are formal. There is a general theory of
intrinsically formal algebras, which goes back to the work of Halperin and Stasheff
\cite{halperin-stasheff79} in the commutative case; the Hochschild cohomology computation
necessary to apply this theory to $A_{m,n}$ is the final step in the proof of Theorem
\ref{th:main}.

{\em Acknowledgments.} Although he does not figure as an author, the paper was originally
conceived jointly with Mikhail Khovanov, and several of the basic ideas are his. At an early
stage of this work, we had a stimulating conversation with Maxim Kontsevich. We would also
like to thank Mark Gross for discussions about mirror symmetry and singularities, and Brian
Conrad, Umar Salam, and Balazs Szendroi for helpful comment. As mentioned earlier, Kontsevich,
Bridgeland and Maciocia also knew about the invertibility of the twist functors. Financial
support came from Max Planck Institute (Bonn) and Hertford College (Oxford).

{\em Addendum.} The results here were first announced at the Harvard
Winter School on Mirror Symmetry in January of 1999 (published in
\cite{thomas99}). In the meantime, a preprint by Horja \cite{horja99}
has appeared which is inspired by similar mirror symmetry considerations.
While there is little actual overlap (\cite{horja99} does not operate in
the derived category) Horja uses
monodromy calculations to predict corresponding conjectural
mirror Fourier-Mukai transforms that ought to be connected to our
work, linking it to the toric construction of mirror manifolds.

\section{Braid group actions\label{sec:braid-group-actions}}

\subsection{Generalities\label{subsec:homological-algebra}}

Fix a field $k$; all categories are assumed to be $k$-linear. If $\ab$ is an abelian category,
$Ch(\ab)$ is the category of cochain complexes in $\ab$ and cochain maps, $K(\ab)$ the
corresponding homotopy category (morphisms are homotopy classes of cochain maps), and $D(\ab)$
the derived category. The variants involving bounded (below, above, or on both sides)
complexes are denoted by $Ch^+(\ab)$, $Ch^-(\ab)$, $Ch^b(\ab)$ and so on. Let
$(C_j,\delta_j)_{j \in \Z}$ be a cochain complex of objects and morphisms in $Ch(\ab)$, that
is to say $C_j \in Ch(\ab)$ and $\delta_j \in \Hom_{Ch(\ab)}(C_j,C_{j+1})$ satisfying
$\delta_{j+1}\delta_j = 0$. Such a complex is exactly the same as a bicomplex in $\ab$. In
this case we will write $\{\dots C_{-1} \rightarrow C_0 \rightarrow C_1 \dots\}$ for the
associated total complex, obtained by collapsing the bigrading; this is a single object in
$Ch(\ab)$. The same notation will be applied to bicomplexes of objects of $Ch(\ab)$ (which are
triple complexes in $\ab$).

For $C, D \in Ch(\ab)$, let $hom(C,D)$ be the standard cochain complex of $k$-vector spaces whose
cohomology is $H^i hom(C,D) = \Hom^i_{K(\ab)}(C,D)$, that is, $hom^i(C,D) = \prod_{j \in \Z}
\Hom_{\ab}(C^j,D^{j+i})$ with $d_{hom(C,D)}^i(\phi) = d_D\phi - (-1)^i \phi d_C$. Now suppose that
$\ab$ contains infinite direct sums and products. Given an object $C \in Ch(\ab)$ and a cochain
complex $b$ of $k$-vector spaces, one can form the tensor product $b \otimes C$ and the complex of
linear maps $lin(b,C)$, both of which are again objects of $Ch(\ab)$. They are defined by choosing a
basis of $b$ and taking a corresponding direct sum (for $b \otimes C$) or product (for $lin(b,C)$)
of shifted copies of $C$, with a differential which combines $d_b$ and $d_C$. The outcome is
independent of the chosen basis up to canonical isomorphism. The definition of $b \otimes C$ is
clear, but for $lin(b,C)$ there are two possible choices of signs. Ours is fixed to fit in with an
evaluation map $b \otimes lin(b,C) \longrightarrow C$. To clarify the issue we will now spell out
the definition. Take a homogeneous basis $(x_i)_{i \in I}$ of the total space $b$, and write
$d_b(x_i) = \sum_j z_{ji} x_j$. Then $lin^q(b,C) = \prod_{i \in I} C_i^q$, where $C_i$ is a copy of
$C$ shifted by $\deg(x_i)$. The differential $d^q: lin^q(b,C) \longrightarrow lin^{q+1}(b,C)$ has
components $d^q_{ji}: C^q_i \longrightarrow C^{q+1}_j$ which are given by
\[
d^q_{ji} =
\begin{cases}
(-1)^{\deg(x_i)} d_C & i = j, \\ (-1)^{\deg(x_i)} z_{ij} \cdot \id_C & \deg(x_i) = \deg(x_j) + 1, \\
0 & \text{otherwise.}
\end{cases}
\]
One can verify that the map $b \otimes lin(b,C) \longrightarrow C$, $x_j \otimes (c_i)_{i \in I}
\longmapsto c_j$, is indeed a morphism in $Ch(\ab)$. Moreover, there are canonical monomorphic
cochain maps
\begin{equation} \label{eq:trivial-morphisms}
\begin{aligned}
b \otimes hom(D,C) &\longrightarrow hom(D,b \otimes C), \\
hom(D,C) \otimes b &\longrightarrow hom(lin(b,D),C), \\
hom(B,lin(b,C)) \otimes D &\longrightarrow lin(b,hom(B,C) \otimes D),
\end{aligned}
\end{equation}
where $b$ is as before and $B,C,D \in Ch(\ab)$. These maps are isomorphisms if $b$ is
finite-dimensional, and quasi-isomorphisms if $b$ has finite-dimensional cohomology.

From now on $\ab$ will be an abelian category and $\ab' \subset \ab$ a full subcategory, such that
the following conditions hold:

\begin{mylist}
\item[(C1)]
$\ab'$ is a Serre subcategory of $\ab$ (this means that any subobject and quotient object of an
object in $\ab'$ lies again in $\ab'$, and that $\ab'$ is closed under extension);
\item[(C2)]
$\ab$ contains infinite direct sums and products;
\item[(C3)]
$\ab$ has enough injectives, and any direct sum of injectives is again injective (this is {\em not}
a trivial consequence of the definition of an injective object);
\item[(C4)]
for any epimorphism $f: A \epi A'$ with $A \in \ab$ and $A' \in \ab'$, there is a $B' \in \ab'$ and
a $g: B' \longrightarrow A$ such that $fg$ is an epimorphism (because $\ab'$ is a Serre subcategory,
$g$ may be taken to be mono):
\[
\xymatrix{ && {A} \ar@{->>}[d]^{f} \\ {B'} \ar@{^{(}-->}[rru]^{g} \ar@{-->>}[rr] && {A'}. }
\]
\end{mylist}

\begin{lemma} \label{ex:qc}
Let $X$ be a noetherian scheme over $k$ and $\ab = \Qco(X)$, $\ab' = \Coh(X)$ the categories of
quasi-coherent resp.\ coherent sheaves. Then properties (C1)--(C4) are satisfied.
\end{lemma}

\proof (C1) and (C2) are obvious. $\ab$ has enough injectives by \cite[II 7.18]{hartshorne66}.
Moreover, it is locally noetherian, which implies that direct sums of injectives are again
injective; see \cite[p.\ 121]{hartshorne66} and the references quoted there. This proves (C3).
Finally, we need to verify that a diagram as in (C4) with $A$ quasi-coherent and $A'$ coherent, can
be completed with a coherent sheaf $B'$. Such a $B'$ certainly exists locally, and replacing it by
its image in $A$ (which is also coherent) we may extend it to be a coherent subsheaf on all of $X$
(see EGA I 9.4.7). Since $X$ is quasi-compact, repeating this a finite number of times and taking
the union yields a $B'$ whose map to $A'$ is globally onto. \qed

As indicated by this example, our main interest is in $D^b(\ab')$. However we find it convenient to
replace all complexes by injective resolutions. These resolutions may exist only in $\ab$, and they
are not necessarily bounded. The precise category we want to work with is this:

\begin{defn} \label{def:k}
${\mathfrak K} \subset K^+(\ab)$ is the full subcategory whose objects are those bounded below
cochain complexes $C$ of $\ab$-injectives which satisfy $H^i(C) \in \ab'$ for all $i$, and
$H^i(C) = 0$ for $i \gg 0$.
\end{defn}

We will now prove, in several steps, that ${\mathfrak K}$ is equivalent to $D^b(\ab')$. First
of all, let ${\mathfrak D} \subset D^+(\ab)$ be the full subcategory of objects whose
cohomology has the same properties as in Definition \ref{def:k}. The assumption that $\ab$ has
enough injectives implies that the obvious functor ${\mathfrak K} \longrightarrow {\mathfrak
D}$ is an equivalence. Now let $Ch^b_{\ab'}(\ab)$ be the category of bounded cochain complexes
in $\ab$ whose cohomology objects lie in $\ab'$, and $D^b_{\ab'}(\ab)$ the corresponding full
subcategory of $D^b(\ab)$. It is a standard result (proved by truncating cochain complexes)
that the obvious functor $D^b_{\ab'}(\ab) \longrightarrow {\mathfrak D}$ is an equivalence.
The final step (and the only nontrivial one) is to relate $D^b_{\ab'}(\ab)$ and $D^b(\ab')$.

\begin{lemma} \label{th:complexes}
For any $C \in Ch_{\ab'}^b(\ab)$ there is an $E \in Ch^b(\ab')$ and a mono\-morphic cochain map
$\iota: E \longrightarrow C$ which is a quasi-isomorphism.
\end{lemma}

\proof Recall that, as an abelian category, $\ab$ has fibre products. The fibre product of two maps
$f_1: A_1 \longrightarrow A$, $f_2: A_2 \longrightarrow A$ is the kernel of $f_1 \oplus 0 - 0 \oplus
f_2:\,A_1 \oplus A_2 \longrightarrow A$. If $f_1$ is mono (thought of as an inclusion) we write
$f_2^{-1}(A_1)$ for the fibre product, and if both $f_1$ and $f_2$ are mono we write $A_1 \cap A_2$.
In the latter case one can also define the sum $A_1 + A_2$ as the image (kernel of the map to the
cokernel) of $f_1 \oplus 0 - 0 \oplus f_2$.

Let $N$ be the largest integer such that $C^N \neq 0$. Set $E^n = 0$ for all $n > N$. For $n \leq N$
define $E^n \subset C^n$ (for brevity, we write the monomorphisms as inclusions) inductively as
follows. By invoking (C4) one finds subobjects $F^n,G^n \subset C^n$ which lie in $\ab'$ and
complete the diagrams
\[
\xymatrix{ && {(d^n_C)^{-1}(E^{n+1})} \ar@{->>}[d]^{d^n_C} &&&& {\ker\;d^n_C} \ar@{->>}[d] \\ {F^n}
\ar@{^{(}-->}[rru] \ar@{-->>}[rr] && {E^{n+1} \cap \im\;d^n_C} & {\text{and}} & {G^n}
\ar@{^{(}-->}[rru] \ar@{-->>}[rr]  && {H^n(C).} }
\]
Set $E^n = F^n + G^n$ (this is again in $\ab'$) and define $d^n_E = d^n_C|E^n$. Since $E^n$ is a
subobject of $C^n$ for any $n$, $E$ is a bounded complex. Consider the obvious map $j^n: \ker\;
d_E^n = E^n \cap \ker\; d_C^n \longrightarrow H^n(C)$ . The definition of $G^n$ implies that $j^n$
is an epimorphism, and the definition of $F^{n-1}$ yields $\ker\; j^n = E^n \cap \im\;d^{n-1}_C =
\im\;d^{n-1}_E$. It follows that the inclusion induces an isomorphism $H^*(E) \iso H^*(C)$. \qed

From this Lemma it now follows by standard homological algebra \cite[Proposition III.2.10]{gelfand-manin} that
the obvious functor $D^b(\ab') \longrightarrow D^b_{\ab'}(\ab)$ is an equivalence of categories.
Combining this with the remarks made above, one gets

\begin{prop} \label{th:d-and-k}
There is an exact equivalence (canonical up to natural isomorphism) $D^b(\ab') \iso {\mathfrak
K}$. \qed
\end{prop}

%
%Note: the object lin(b,C) is defined to be (*b \otimes C). Here *b is the "right dual"
%of the complex b, in the sense that there is a natural map b \otimes *b ----> k. the
%right dual and left dual differ by a single sign; taking this sign wrong messes up
%everything. Note for instance that (*b)* = b in the finite-dimensional case, whereas
%b**$ is $b$ with the differential reversed. Note also that hom(C,D) = D \otimes C*,
%with the left dual! As another example take the second of the
%tautological equations;
%
% hom(D,C) \otimes b = C \otimes D* \otimes b = (NOW!!!) C \otimes (*b \otimes D)*
% = hom(lin(b,D),C).
%

\subsection{Twist functors and spherical objects\label{subsec:inverse}}

\begin{defn} \label{def:twists}
Let $E \in {\mathfrak K}$ be an object satisfying the following finiteness conditions:
\begin{mylist}
\item[(K1)] $E$ is a bounded complex,
\item[(K2)] for any $F \in {\mathfrak K}$, both $\Hom_{\mathfrak K}^*(E,F)$ and $\Hom^*_{\mathfrak K}(F,E)$ have finite
(total) dimension over $k$.
\end{mylist}
Then we define the {\em twist functor} $T_E: {\mathfrak K} \longrightarrow {\mathfrak K}$ by
\begin{equation} \label{eq:the-twist}
T_E(F) = \{hom(E,F) \otimes E \stackrel{ev}{\longrightarrow} F\}.
\end{equation}
\end{defn}

This expression requires some explanation. $ev$ is the obvious evaluation map. The grading is
such that if one ignores the differential, $T_E(F) = F \oplus (hom(E,F) \otimes E)[1]$. In
other words $T_E(F)$ is the cone of $ev$. Since $E$ is bounded and $F$ is bounded below,
$hom(E,F)$ is again bounded below. Hence $hom(E,F) \otimes E$ is a bounded below complex of
injectives in $\ab$ (this uses property (C3) of $\ab$). Its cohomology $H^*(hom(E,F) \otimes
E)$ is isomorphic to $\Hom_{\mathfrak K}^*(E,F) \otimes H^*(E)$ (for
instance because $hom(E,F)$ is quasi-isomorphic to $Hom_{\mathfrak
  K}^*(E,F)$, which is finite dimensional), and so
 is bounded, and the finiteness conditions
imply that each cohomology group lies in $\ab'$. Therefore $hom(E,F) \otimes E$ lies in
${\mathfrak K}$, and the same holds for $T_E(F)$. The functoriality of $T_E$ is obvious, and
one sees easily that it is an exact functor. Actually, for any $F,G \in {\mathfrak K}$ there
is a canonical map of complexes $(T_E)_*: hom(F,G) \longrightarrow hom(T_E(F), T_E(G))$. In
fancy language, this means that $T_E$ is functorial on the differential graded category which
underlies ${\mathfrak K}$.

\begin{prop} \label{th:iso-implies-eq}
If two objects $E_1,E_2 \in {\mathfrak K}$ satisfying (K1), (K2) are isomorphic, the
corresponding functors $T_{E_1},T_{E_2}$ are isomorphic.
\end{prop}

\proof Take cones of the rows of the following commutative diagram,
\[
\xymatrix{
 {hom(E_1,F) \otimes E_1} \ar[r] & {F} \\
 {hom(E_2,F) \otimes E_1} \ar[r] \ar[u] \ar[d] & {F} \ar@{=}[u] \ar@{=}[d] \\
 {hom(E_2,F) \otimes E_2} \ar[r] & {F.\!} }
\]
Here the vertical arrows are induced by a quasi-isomorpism
of complexes $E_1\to E_2$. \qed

Note also that $T_{E[j]}$ is isomorphic to $T_E$ for any $j \in \Z$.

\begin{defn}
For an object $E$ as in Definition \ref{def:twists} we define the {\em dual twist functor}
$T_E': {\mathfrak K} \longrightarrow {\mathfrak K}$ by $T_E'(F) = \{ev': F \longrightarrow
lin(hom(F,E), E)\}$.
\end{defn}

Here the grading is such that $F$ lies in degree zero. $ev'$ is again some kind of evaluation
map. To write it down explicitly, choose a homogeneous basis $(\psi_i)$ of $hom(F,E)$. Then
$lin^q(hom(F,E),E) = \prod_i E_i^q$, where $E_i$ is a copy of $E[\deg(\psi_i)]$, and the
$i$-th component of $ev'$ is simply $\psi_i$ itself. $T_E'$ is again an exact functor from
${\mathfrak K}$ to itself.

\begin{lemma} \label{th:adjointness}
$T_E'$ is left adjoint to $T_E$. \end{lemma}

\proof Using the maps from \eqref{eq:trivial-morphisms} and condition (K2) one constructs a
chain of natural (in $F,G \in {\mathfrak K}$) quasi-isomorphisms
\begin{align*}
hom(F,T_E(G)) =&\; \{hom(F, hom(E,G) \otimes E) \longrightarrow hom(F,G)\} \\
\longleftarrow&\; \{hom(E,G) \otimes hom(F,E) \longrightarrow hom(F,G)\} \\ \longrightarrow&\;
\{hom(lin(hom(F,E),E),G) \longrightarrow hom(F,G)\} \\ =&\; hom(T_E'(F),G).
\end{align*}
Here the chain map $hom(E,G) \otimes hom(F,E) \longrightarrow hom(F,G)$ is just composition.
The reader may easily check that the required diagrams commute. Taking $H^0$ on both sides
yields a natural isomorphism $\Hom_{\mathfrak K}(F,T_E(G)) \iso \Hom_{\mathfrak
K}(T_E'(F),G)$. \qed

\begin{defn} \label{def:k-spherical}
An object $E \in {\mathfrak K}$ is called {\em $n$-spherical} for some $n > 0$ if it satisfies
(K1), (K2) above and in addition,
\begin{mylist}
\item[(K3)] $\Hom^i_{\mathfrak K}(E,E)$ is equal to $k$ for $i = 0,n$ and zero in all other degrees;
\item[(K4)] The composition $\Hom^j_{\mathfrak K}(F,E) \times \Hom^{n-j}_{\mathfrak K}(E,F) \longrightarrow
\Hom^n_{\mathfrak K}(E,E) \iso k$ is a nondegenerate pairing for all $F \in {\mathfrak K}$, $j
\in \Z$.
\end{mylist}
\end{defn}

One can also define $0$-spherical objects: these are objects $E$ for which $\Hom^*_{\mathfrak
K}(E,E)$ is two-dimensional and concentrated in degree zero, and such that the pairings
$\Hom^j_{\mathfrak K}(E,F) \times \Hom^{-j}_{\mathfrak K}(F,E) \longrightarrow
\Hom^0_{\mathfrak K}(E,E)/k \cdot \id_E$ are nondegenerate (this means in particular that
$\Hom^0_{\mathfrak K}(E,E)$ is isomorphic to $k[t]/t^2$ as a $k$-algebra). We will not pursue
this further; the interested reader can easily verify that the proof of the next Proposition
extends to this case.

\begin{prop} \label{th:invertibility}
If $E$ is $n$-spherical for some $n>0$, both $T_E' T_E$ and $T_E T_E'$ are naturally
isomorphic to the identity functor $\Id_{\mathfrak K}$. In particular, $T_E$ is an exact
self-equivalence of ${\mathfrak K}$.
\end{prop}

\proof\hspace{-3mm} \footnote{We thank one of the referees for simplifying our
original proof of this result.}\hspace{2mm} $T_E T_E'(F)$ is a total complex
\begin{equation} \label{eq:first-complex}
\left\{
\begin{CD}
{hom(E,F) \otimes E} @>{\delta}>> {hom(E,lin(hom(F,E),E)) \otimes E} \\ @V{\alpha}VV
@V{\gamma}VV \\ {F} @>{\beta}>> lin(hom(F,E),E)
\end{CD}
\right\}
\end{equation}
Here $\alpha = ev$, $\beta = ev'$, $\gamma$ is a map induced by $ev$, and $\delta$ a map
induced by $ev'$. We shall need to know a little more about $\delta$. By the very
definition of $ev'$ by duality, $\delta$'s induced map on cohomology
\begin{equation} \label{eq:second-complex}
\Hom^*_{\mathfrak K}(E,F) \otimes H^*(E) \longrightarrow \Hom^*_{\mathfrak K}(F,E)^{\vee}
\otimes \Hom^*_{\mathfrak K}(E,E) \otimes H^*(E)
\end{equation}
is dual to the the composition $\Hom^*_{\mathfrak K}(F,E) \otimes
\Hom^*_{\mathfrak K}(E,F) \longrightarrow \Hom^*_{\mathfrak K}(E,E)$,
tensored with the identity map on $H^*(E)$. This second pairing is, by
the conditions (K3) and (K4) on $E$, perfect when we divide
$\Hom^*_{\mathfrak K}(E,E)$ by its degree zero piece $(k\cdot
\id_E)$. Thus the following modification of the map (\ref{eq:second-complex}),
\begin{equation} \label{eq:third-complex}
\Hom^*_{\mathfrak K}(E,F) \otimes H^*(E) \longrightarrow \Hom^*_{\mathfrak K}(F,E)^{\vee}
\otimes \frac{\Hom^*_{\mathfrak K}(E,E)}{k \cdot \id_E} \otimes H^*(E),
\end{equation}
is an isomorphism.

We now enlarge slightly the object in the top right hand
corner of (\ref{eq:first-complex}) to produce a new, quasi-isomorphic, complex $Q_E(F)$. The last
equation in \eqref{eq:trivial-morphisms} gives a map $hom(E,lin(hom(F,E),E))
\otimes E \hookrightarrow lin(hom(F,E),hom(E,E) \otimes E)$. Since $hom(F,E)$ has
finite-dimensional cohomology, this is a quasi-isomorphism. $\gamma$
extends naturally to $\bar{\gamma}: lin(hom(F,E), hom(E,E) \otimes E) \longrightarrow lin(hom(F,E),E)$; it is just the map induced by $ev: hom(E,E) \otimes E
\longrightarrow E$.

In fact $\bar{\gamma}$ splits canonically: define
the map $\phi: lin(hom(F,E),E) \longrightarrow lin(hom(F,E), hom(E, E) \otimes E$
induced by $k \longrightarrow hom(E,E)$, $1 \mapsto \id_E$. From the definition of
$\bar{\gamma}$ it follows that $\bar{\gamma} \circ \phi = \id$. This splitting gives a way of
embedding an acyclic complex $\{\id: lin(hom(F,E),E) \longrightarrow lin(hom(F,E),E)\}$ into
our enlarged complex $Q_E(F)$; the cokernel is
\[
\{hom(E,F) \otimes E \xrightarrow{\delta\oplus\alpha}
lin(hom(F,E),\frac{hom(E,E)}{k \cdot \id_E} \otimes E)\,\oplus\,F\}.
\]
There is an obvious map of $F$ to this, and everything we have done is
functorial in $F$; thus to prove that $T_E T_E' \iso \Id_{\mathfrak
  K}$ we are left with showing that the cokernel
\begin{equation}
\{hom(E,F) \otimes E \xrightarrow{\delta}
lin(hom(F,E),\frac{hom(E,E)}{k \cdot \id_E} \otimes E)\}
\end{equation}
is acyclic, i.e. the arrow induces an isomorphism on cohomology. But passing to
cohomology yields (\ref{eq:third-complex}), which we already noted was an
isomorphism.

The proof that $T_E' T_E \iso \Id_{\mathfrak K}$ is similar; one passes from
$T_E' T_E(F)$ to a quasi-isomorphic but slightly smaller object, which then has
a natural map to $F$. The details are almost the same as before, and we leave
them to the reader. \qed

\subsection{The braid relations\label{subsec:braiding}}

\begin{lemma} \label{th:braid-zero}
Let $E_1,E_2 \in {\mathfrak K}$ be two objects such that $E_1$ satisfies the conditions (K1),
(K2) of Definition \ref{def:twists}, and $E_2$ is $n$-spherical for some $n>0$. Then
$T_{E_2}(E_1)$ also satisfies (K1), (K2) and $T_{E_2}T_{E_1}$ is naturally isomorphic to
$T_{T_{E_2}(E_1)}T_{E_2}$.
\end{lemma}

\proof Since $E_1$ and $E_2$ are bounded complexes, so are $hom(E_1,E_2)$ and $T_{E_2}(E_1)$.
Lemma \ref{th:adjointness} says that $\Hom^*_{\mathfrak K}(F,T_{E_2}(E_1)) \iso
\Hom^*_{\mathfrak K}(T_{E_2}'(F),E_1)$. By assumption on $E_1$, this implies that
$\Hom^*_{\mathfrak K}(F,T_{E_2}(E_1))$ is always finite-dimensional. Similarly, the
finite-dimensionality of $\Hom^*_{\mathfrak K}(T_{E_2}(E_1),F)$ follows from Proposition
\ref{th:invertibility} since $\Hom^*_{\mathfrak K}(T_{E_2}(E_1),F) \iso \Hom^*_{\mathfrak
K}(E_1,T_{E_2}'(F))$. We have now proved that $T_{E_2}(E_1)$ satisfies (K1), (K2).
$T_{E_2}T_{E_1}(F)$ is a total complex
\[
\left\{ \begin{CD} hom(E_2,hom(E_1,F) \otimes E_1) \otimes E_2 @>>> hom(E_2,F) \otimes E_2
\\ @VVV @VVV \\ hom(E_1,F) \otimes E_1 @>>> F \end{CD} \right\}
\]
where all arrows are evaluation maps or induced by them. We will argue as in the proof of
Proposition \ref{th:invertibility}. Using \eqref{eq:trivial-morphisms} one sees that the
object in the top left hand corner can be replaced by the smaller quasi-isomorphic one
$hom(E_1,F) \otimes hom(E_2,E_1) \otimes E_2$. More precisely, this modification defines
another functor $R_{E_1,E_2}$ on ${\mathfrak K}$ which is naturally isomorphic to
$T_{E_2}T_{E_1}$. One can rewrite the definition of this functor as
\begin{equation} \label{eq:r-functor}
R_{E_1,E_2}(F) = \{ hom(E_1,F) \otimes T_{E_2}(E_1) \longrightarrow T_{E_2}(F) \}.
\end{equation}
The arrow in \eqref{eq:r-functor} is obtained by composing
\[
hom(E_1,F) \otimes T_{E_2}(E_1) \xrightarrow{(T_{E_2})_* \otimes \id}
hom(T_{E_2}(E_1),T_{E_2}(F)) \otimes T_{E_2}(E_1)
\]
with the evaluation map $ev: hom(T_{E_2}(E_1),T_{E_2}(F)) \otimes T_{E_2}(E_1) \longrightarrow
T_{E_2}(F)$. This means that one has a natural map from $R_{E_1,E_2}(F)$ to $T_{T_{E_2}(E_1)}
T_{E_2}(F)$, given by $(T_{E_2})_* \otimes \id$ on the first component and by the identity on
the second one. Since $(T_{E_2})_*$ is a quasi-isomorphism by Proposition
\ref{th:invertibility}, this natural transformation is an isomorphism. \qed

\begin{prop} \label{th:braid-one}
Let $E_1,E_2$ be as before, and assume in addition that $\Hom^i_{\mathfrak K}(E_2,E_1) = 0$
for all $i$. Then $T_{E_1}T_{E_2} \iso T_{E_2}T_{E_1}$.
\end{prop}

\proof The assumption implies that $T_{E_2}(E_1)$ is isomorphic to $E_1$. Hence the result
follows directly from Lemma \ref{th:braid-zero} and Proposition \ref{th:iso-implies-eq} (one
can also prove this by a direct computation, without using Lemma \ref{th:braid-zero}). \qed

\begin{prop} \label{th:braid-two}
Let $E_1,E_2 \in {\mathfrak K}$ be two $n$-spherical objects for some $n>0$. Assume that the
total dimension of $\Hom^*_{\mathfrak K}(E_2,E_1)$ is one. Then $T_{E_1}T_{E_2}T_{E_1} \iso
T_{E_2}T_{E_1}T_{E_2}$.
\end{prop}

\proof Since the twists are not affected by shifting, we may assume that $\Hom^i_{\mathfrak
K}(E_2,E_1)$ is one-dimensional for $i = 0$ and zero in all other dimensions. A simple
computation shows that
\[
T_{E_2}(E_1) \iso \{E_2 \stackrel{g}{\rightarrow} E_1\}, \quad T_{E_1}'(E_2) \iso \{E_2
\stackrel{h}{\rightarrow} E_1\}
\]
where $g$ and $h$ are nonzero maps. As $\Hom_{\mathfrak K}(E_2,E_1)$ is one-dimensional it
follows that $T_{E_2}(E_1)$ and $T_{E_1}'(E_2)$ are isomorphic up to the shift $[1]$. By
applying Lemma \ref{th:braid-zero} and Proposition \ref{th:iso-implies-eq} one finds that
\[
T_{E_1}T_{E_2}T_{E_1} \iso T_{E_1} T_{T_{E_2}(E_1)} T_{E_2} \iso T_{E_1} T_{T_{E_1}'(E_2)}
T_{E_2}.
\]
On the other hand, applying Lemma \ref{th:braid-zero} to $T_{E_1}'(E_2)$ and $E_1$, and using
Proposition \ref{th:invertibility}, shows that $T_{E_1} T_{T_{E_1}'(E_2)} T_{E_2} \iso T_{E_2}
T_{E_1} T_{E_2}$. \qed

We will now carry over the results obtained so far to the derived category $D^b(\ab')$. During
the rest of this section, $\Hom$ always means $\Hom_{D^b(\ab')}$.

\begin{definition} \label{def:d-spherical}
An object $E \in D^b(\ab')$ is called $n$-spherical for some $n>0$ if it has the following
properties:
\begin{mylist}
\item[(S1)]
$E$ has a finite resolution by injective objects in $\ab$;
\item[(S2)]
$\Hom^*(E,F)$, $\Hom^*(F,E)$ are finite-dimensional for any $F \in D^b(\ab')$.
\item[(S3)]
$\Hom^i(E,E)$ is equal to $k$ for $i = 0, n$ and zero in all other dimensions;
\item[(S4)]
The composition map $\Hom^i(F,E) \times \Hom^{n-i}(E,F) \longrightarrow \Hom^n(E,E) \iso k$ is
a nondegenerate pairing for all $F \in {\mathfrak K}$ and $i \in \Z$.
\end{mylist}
\end{definition}

Clearly, if $E$ is such an object, any finite resolution by $\ab$-injectives is an
$n$-spherical object of ${\mathfrak K}$ in the sense of Definition \ref{def:k-spherical}.
Using such a resolution, and the equivalence of categories from Proposition \ref{th:d-and-k},
one can associate to $E$ a twist functor $T_E$ which, by Proposition \ref{th:invertibility},
is an exact self-equivalence of $D^b(\ab')$. This will be independent of the choice of
resolution up to isomorphism, thanks to Proposition \ref{th:iso-implies-eq}.

\begin{lemma} \label{th:spherical-check}
In the presence of (S2) and (S3), condition (S4) is equivalent to the following apparently
weaker one:
\begin{mylist}
\item[(S4')]
There is an isomorphism $\Hom(E,F) \iso \Hom^n(F,E)^\vee$ which is natural in $F \in
D^b(\ab')$.
\end{mylist}
\end{lemma}

\proof The proof is by a `general nonsense' argument. Take any natural isomorphism as in (S4')
and let $q_F: \Hom(E,F) \times \Hom^n(F,E) \longrightarrow k$ be the family of nondegenerate
pairings induced by it. Because of the naturality, these pairings satisfy $q_F(\phi,\psi) =
q_F(\phi \circ \id_E,\psi) = q_E(\id_E, \phi \circ \psi)$. Since the pairings are all
nondegenerate, $q_E(\id_E,-): \Hom^n(E,E) \longrightarrow k$ is nonzero, hence by (S3) an
isomorphism. We have therefore shown that
\[
\Hom(E,F) \times \Hom^n(F,E) \xrightarrow{\text{composition}} \Hom^n(E,E) \iso k
\]
is a nondegenerate pairing for any $F$, which is the special case $i = 0$ of (S4). The other
cases follow by replacing $F$ by $F[i]$. \qed

\begin{lemma} \label{th:finite-resolution-check}
Let $X$ be a noetherian scheme over $k$ and $\ab = \Qco(X)$, $\ab' = \Coh(X)$. Then condition
(S4) or (S4') for an object of $D^b(\ab')$ implies condition (S1).
\end{lemma}

\proof Let ${\mathcal E}$ be an object of $D^b(\ab')$ and ${\mathcal F} \in \ab'$ a coherent
sheaf. Since ${\mathcal E}$ is bounded, and ${\mathcal F}$ has a bounded below resolution by
$\ab$-injectives, one has $\Hom^i(\mathcal E,\mathcal F) = 0$ for $i \ll 0$. Using (S4) or
(S4') it follows that $\Hom^i(\mathcal F,\mathcal E) = 0$ for $i \gg 0$, and \cite[Proposition
II.7.20]{hartshorne66} completes the proof. \qed

Now define an $(A_m)$-configuration $(m>0)$ of $n$-spherical objects in $D^b(\ab')$ to be a
collection $(E_1,\dots,E_m)$ of such objects, satisfying
\begin{equation} \label{eq:chain}
\dim_k \Hom^*_{D^b(\ab')}(E_i,E_j) = \begin{cases} 1 & |i-j| = 1, \\ 0 & |i-j| \geq 2.
\end{cases}
\end{equation}

\begin{thm} \label{th:braid-group-action}
Let $(E_1,\dots,E_m)$ be an $(A_m)$-configuration of $n$-spherical objects in $D^b(\ab')$.
Then the twists $T_{E_1},\dots,T_{E_m}$ satisfy the relations of the braid group $B_{m+1}$ up
to graded natural isomorphism. That is to say, they generate a homomorphism $\rho: B_{m+1}
\longrightarrow \Auteq(D^b(\ab'))$.
\end{thm}

This follows immediately from the corresponding results for ${\mathfrak K}$ (Propositions
\ref{th:braid-one} and \ref{th:braid-two}). One minor point remains to be cleared up: the
Theorem states that the braid relations hold up to graded natural isomorphism, whereas before
we have only talked about ordinary natural isomorphism. But one can easily see all the natural
isomorphisms which we have constructed are graded ones, essentially because everything
commutes with the translation functors. We can now state the main result of this paper:

\begin{theorem} \label{th:main}
Suppose that $n \geq 2$. Then the homomorphism $\rho$ defined in Theorem
\ref{th:braid-group-action} is injective, and in fact the following stronger
statement holds: if $g \in B_{m+1}$ is not the identity element, then
$\rho(g)(E_i) \not\iso E_i$ for some $i \in \{1,\dots,m\}$.
\end{theorem}

\section{Applications}

\subsection{Smooth projective varieties\label{subsec:sheaves}}
We now return to the concrete situation of derived categories of coherent sheaves. The main
theme will be the use of suitable duality theorems to simplify condition (S4') in the
definition of spherical objects. Throughout, all varieties will be over an algebraically
closed field $k$.

For the moment we consider only smooth projective varieties $X$, of dimension $n$. Let us
recall some facts about duality on such varieties. Serre duality says that for any $\mathcal G
\in D^b(X)$ the composition
\begin{equation} \label{eq:serre}
 \Hom^{n-*}(\mathcal G,\o_X) \otimes \Hom^*(\O,\mathcal G)
 \longrightarrow \Hom^n(\O,\o_X) = H^n(\o_X) \iso k
\end{equation}
is a nondegenerate pairing (the classical form is for a single sheaf ${\mathcal G}$; the
general case can be derived from this by induction on the length, using the Five-Lemma). Now
let $\mathcal E$ be a bounded complex of locally free coherent sheaves on $X$. For all
$\mathcal G_1, \mathcal G_2 \in D^+(X)$ there is a natural isomorphism
\begin{equation} \label{eq:locally-free}
\Hom^*({\mathcal G}_1 \otimes {\mathcal E},{\mathcal G_2}) \iso \Hom^*({\mathcal
G}_1,{\mathcal G}_2 \otimes {\mathcal E}^\vee).
\end{equation}
This is proved using a resolution ${\mathcal G}_2'$ of ${\mathcal G}_2$ by injective
quasi-coherent sheaves; the point is that ${\mathcal G}_2' \otimes {\mathcal E}^\vee$ is an
injective resolution of ${\mathcal G}_2 \otimes {\mathcal E}^\vee$ \cite[Proposition
7.17]{hartshorne66}. Setting $\mathcal G = \mathcal F \otimes \mathcal E^\vee$ in
\eqref{eq:serre} for some $\mathcal F \in D^b(X)$ and using \eqref{eq:locally-free} shows that
there is an isomorphism, natural in $\mathcal F$,
\begin{equation} \label{eq:serre-two}
\Hom^*(\mathcal E, \mathcal F) \iso \Hom^{n-*}(\mathcal F, \mathcal E \otimes \o_X)^\vee.
\end{equation}
Again by \eqref{eq:locally-free} and the standard finiteness theorems, $\Hom^*(\mathcal E,
\mathcal F) \iso {\mathbb H}^{\,*}(\mathcal E^\vee \otimes \mathcal F)$ is of finite total
dimension; hence so is $\Hom^*(\mathcal F,\mathcal E)$ by \eqref{eq:serre-two}. Finally,
because of the existence of finite locally free resolutions, everything we have said holds for
an arbitrary ${\mathcal E} \in D^b(X)$.
\begin{lemma} \label{th:nonsingular-projective}
An object ${\mathcal E} \in D^b(X)$ is spherical, in the sense of Definition
\ref{def:d-spherical}, iff it satisfies the following two conditions: $\Hom^j(\mathcal
E,\mathcal E)$ is one-dimensional for $j = 0, n$ and zero for all other $j$; and $\mathcal E
\otimes \o_X \iso \mathcal{E}$.
\end{lemma}
\proof It follows from \eqref{eq:serre-two} and the previous discussion that the conditions
are sufficient. Conversely, assume that ${\mathcal E}$ is a spherical object. Then property
(S4) and \eqref{eq:serre-two} imply that the functors $\Hom(-,{\mathcal E} \otimes \o_X)$ and
$\Hom(-,\mathcal E)$ are isomorphic. By a general nonsense argument $\mathcal E$ must be
isomorphic to $\mathcal E \otimes \o_X$. \qed

This shows that the abstract definition of spherical objects specializes to the one in section
\ref{subsec:introduction}. We will now prove the corresponding statement for twist functors.
\begin{lemma} \label{th:equivalent}
Let $\mathcal E \in D^b(X)$ be a bounded complex of locally free sheaves, which is a spherical
object. Then the twist functor $T_{\mathcal E}$ as defined in section \ref{subsec:braiding} is
isomorphic to the FMT by ${\mathcal P} = \Cone(\eta: \mathcal E^\vee \boxtimes \mathcal E
\longrightarrow \O_\Delta)$.
\end{lemma}
\proof Let ${\mathcal E'} \in {\mathfrak K}$ be a bounded resolution of $\mathcal E$ by
injective quasi-coherent sheaves. Let $T: {\mathfrak K} \longrightarrow D^+(X)$ be the functor
which sends ${\mathcal F}$ to $\Cone(ev: hom({\mathcal E},{\mathcal F}) \otimes {\mathcal E}
\longrightarrow {\mathcal F})$. We will show that the diagram
\begin{equation} \label{eq:diagram-of-twists}
\xymatrix{
 {\mathfrak K} \ar[d] \ar[dr]^{T} \ar[r]^-{T_{\mathcal E'}} &
 {\mathfrak K} \ar[d] \\
 {D^b(X)} \ar[r]^-{\Phi_{\mathcal P}} &
 {D^+(X),}
}
\end{equation}
where the unlabeled arrows are the equivalence ${\mathfrak K} \iso D^b(X)$ and its inclusion
into $D^+(X)$, commutes up to isomorphism. Since $T_{\mathcal E}$ is defined using the twist
functor $T_{\mathcal E'}$ on ${\mathfrak K}$ and ${\mathfrak K} \iso D^b(X)$, the
commutativity of \eqref{eq:diagram-of-twists} implies that $\Phi_{\mathcal P} \iso T_{\mathcal
E}$.
Take an object ${\mathcal F} \in D^b(X)$ and a resolution ${\mathcal F'} \in {\mathfrak K}$.
Then
\begin{align*}
 \Phi_{\mathcal P}({\mathcal F})
 & = \RR\pi_{2\,*}\left\{\pi_1^*\mathcal F \otimes \pi_1^*\mathcal E^\vee \otimes
 \pi_2^*\mathcal E \longrightarrow \O_\Delta \otimes \pi_1^*\mathcal F \right\} \\
 & \iso \RR\pi_{2\,*}\left\{\pi_1^*\mathcal F' \otimes \pi_1^*\mathcal E^\vee \otimes
 \pi_2^*\mathcal E \longrightarrow \O_\Delta \otimes \pi_1^*\mathcal F' \right\} \\
 & \iso \pi_{2,*}\left\{\pi_1^*\SHom(\mathcal E,\mathcal F') \otimes \pi_2^*\mathcal E
 \longrightarrow \O_\Delta \otimes \pi_1^*\mathcal F' \right\} \\
 & \iso \left\{hom(\mathcal E,\mathcal F') \otimes \mathcal E \longrightarrow \mathcal F'
 \right\} = T({\mathcal F'}),
\end{align*}
where the arrow in the last line is evaluation. This provides a natural isomorphism which
makes the left lower triangle in \eqref{eq:diagram-of-twists} commute. To deal with the other
triangle, set up a diagram as in the proof of Proposition
\ref{th:iso-implies-eq}. \qed
\begin{example}
Let $X$ be a variety which is Calabi-Yau in the strict sense, that is to say $\o_X \iso \O$
and $H^i(X,\O) = 0$ for $0 < i < n$. Then any invertible sheaf on $X$ is spherical. For the
trivial sheaf, the twist $T_{\O}$ is the FMT given by the object on $X\times X$
which is the ideal sheaf of the diagonal shifted by $[1]$. This is what Mukai
\cite{mukai87} calls the `reflection functor'.
\end{example}
\begin{lemma} \label{th:subvarieties}
Let $Y \subset X$ be a connected subscheme which is a local complete intersection, with
(locally free) normal sheaf $\nu = ({\mathcal J_Y}/{\mathcal J_Y^2})^\vee$. Assume that $\o_X|Y$
is trivial, and that $H^i(Y,\Lambda^j\nu) = 0$ for all $0 < i+j < n$. Then $\O_Y \in D^b(X)$
is a spherical object.
\end{lemma}
\proof Denote by $\iota$ the embedding of $Y$ into $X$. The local Koszul resolution of
$\iota_*\O_Y$ gives the well-known formula for the sheaf Exts, $\SExt^j(\iota_*\O_Y,
\iota_*\O_Y) \iso \iota_*(\Lambda^j\nu)$. The assumptions and the spectral sequence
$H^i(\SExt^j) \Rightarrow \Ext^{i+j}$ (i.e.\ the hypercohomology spectral sequence of $\mathbb
H\,(\RR\SHom) = \Ext$) give $\Ext^r(\iota_*\O_Y,\iota_*\O_Y) = 0$ for $0 < r < n$. We have
$\Hom(\iota_*\O_Y, \iota_*\O_Y) \iso k$, hence $\Ext^n(\iota_*\O_Y,\iota_*\O_Y) \iso k$ by
duality. \qed
\begin{example} \label{ex:-2spheres}
Let $X$ be a surface. Then any smooth rational curve $C \subset X$ with $C \cdot C = -2$
satisfies the conditions of Lemma \ref{th:subvarieties}. Now take a chain $C_1,\dots,C_m$ of
such curves such that $C_i \cap C_j = \emptyset$ for $|i-j| \geq 2$, and $C_i \cdot C_{i+1} =
1$ for $i = 1,\dots,m-1$. Then $(\O_{C_1},\dots,\O_{C_m})$ is an $(A_m)$-configuration of
spherical objects.
\end{example}
\begin{remark}
As far as Lemma \ref{th:nonsingular-projective} is concerned, one could remove the assumption
of smoothness and work with arbitrary projective varieties $X$. Serre duality must then be
replaced by the general duality theorem \cite[Theorem III.11.1]{hartshorne66} applied to the
projection $\pi: X \rightarrow \mathrm{Spec}\, k$. This yields a natural isomorphism, for
$\mathcal G \in D^-(X)$,
\[
\Ext^{n-*}(\mathcal{G},\o_X) \iso \Ext^*(\O_X,\mathcal{G})^\vee,
\]
where now $\o_X = \pi^!(\O_{\mathrm{Spec}\,k}) \in D^+(X)$ is the dualizing complex. With this
replacing \eqref{eq:serre} one can essentially repeat the same discussion as in the smooth
case, leading to an analogue of Lemma \ref{th:nonsingular-projective}. The only difference is
that the condition that ${\mathcal E}$ has a finite locally free resolution must be included
as an assumption. We do not pursue this further, for lack of a really relevant application.
\end{remark}

\subsection{Two generalisations\label{subsec:generalisations}}
We will now look at smooth quasi-projective varieties. Rather than aiming at a comprehensive
characterisation of spherical objects, we will just carry over Lemma \ref{th:subvarieties}
which provides one important source of examples.

Let $X$ be a smooth quasi-projective variety of dimension $n$, and $Y
\subset X$ a complete subscheme, of codimension $c$. $\iota$ denotes the
embedding $Y \hookrightarrow X$. Complete $X$ to a projective variety
$\bar X$. Then $Y\subset\bar X$ is closed, and $X$ is smooth, so
$\iota_*\O_Y$ has a finite locally free resolution; thus we may use
Serre duality \cite[Theorem III.11.1]{hartshorne66} on $\bar X$, and
the methods of (\ref{eq:serre-two}), to conclude that
\[
 \Hom(\iota_*\O_Y,\mathcal{F}) \iso
 \Hom^n({\mathcal F},\iota_*\O_Y \otimes \o_X)^\vee,
\]
on $X$. By continuing as in the projective case, and using the same
spectral sequence as in Lemma
\ref{th:subvarieties}, one obtains the following result\footnote{We
thank one of the referees for simplifying our original version of
the above proof.}:
\begin{lemma} \label{th:quasi-projective}
Assume that $H^i(Y,\Lambda^j\nu) = 0$ for all $0 < i+j < n$, and that $\iota^*\o_X$ is
trivial. Then $\iota_*\O_Y$ is a spherical object in $D^b(X)$. \qed
\end{lemma}
One can now e.g.\ extend Example \ref{ex:-2spheres} to quasi-projective surfaces. For
subschemes of codimension one, we will later on provide a stronger result, Proposition
\ref{th:pushforward}, which can be used to construct more interesting spherical objects.

The other generalisation which we want to look at is technically much simpler. Let $X$ be a
smooth $n$-dimensional projective variety over $k$ with an action of a finite group $G$. We
will assume that $\mathrm{char}(k) = 0$; this implies the complete reducibility of
$G$-representations, which will be used in an essential way. Let $Qco_G(X)$ be the category
whose objects are $G$-equivariant quasi-coherent sheaves, and whose morphisms are the
$G$-equivariant sheaf homomorphisms. One can write
\begin{equation} \label{eq:g-invariants}
\Hom_{Qco_G(X)}(\mathcal E_1,\mathcal E_2) = \Hom_{Qco(X)}(\mathcal E_1,\mathcal E_2)^G
\end{equation}
with respect to the obvious $G$-action on $\Hom_{Qco(X)}(\mathcal E_1,\mathcal E_2)$. Because
taking the invariant part of a $G$-vector space is an exact functor, it follows that a
$G$-sheaf is injective in $Qco_G(X)$ iff it is injective in $Qco(X)$. This can be used to show
that $Qco_G(X)$ has enough injectives, and also that $\ab = Qco_G(X)$ and its Serre
subcategory $\ab' = Coh_G(X)$ of coherent $G$-sheaves satisfy the conditions (C1)--(C4) from
section \ref{subsec:homological-algebra}. As a further application one derives a formula
similar to \eqref{eq:g-invariants} for the derived category:
\begin{equation} \label{eq:derived-invariants}
\Hom_{D^+(Qco_G(X))}(\mathcal F_1,\mathcal F_2) = \Hom_{D^+(Qco(X))}(\mathcal F_1,\mathcal
F_2)^G
\end{equation}
for all $\mathcal F_1,\mathcal F_2 \in D^+(Qco_G(X))$. This allows one to carry over the usual
finiteness results for coherent sheaf cohomology, as well as Serre duality, to the equivariant
context. The same argument as in the non-equivariant case now leads to
\begin{lemma} \label{th:equivariant}
An object $\mathcal E$ in the derived category $D^b_G(X) = D^b(Coh_G(X))$ of coherent
equivariant sheaves is spherical iff the following two conditions are satisfied:
$\Hom^j_{D^b_G(X)}(\mathcal E,\mathcal E)$ is one-dimensional for $j = 0,n$ and zero in other
degrees; and $\mathcal E \otimes \o_X$ is equivariantly isomorphic to
$\mathcal E$. \qed
\end{lemma}
Finally, one can combine the two generalisations and obtain an equivariant version of Lemma
\ref{th:quasi-projective}. This is useful in examples which arise in connection with the McKay
correspondence. We will concentrate on the simplest of these examples, which also happens to
be particularly relevant for our purpose.

Consider the diagonal subgroup $G \iso \Z/(m+1)$ of $SL_2(k)$. Write $R$ for its regular
representation and $V_1,\dots,V_m$ for its (nontrivial) irreducible representations. Let $X$
be a smooth quasiprojective surface with a complex symplectic form, carrying an effective
symplectic action of $G$. Choose a fixed point $x \in X$; the tangent space $T_xX$ must
necessarily be isomorphic to $R$ as a $G$-vector space. For $i = 1,\dots,m$ set ${\mathcal
E_i} = \O_x \otimes V_i \in Coh_G(X)$. The Koszul resolution of $\O_x$ together with
\eqref{eq:derived-invariants} shows that
\[
\Hom_{D^b_G(X)}^r(\mathcal E_i,\mathcal E_j) \iso (\Lambda^r R \otimes V_i^\vee \otimes
V_j)^G.
\]
This implies that each ${\mathcal E_i}$ is a spherical object, and that these objects form an
$(A_m)$-configuration, so that we obtain a braid group action on $D^b_G(X)$.
\begin{example}
In particular, we have a braid group action on the equivariant derived category of coherent
sheaves over ${\mathbb A}^2$, with respect to the obvious linear action of $G$ (this is
probably the simplest example of a braid group action on a category in the present paper).
\end{example}
Let $\pi: Z \rightarrow X/G$ be the minimal resolution. This is again a quasiprojective
surface with a symplectic form; it can be constructed as Hilbert scheme of $G$-clusters on
$X$. The irreducible components of $\pi^{-1}(x)$ are smooth rational curves $C_1,\dots,C_m$
which are arranged as in Example \ref{ex:-2spheres}, so that their structure sheaves generate
a braid group action on $D^b(Z)$. A theorem of Kapranov and Vasserot
\cite{kapranov-vasserot98} provides an equivalence of categories
\begin{equation} \label{eq:kv}
D^b_G(X) \iso D^b(Z),
\end{equation}
which takes ${\mathcal E}_j$ to $\O_{C_j}$ up to tensoring by a line bundle \cite[p.\
7]{kapranov-vasserot98}. This means that the braid group actions on the two categories
essentially correspond to each other. Adding the trivial one-dimensional representation $V_0$,
and the corresponding equivariant sheaf ${\mathcal E}_0 = \O_x = \O_x \otimes V_0$, extends
the action on $D^b_G(X)$ to an action of the affine braid group, except for $m = 1$.
Interestingly, the cyclic symmetry between $V_0,V_1,\dots,V_m$ is not immediately visible on
$D^b(Z)$; the equivalence \eqref{eq:kv} takes ${\mathcal E_0}$ to the structure sheaf of the
whole exceptional divisor $\pi^{-1}(x)$. Finally, everything we have said carries over to the
other finite subgroups of $SL(2,k)$ with the obvious modifications: the Dynkin diagram of type
$(A_m)$ which occurs implicitly several times in our discussion must be replaced by those of
type D/E, and one obtains actions of the corresponding (affine) generalized braid groups.

A recent deep theorem of Bridgeland, King and Reid \cite{bridgeland-king-reid99} extends the
equivalence \eqref{eq:kv} to certain higher-dimensional quotient singularities. We consider
only one very concrete case.
\begin{example}
Let $X$ be the Fermat quintic in ${\mathbb P}^4$ with the diagonal action of $G = (\Z/5)^3$
familiar from mirror symmetry. The fixed point set $X^H$ of the subgroup $H = (\Z/5)^2 \times
1$ consists of a single $G$-orbit $\Sigma$, whose structure sheaf is a spherical object in
$D^b_G(X)$. By considering other subgroups of the same kind one finds a total of ten spherical
objects, with no Homs between any two of them. Now let $\pi: Z \rightarrow X/G$ be the
crepant resolution given by the Hilbert scheme of $G$-clusters. Then $D^b_G(X) \iso D^b(Z)$ by
\cite{bridgeland-king-reid99} so that one gets corresponding spherical objects on $Z$. Because
of the nature of the equivalence, the object corresponding to $\O_\Sigma$ must be supported on the exceptional set
$p^{-1}(\Sigma)$ of the resolution. We have not determined its precise
nature, but this is clearly related to
Proposition \ref{th:pushforward} and Examples \ref{ex:threefolds} of the
next sections.
\end{example}

\subsection{Spherical and exceptional objects\label{subsec:exceptional}}
The reader familiar with the theory of exceptional sheaves \cite{rudakov90}, or with certain
aspects of tilting theory in representation theory, will have noticed a similarity between our
twist functors and mutations of exceptional objects. (See also
\cite{bondal-polishchuk93}, and note their `elliptical exceptional'
objects are
examples of $1$-spherical objects.) The braid group also occurs in the
mutation context, but there it acts on collections of exceptional objects in a triangulated
category instead of on the category itself. The relation of the two kinds of braid group
actions is not at all clear. We will here content ourselves with two observations, the first
of which is motivated by examples in \cite{kuleshov90}.
\begin{defn} \label{def:simple}
Let $X,Y$ be smooth projective varieties, with $\o_X$ trivial. A morphism $f: X
\longrightarrow Y$ (of codimension $c = \dim X - \dim Y$) is called \emph{simple} if there is
an exact triangle
\[
\O_Y \rightarrow \RR f_* \O_X \rightarrow \o_Y[-c].
\]
\end{defn}
In most applications $Y$ would be Fano, because one could then use the wealth of known results
about exceptional sheaves on such varieties. However, the general theory does not require this
assumption on $Y$.
\begin{lemma}
Suppose that $c>0$ and
\[
R^i f_* \O_X \iso
\begin{cases}
 \O_Y & \text{\em for } i = 0, \\
 0    & \text{\em for } 0 < i < c, \\
 \o_Y & \text{\em for } i = c. \\
\end{cases}
\]
Then $f$ is simple.
\end{lemma}
\proof  $\RR f_* \O_X$ is a complex of sheaves whose cohomology is nonzero only in two
degrees; a general argument, valid in any derived category, shows that there is an exact
triangle $R^0 f_* \O_X \rightarrow \RR f_* \O_X \rightarrow (R^c f_* \O_X)[-c]$. \qed
\begin{prop}
Suppose that $f$ is simple, and $\mathcal F \in D^b(Y)$ is an exceptional object, in the sense
that $\Hom(\mathcal F,\mathcal F) \iso k$ and $\Hom^i(\mathcal F,\mathcal F) = 0$ for all $i
\neq 0$. Then $\mathbf L f^*\mathcal F \in D^b(X)$ is a spherical object.
\end{prop}
\proof One can easily show, using e.g.\ a finite locally free resolution of $\mathcal F$ and a
finite injective quasi-coherent resolution of $\O_X$, that $\mathbf R f_* \mathbf L f^*
\mathcal F \iso \mathcal F \otimes^{\mathbf L} (\mathbf R f_* \O_X)$. Hence, by tensoring the
triangle in Definition \ref{def:simple} with $\mathcal F$, one obtains another exact triangle
$\mathcal F \rightarrow \mathbf R f_* \mathbf L f^* \mathcal F \rightarrow \mathcal F \otimes
\o_Y[-c]$. This yields a long exact sequence
\[
\dots \Hom^*(\mathcal F,\mathcal F) \longrightarrow \Hom^*(\mathcal F, \mathbf R f_* \mathbf L
f^* \mathcal F) \longrightarrow \Hom^{*-c}(\mathcal F,\mathcal F \otimes \o_Y) \dots
\]
The second and third group are $\Hom^*(\mathbf L f^* \mathcal F, \mathbf L f^* \mathcal F)$
and $\Hom^{\dim X-*}(\mathcal F, \mathcal F)$ by, respectively, adjointness and Serre duality.
From the assumption that ${\mathcal F}$ is exceptional, one now immediately obtains the
desired result. \qed
\begin{examples}
\begin{theoremlist}
\item
(This assumes $\text{char}(k) = 0$.) Consider a Calabi-Yau $X$ with a fibration $f: X
\longrightarrow Y$ over a variety $Y$ such that the generic fibres are elliptic curves or $K3$
surfaces. Clearly $f_*\O_X \iso \O_Y$; relative Serre duality shows that $R^cf_*\O_X \iso
\o_Y$; and in the $K3$ fibred case one has also $R^1f_*\O_X = 0$. Hence $f$ is simple.
\item
(This assumes $\text{char}(k) \neq 2$.) Let $f: X \longrightarrow Y$ be a twofold covering
branched over a double anticanonical divisor. One can use the $\Z/2$-action on $X$ to split
$f_*\O_X$ into two direct summands, which are isomorphic to $\O_Y$ and $\o_Y$ respectively;
this implies that $f$ is simple. An example, already considered in \cite{kuleshov90}, is a
$K3$ double covering of $\mathbb P^2$ branched over a sextic. Another example, which is
slightly degenerate but still works, is the unbranched covering map from a $K3$ surface to an
Enriques surface.
\item
Examples with $c=-1$ come from taking $X$ to be a smooth anticanonical divisor in $Y$, and $f$
the embedding. Then $\RR f_*\O_X = f_*\O_X \iso \{\o_Y \rightarrow \O_Y\}$ with the map given
by the section of $\o_Y^{-1}$ defining $X$. Quartic surfaces in $\mathbb P^3$ are an example
considered in \cite{kuleshov90}.
\end{theoremlist}
\end{examples}
We will now describe a second connection between spherical and exceptional objects, this time
using pushforwards instead of pullbacks. The result applies to quasi-projective varieties as
well, but it is limited to embeddings of divisors. Let $X\subset\mathbb P^N$ be a smooth quasi-projective variety
and $\iota: Y \hookrightarrow X$ an embedding of a complete connected hypersurface
$Y$. As in the parallel argument in the previous section,
we work on the projective completion $\bar X$ of $X$, in which $Y$ is
closed. By the smoothness of $X$, given $\mathcal F\in D^b(Y)$,
$\iota_*\mathcal F$ has a finite locally free resolution, and
Serre duality on $\bar X$ \cite[Theorem III.11.1]{hartshorne66} yields
\[
\Hom(\iota_*\mathcal F,\mathcal G) \iso
\Hom^{\dim X}(\mathcal G, \iota_* \mathcal F \otimes \o_X)^\vee,
\]
on $X$.
\begin{prop} \label{th:pushforward}
Assume that $\iota^*\o_X$ is trivial. If $\mathcal F \in
D^b(Y)$ is an exceptional object \emph{with a finite locally free
resolution}, then $\iota_*\mathcal F$ is spherical in $D^b(X)$.
\end{prop}
\proof In view of the previous discussion, what remains to be done is
to compute $\Hom^i(\iota_* \mathcal F, \iota_* \mathcal F)$, which by
\cite[Theorem III.11.1]{hartshorne66} applied to $\iota_*$, is isomorphic
to $\Hom^{i-1}(\mathcal F,{\mathbf L}\iota^*\iota_* \mathcal F \otimes
\o_Y)$. We will need the following result (which, perhaps
surprisingly, need not be true without the $\iota_*$\,s).
\begin{lemma} \label{th:splitting}
$\iota_*\mathbf L\iota^*(\iota_*\mathcal F)\,\iso\,\iota_*(\mathcal F\otimes\o_Y^{-1})[1]
\oplus \iota_*\mathcal F$.
\end{lemma}
\proof Replacing $\iota_*\mathcal F$ by a quasi-isomorphic complex
$\mathcal F'$ of locally free sheaves, the left hand side of the above equation is
$\iota_*(\mathcal F'|_Y)=\mathcal F' \otimes \O_Y$ which is quasi-isomorphic to
\[
 \mathcal F' \otimes \{\O(-Y) \rightarrow \O\} \htp
 \iota_* \mathcal F \otimes\{\O(-Y)\to\O\},
\]
where the arrow is multiplication by the canonical section of $\O(Y)$. Since this vanishes on
$Y$, which contains the support of $\iota_*\mathcal F$, we obtain $\iota_*(\mathcal F \otimes
\O(-Y)|_Y)[1]\ \oplus\ \iota_*\mathcal F$ as required. \qed

By hypothesis we may assume that $\mathcal F$ is a finite complex of
locally frees on $Y$, so that $\SHom(F,F)\iso\mathcal F\otimes\mathcal
F^\vee$. Thus, computing $\Hom^i(\iota_*\mathcal F,\iota_*\mathcal F)$
as the $(i-1)$th (derived/hyper) sheaf cohomology of the complex of
$\O_Y$-module sheaves ${\mathbf L}\iota^* \iota_* \mathcal F \otimes
\mathcal F^\vee \otimes \o_Y$, we may push forward to $X$ and there
use the Lemma above. That is, pick an injective resolution $\O_Y \to
I$ on $Y$, so that $\Hom^i(\iota_*\mathcal F,\iota_*\mathcal F)$ is
the $(i-1)$th cohomology of
\[
\Gamma_Y({\mathbf L}\iota^*\iota_*\mathcal F\otimes\mathcal F^\vee
\otimes\o_Y\otimes I),
\]
where $\Gamma$ is the global section functor. Pushing forward to $X$, this is
\[
\Gamma_X(\iota_*({\mathbf L}\iota^*\iota_*\mathcal F)\otimes\iota_*((\mathcal
F)^\vee\otimes\o_Y\otimes I)),
\]
which by Lemma \ref{th:splitting} is
\[
\Gamma_X(\iota_*(\mathcal F\otimes\mathcal F^\vee\otimes I))[1]\ \oplus\
\Gamma_X(\iota_*(\mathcal F\otimes\mathcal F^\vee\otimes\o_Y\otimes I)).
\]
This may be brought back onto $Y$ to give the $(i-1)$th cohomology of $\Gamma_Y(\mathcal
F\otimes\mathcal F^\vee\otimes I)[1]\ \oplus\ \Gamma_Y(\mathcal F\otimes\mathcal
F^\vee\otimes\o_Y\otimes I)$. This is $\Hom^i(\mathcal F,\mathcal F) \oplus
\Hom^{n-i}(\mathcal F,\mathcal F)^\vee$, where for the second term we have used Serre duality
on $Y$. Since $\mathcal F$ is exceptional this completes the proof. \qed

\subsection{Elliptic curves\label{subsec:elliptic-curve}}
The homological mirror conjecture for elliptic curves has been studied extensively by
Polishchuk and Zaslow \cite{polishchuk-zaslow98} \cite{polishchuk98} (unfortunately, their
formulation of the conjecture differs somewhat from that in section \ref{subsec:mirror}, so
that their results cannot be applied directly here). Polishchuk
\cite{polishchuk96} and Orlov \cite{orlov97}, following earlier
work of Mukai \cite{mukai81}, have completely determined the automorphism group of the derived category of
coherent sheaves. These are difficult results, to which we have little to add. Still, it is maybe
instructive to see how things work out in a well-understood case.

We begin with the symplectic side of the story. Let $(M,\beta)$ be the torus $M = \R/\Z \times
\R/\Z$ with its standard volume form $\beta = ds_1 \wedge ds_2$. Matters are slightly more
complicated than in section \ref{subsec:mirror}, because the fundamental group is nontrivial.
In particular, the $\smooth$-topology on $\Symp(M,\beta)$ is no longer the correct one; this
is due to the fact that Floer cohomology is not invariant under arbitrary isotopies, but only
under Hamiltonian ones. There is a bi-invariant foliation ${\mathfrak F}$ of codimension two
on $\Symp(M,\beta)$, and the Hamiltonian isotopies are precisely those which are tangent to
the leaves. To capture this idea one introduces a new topology, the Hamiltonian topology, on
$\Symp(M,\beta)$. This is the topology generated by the leaves of ${\mathfrak F}|U$, where $U
\subset \Symp(M,\beta)$ runs over all $\smooth$-open subsets. To avoid confusion, we will
write $\Symph(M,\beta)$ whenever we have the Hamiltonian topology in mind, and call this the
Hamiltonian automorphism group; this differs from the terminology in most of the literature
where the name is reserved for what, in our terms, is the connected component of the identity
in $\Symph(M,\beta)$. The difference between the two topologies becomes clear if one considers
the group $\Aff(M) = M \rtimes SL(2,\Z)$ of oriented affine diffeomorphisms of $M$. As a
subgroup of $\Symp(M,\beta)$ this has its Lie group topology, in which the translation
subgroup $M$ is connected. In contrast, as a subgroup of $\Symph(M,\beta)$ it has the discrete
topology.
\begin{lemma} \label{th:moser}
The embedding of $\Aff(M)$ into $\Symph(M,\beta)$ as a discrete subgroup is a homotopy
equivalence.
\end{lemma}
The proof consists of combining the known topology of $\Diff^+(M)$, Moser's theorem that
$\Symp(M,\beta) \subset \Diff^+(M)$ is a homotopy equivalence, and the flux homomorphism which
describes the global structure of the foliation ${\mathfrak F}$. We omit the details.

Let $\pi: \R \longrightarrow \RP{1}$, $s \mapsto [\cos(\pi s): \sin(\pi s)]$ be the universal
covering of $\RP{1}$. Consider the subgroup $\widetilde{SL}(2,\R) \subset SL(2,\R) \times
\Diff(\R)$ of pairs $(g,\tilde{g})$ such that $\tilde{g}$ is a lift of the action of $g$ on
$\RP{1}$. $\widetilde{SL}(2,\R)$ is a central extension of $SL(2,\R)$ by $\Z$ (topologically,
it consists of two copies of the universal cover). We define a graded symplectic automorphism
of $(M,\beta)$ to be a pair
\[
(\phi,\tilde{\phi}) \in \Symph(M,\beta) \times \smooth(M,\widetilde{SL}(2,\R))
\]
such that $\tilde{g}$ is a lift of $Dg: M \longrightarrow SL(2,\R)$; here we have used the
standard trivialisation of $TM$. The graded symplectic automorphisms form a group under the
composition $(\phi,\tilde{\phi}) (\psi,\tilde{\psi}) = (\phi\psi, (\tilde{\phi} \circ \psi)
\tilde{\psi})$. We denote this group by $\Sympgrh(M,\beta)$, and equip it with the topology
induced from $\Symph(M,\beta) \times \smooth(M,\widetilde{SL}(2,\R))$. It is a central
extension of $\Symph(M,\beta)$ by $\Z$. One can easily verify that the definition is
equivalent to that in \cite{seidel99}, which in turn goes back to ideas of Kontsevich
\cite{kontsevich94}.

Even in this simplest example, the construction of the derived Fukaya category
$D^b\Fuk(M,\beta)$ has not yet been carried out in detail, so we will proceed on the basis of
guesswork in the style of section \ref{subsec:mirror}. The basic objects of $D^b\Fuk(M,\beta)$
are pairs $(L,E)$ consisting of a Lagrangian submanifold and a flat unitary bundle on it.
Thus, in addition to symplectic automorphisms, the category should admit another group of
self-equivalences, which act on all objects $(L,E)$ by tensoring $E$ with some fixed flat
unitary line bundle $\xi \longrightarrow M$. The two kinds of self-equivalence should give a
homomorphism
\begin{equation} \label{eq:sympgrh}
\gamma: G \stackrel{\mathrm{def}}{=} M^\vee \rtimes \pi_0(\Sympgrh(M,\beta)) \longrightarrow
\Auteq(D^b\Fuk(M,\beta)),
\end{equation}
where $M^\vee = H^1(M;\R/\Z)$ is the Jacobian, or dual torus. In order to make the picture
more concrete, we will now write down the group $G$ explicitly. Take the standard presentation
of $SL(2,\Z)$ by generators
 $g_1 = \left(\begin{smallmatrix}1 & 1 \\ 0 & 1 \end{smallmatrix}\right)$,
 $g_2 = \left(\begin{smallmatrix}1 & 0 \\ -1 & 1 \end{smallmatrix}\right)$
and relations $g_1g_2g_1 = g_2g_1g_2$, $(g_1g_2)^6 = 1$. Let $\widetilde{SL}(2,\Z) \subset
\widetilde{SL}(2,\R)$ be the preimage of $SL(2,\Z)$. One can lift $g_1,g_2$ to elements $a_1 =
(g_1,\tilde{g}_1)$ and $a_2 = (g_2,\tilde{g}_2)$ in $\widetilde{SL}(2,\Z)$ which satisfy
$\tilde{g}_1(1/2) = 1/4$ and $\tilde{g}_2(1/4) = 0$. Together with the central element $t =
(\id, s \mapsto s-1)$ these generate $\widetilde{SL}(2,\Z)$, and one can easily work out what
the relations are:
\[
\widetilde{SL}(2,\Z) = \leftsc a_1,a_2,t \suchthat a_1a_2a_1 = a_2a_1a_2,\; (a_1a_2)^6 =
t^2,\; [a_1,t] = [a_2,t] = 1 \rightsc.
\]
Any element of $(g,\tilde{g}) \in \widetilde{SL}(2,\Z)$ defines a graded symplectic
automorphism of $(M,\beta)$: one simply takes $\phi = g$ and $\tilde{\phi}$ to be the constant
map with value $\tilde{g}$. Moreover, any translation of $M$ has a canonical lift to a graded
symplectic automorphism, by taking $\tilde{\phi}$ to be the constant map with value $1 \in
\widetilde{SL}(2,\R)$. These two observations together give a subgroup $\widetilde{\Aff}(M) =
M \rtimes \widetilde{SL}(2,\Z)$ of $\Sympgrh(M,\beta)$, which fits into a commutative diagram
\[
\xymatrix{
 {1} \ar[r] &
 {\Z} \ar[d]_{=} \ar[r] &
 {\widetilde{\Aff}(M)} \ar@{^(->}[d] \ar[r] &
 {\Aff(M)} \ar[r] \ar@{^(->}[d] &
 {1} \\
 {1} \ar[r] &
 {\Z} \ar[r] &
 {\Sympgrh(M,\beta)} \ar[r] &
 {\Symph(M,\beta)} \ar[r] &
 {1.}
}
\]
Hence, in view of Lemma \ref{th:moser}, $\pi_0(\Sympgrh(M,\beta)) \iso \widetilde{\Aff}(M)$.
After spelling out everything one finds that $G$ is the semidirect product $(\R/\Z)^4 \rtimes
\widetilde{SL}(2,\Z)$, with respect to the action of $\widetilde{SL}(2,\Z)$ on $\R^4$ given by
\begin{equation} \label{eq:real-action}
 a_1 \mapsto
 \begin{pmatrix}
 1 & 1 & 0  & 0 \\
 0 & 1 & 0  & 0 \\
 0 & 0 & 1  & 0 \\
 0 & 0 & -1 & 1
 \end{pmatrix}, \quad
 a_2 \mapsto
 \begin{pmatrix}
 1  & 0 & 0 & 0 \\
 -1 & 1 & 0 & 0 \\
 0  & 0 & 1 & 1 \\
 0  & 0 & 0 & 1
 \end{pmatrix}, \quad
 t \mapsto \id.
\end{equation}
We now pass to the mirror dual side. Let $X$ be a smooth elliptic curve over $\C$. We choose a
point $x_0 \in X$ which will be the identity for the group law on $x$. The derived category
$D^b(X)$ has self-equivalences
\[
T_\O,\;S \quad \text{and} \quad R_x, \; L_x, \; T_{\O_x} \; (x \in X)
\]
defined as follows: $T_\O$ is the twist by $\O$, which is spherical for obvious reasons. $S$
is the original example of a Fourier-Mukai transform, $S = \Phi_{\mathcal L}$ with ${\mathcal
L} = \O(\Delta - \{x_0\} \times X - X \times \{x_0\})$ the Poincar{\'e} line bundle. It maps
the structure sheaves of points $\O_x$ to the line bundles $\O(x-x_0)$, and was shown to be an
equivalence by Mukai \cite{mukai81}. $R_x$ is the self-equivalence induced by the translation
$y \mapsto y + x$; $L_x$ is the functor of tensoring with the degree zero line bundle
$\O(x-x_0)$; and $T_{\O_x}$ is the twist along $\O_x$ which is spherical by Lemma
\ref{th:subvarieties}. These functors have the following properties:
\begin{align}
 \label{eq:commute}
 & [L_x,R_y] \iso \id \quad \text{for all $x,y$}, \\
 \label{eq:twisting-is-tensoring}
 & T_{\O_x} \text{ is isomorphic to } \O(x) \otimes -, \\
 \label{eq:mukai} & S^4 \iso [-2], \\
 \label{eq:braid-mukai}
 & T_{\O_{x_0}} T_\O T_{\O_{x_0}} \iso T_\O T_{\O_{x_0}} T_\O \iso S^{-1}, \\
 \label{eq:conjugation-one}
 & T_{\O_{x_0}} R_x T_{\O_{x_0}}^{-1} \iso R_x L_x^{-1}, \\
 \label{eq:conjugation-two}
 & T_{\O_{x_0}} L_x T_{\O_{x_0}}^{-1} \iso L_x, \\
 \label{eq:conjugation-three}
 & T_\O R_x T_\O^{-1} \iso R_x, \\
 \label{eq:conjugation-four}
 & T_\O L_x T_\O^{-1} \iso R_x L_x.
\end{align}
Most of these isomorphisms are easy to prove; those which present any difficulties are
\eqref{eq:twisting-is-tensoring}, \eqref{eq:mukai}, and \eqref{eq:braid-mukai}. The first and
third of these are proved below, and the second one is a consequence of \cite[Theorem
3.13(1)]{mukai81}.
\proof[Proof of \eqref{eq:twisting-is-tensoring}] (This argument is valid for the structure
sheaf of a point on any algebraic curve.) A simple computation shows that the dual in the
derived sense is $\O_x^\vee \iso \O_x[-1]$. The formula for inverses of FMTs, for which see
e.g.\ \cite[Lemma 4.5]{bridgeland98}, shows that $T_{\O_x}^{-1} \iso \Phi_{\mathcal Q}$ for
some object ${\mathcal Q}$ fitting into an exact triangle
\[
\mathcal Q \longrightarrow \O_\Delta \stackrel{f}{\longrightarrow} \O_{(x,x)}.
\]
When following through the computation it is not easy to keep track of the map $f$, but that
is not really necessary. All we need to know is that $f \neq 0$, which is true because the
converse would violate the fact that $\Phi_{\mathcal Q}$ is an equivalence. Then, since any
morphism $\O_\Delta \rightarrow \O_{(x,x)}$ in the derived category is represented by a
genuine map of sheaves, $f$ must be some nonzero multiple of the obvious restriction map. It
follows ${\mathcal Q}$ is isomorphic to the kernel of $f$, which is $\O_{\Delta} \otimes
\pi_1^*\O(-x)$. This means that $T_{\O_x}^{-1}$ is the functor of tensoring with $\O(-x)$.
Passing to inverses yields the desired result. \qed
\proof[Proof of \eqref{eq:braid-mukai}] The equality between the first two terms follows from
Theorem \ref{th:braid-group-action}, because $\O_{x_0},\O$ form an $(A_2)$-configuration of
spherical objects. By the standard formula for the adjoints of a FM transform, the inverse of
$S$ is the FMT with ${\mathcal L^\vee}[1]$. By definition $T_\O$ is the FMT with
$\O(-\Delta)[1]$. Using \eqref{eq:twisting-is-tensoring} it follows that $T_{\O_{x_0}} T_\O
T_{\O_{x_0}}$ is the FMT with $\pi_1^*\O(x_0) \otimes \O(-\Delta)[1] \otimes \pi_2^*\O(x_0)
\iso {\mathcal L^\vee}[1]$. \qed

Equations \eqref{eq:mukai} and \eqref{eq:braid-mukai} show that $(T_\O T_{\O_{x_0}})^6 \iso
[2]$. Therefore one can define a homomorphism
\[
\widetilde{SL}(2,\Z) \longrightarrow \Auteq(D^b(X))
\]
by mapping the generators $a_1,a_2,t$ to $T_\O, T_{\O_{x_0}}$ and the translation $[1]$; this
already occurs in Mukai's paper \cite{mukai81}, slightly disguised by the fact that he uses a
different presentation of $SL(2,\Z)$. The functors $L_x,R_x$ yield another homomorphism $X
\times X \longrightarrow \Auteq(D^b(X))$; and one can combine the two constructions into a map
\begin{equation} \label{eq:elliptic-autos}
\gamma' : G' \stackrel{\text{def}}{=} (X \times X) \rtimes \widetilde{SL}(2,\Z)
\longrightarrow \Auteq(D^b(X)).
\end{equation}
Here the semidirect product is taken with respect to the $\widetilde{SL}(2,\Z)$-action on $X
\times X$ indicated by \eqref{eq:conjugation-one}--\eqref{eq:conjugation-four}; explicitly, it
is given by the matrices
\begin{equation} \label{eq:complex-action}
 a_1 \mapsto \begin{pmatrix} 1 & 1 \\ 0 & 1 \end{pmatrix}, \quad
 a_2 \mapsto \begin{pmatrix} 1 & 0 \\ -1 & 1 \end{pmatrix}, \quad
 t \mapsto \id.
\end{equation}
\begin{lemma} The group $G$ in \eqref{eq:sympgrh} is isomorphic to the group $G'$ in
\eqref{eq:elliptic-autos}.
\end{lemma}
\proof Introduce complex coordinates $z_1 = r_1 + ir_4$, $z_2 = r_2-ir_3$ on $\R^4/\Z^4$. Then
the action of $\widetilde{SL}(2,\Z)$ described in \eqref{eq:real-action} becomes $\C$-linear,
and is given by the same matrices as in \eqref{eq:complex-action}. This is sufficient to
identify the two semidirect products which define $G$ and $G'$. We should point out that
although the argument is straightforward, the change of coordinates is by no means obvious
from the geometric point of view: a look back at the definition of $G$ shows that $z_1,z_2$
mix genuine symplectic automorphisms with the extra symmetries of $D^b\Fuk(M,\beta)$ which
come from tensoring with flat line bundles. \qed

The way in which this fits into the general philosophy is that one expects to have a
commutative diagram, with the right vertical arrow given by Kontsevich's conjecture,
\begin{equation} \label{eq:big-picture}
\xymatrix{
 {G} \ar[r]^-{\gamma} \ar[d]_{\iso} &
 {\Auteq(D^b\Fuk(M,\beta))} \ar[d]_{\iso} \\
 {G'} \ar[r]^-{\gamma'} &
 {\Auteq(D^b(X)).}
}
\end{equation}
To be accurate, one should adjust the modular parameter of $X$ and the volume of $(M,\beta)$,
eventually introducing a complex part $\beta_\C$ as in Remark \ref{remark:imaginary}, so that
they are indeed mirror dual. This has not played any role up to now, since the groups $G$ and
$G'$ are independent of the parameters, but it would become important in further study. A
theorem of Orlov \cite{orlov97} says that $\gamma'$ is always injective, and is an isomorphism
iff $X$ has no complex multiplication. Only the easy part of the theorem is important for us
here: if $X$ has complex multiplication then its symmetries induce additional automorphisms of
$D^b(X)$, which are not contained in the image of $\gamma'$. Therefore, if the picture
\eqref{eq:big-picture} is correct, the derived Fukaya category for the corresponding values of
$\beta_\C$ must admit exotic automorphisms which do not come from symplectic geometry or from
flat line bundles. It would be interesting to check this claim, especially because similar
phenomena may be expected to occur in higher dimensions.

We will now apply the intuition provided by the general discussion to the specific topic of
braid group actions. To a simple closed curve $S$ on $(M,\beta)$ one can associate a Dehn
twist $\tau_S \in \Symph(M,\beta)$ which is unique up to Hamiltonian isotopy. This is defined
by taking a symplectic embedding $\iota$ of $(U,\theta) = ([-\epsilon;\epsilon] \times \R/\Z,
ds_1 \wedge ds_2)$ into $M$ for some $\epsilon>0$, with $\iota(\{0\} \times \R/\Z) = S$, and
using a local model
\[
\tau: U \longrightarrow U, \quad \tau(s_1,s_2) = (s_1,s_2 - h(s_1))
\]
where $h \in \smooth(\R,\R)$ is some function with $h(s) = 0$ for $s \leq -\epsilon/2$, $h(s)
= 1$ for $s \geq \epsilon/2$, and $h(s) + h(-s) = 1$ for all $s$. The interesting fact is that
the Dehn twists along two parallel geodesic lines are not Hamiltonian isotopic: they differ by
a translation which depends on the area lying between the two lines. Now take
\[
S_1 = \R/\Z \times \{0\}, \quad S_2 = \{0\} \times \R/\Z, \quad S_3 = \R/\Z \times \{1/2\}.
\]
This is an $(A_3)$-configuration of circles. Hence the Dehn twists $\tau_{S_1}, \tau_{S_2},
\tau_{S_3}$ define a homomorphism from the braid group $B_4$ to $\pi_0(\Symph(M,\beta))$.
However, this is not injective: $\tau_{S_3}^{-1}\tau_{S_1}$ is Hamiltonian isotopic to a
translation which has order two, so that the nontrivial braid $(g_3^{-1}g_1)^2 \in B_4$ gets
mapped to the identity element. The natural lift of this homomorphism to $\Sympgrh(M,\beta)$
has the same non-injectivity property. Guided by mirror symmetry, one translates this example
into algebraic geometry as follows: ${\mathcal E}_1 = \O_{x_0}, {\mathcal E}_2 = \O, {\mathcal
E}_3 = \O_x \in D^b(X)$, where $x \neq x_0$ is a point of order two on $X$, form an
$(A_3)$-configuration of spherical objects. Hence their twist functors generate a weak action
of $B_4$ on $D^b(X)$. By \eqref{eq:twisting-is-tensoring} $T_{{\mathcal E}_3}^{-1}T_{{\mathcal
E}_1}$ is the functor of tensoring with $\O(x-x_0)$. Since the square of this is the identity
functor, we have the same relation as in the symplectic case, so that the action is not
faithful.

\subsection{$K3$ surfaces\label{subsec:k3}}
Let $X$ be a smooth complex $K3$ surface. Consider, as in Example \ref{ex:-2spheres}, a chain
of embeddings $\iota_1,\dots,\iota_m: {\mathbb P}^1 \longrightarrow X$ whose images $C_i$
satisfy $C_i \cdot C_j = 1$ for $|i-j| = 1$, and $C_i \cap C_j = 0$ for $|i-j| \geq 2$. One
can then use the structure sheaves $\O_{C_i}$ to define a braid group action on $D^b(X)$.
However, this is not the only way:
\begin{prop} \label{th:k3}
For each $i = 1, \dots, m$ choose ${\mathcal E}_i$ to be either $\O(-C_i)$ or else
$\O_{C_i}(-1):=(\iota_i)_*\O_{\mathbb P_1}(-1)$. Then the ${\mathcal E}_i$ form an $(A_m)$-configuration of
spherical objects in $D^b(X)$, and hence generate a weak braid group action on that category.
\qed
\end{prop}
The choice can be made for each ${\mathcal E}_i$ independently. These multiple possibilities
are relevant from the mirror symmetry point of view. This is explained in \cite{thomas99}, so
we will only summarize the discussion here.

Suppose that $X$ is elliptically fibred with a section $S$. Its mirror should
be the
symplectic four-manifold $(M,\beta)$ with $M = X$ and where $\beta$ is the real part of some
holomorphic two-form on $X$ (hyperk\"ahler rotation). The smooth holomorphic
curves in $M$
are precisely the Lagrangian submanifolds in $(M,\beta)$ which are special
(with respect to
the calibration given by the K\"ahler form of a Ricci-flat metric on $X$). In
particular,
the curves $C_i$ turn into an $(A_m)$-configuration of Lagrangian
two-spheres; hence the generalized Dehn
twists along them generate a homomorphism $B_{m+1} \rightarrow
\pi_0(\Sympgr(M,\beta))$. One can wonder what the corresponding
braid group action on $D^b(X)$ should be. This question is not
really meaningful without a {\em distinguished} equivalence between
the derived Fukaya category of $(M,\beta)$ and that of coherent
sheaves on $X$, which is not what is predicted by Kontsevich's conjecture.
But if we adopt the Strominger-Yau-Zaslow \cite{strominger-yau-zaslow96}
picture of mirror symmetry, then conjecturally there should be a distinguished
full and faithful embedding of triangulated categories
\[
D^b\Fuk(M,\beta) \hookrightarrow D^b(X)
\]
induced by the particular special Lagrangian torus fibration of $M$ that comes
from the elliptic fibration of $X$ (this fibration may, of course, not
be distinguished). That this is an embedding, and not an equivalence,
is a feature of even-dimensional mirror symmetry.
This embedding should be an extension of the Fourier-Mukai transform
which takes special
Lagrangian submanifolds of $M$ (algebraic curves in $X$) to coherent
sheaves on $X$ using the relative Poincar\'e sheaf on $X\times_{\mathbb P^1}X$
that comes from considering the elliptic fibres to be self-dual using the
section; see e.g.\ \cite{thomas99}.

Assuming this, it now makes sense to ask what spherical objects of $D^b(X)$
correspond to the special Lagrangian spheres $C_1,\dots,C_m$. The
Fourier-Mukai transform takes any special Lagrangian submanifold
$C$ which is a section of the elliptic fibration to the invertible
sheaf $\O(S-C)$; and if $C$ lies in a fibre of the fibration, it
goes to the structure sheaf $\O_C$. If we assume that all curves
$C_i$ fall into one of these two categories, and that $S$ intersects all those
of
them which lie in one fibre, then the Fourier-Mukai transform takes the special Lagrangian
submanifolds $C_1,\dots,C_m$ in $(M,\beta)$ to sheaves ${\mathcal
E_1},\dots,{\mathcal E}_m$
as in Proposition \ref{th:k3}, tensored with $\O(S)$. Then, up to
the minor difference of tensoring by $\O(S)$, one of the braid group
actions mentioned in that Proposition would be the correct mirror dual of
the symplectic one.

As mentioned in section \ref{subsec:generalisations}, such configurations of curves $C_i$ are
the exceptional loci in the resolution of any algebraic surface with an $(A_m)$-singularity.
Now, $(A_m)$-configurations of Lagrangian two-spheres occur as vanishing cycles in the
smoothing of the same singularity. Thus, in a sense, mirror symmetry interchanges smoothings
and resolutions. A more striking, though maybe less well understood, instance of this
phenomenon is Arnold's strange duality (see e.g.\ \cite{pinkham77}), which has been
interpreted as a manifestation of mirror symmetry by a number of people
(Aspinwall and Morrison, Kobayashi, Dolgachev, Ebeling, etc.). Each of the 14
singular affine surfaces $S(c_1,c_2,c_3)$ on Arnold's list has a natural compactification
$\overline{S}(c_1,c_2,c_3)$ which has four singular points. One of these points is the
original singularity at the origin; the other three are quotient singularities lying on the
divisor at infinity, which is a ${\mathbb P}^1$. One can smooth the singular point at the
origin; the intersection form of the vanishing cycles obtained in this way is $T(c_1,c_2,c_3)
\oplus H$, where $T(c_1,c_2,c_3)$ is the matrix associated to the Dynkin-type diagram

\includefigure{dynkin}{dynkin.eps}{ht}

and $H = \left(\begin{smallmatrix} 0 & 1 \\ 1 & 0 \end{smallmatrix}\right)$. On the other
hand, one can resolve the three singular points at infinity. Inside the resolution, this
yields a configuration of smooth rational curves of the form $T(b_1,b_2,b_3)$ for certain
other numbers $(b_1,b_2,b_3)$. One can also do the two things together: this removes all
singularities, yielding a smooth $K3$ surface $X(c_1,c_2,c_3)$ with a splitting of its
intersection form as
\[
T(c_1,c_2,c_3) \oplus H \oplus T(b_1,b_2,b_3).
\]
Strange duality is the observation that the numbers $(b_1,b_2,b_3)$ associated to one
singularity on the list occur as $(c_1,c_2,c_3)$ for another singularity, and vice versa.
Kobayashi \cite{kobayashi95} (extended by Ebeling \cite{ebeling96} to more general
singularities) explains this by showing that the $K3$s $X(c_1,c_2,c_3)$ and $X(b_1,b_2,b_3)$
belong to mirror dual families. The associated map on homology interchanges the
$T(c_1,c_2,c_3)$ and $T(b_1,b_2,b_3)$ summands of the intersection form (the extra hyperbolic
of vanishing cycles goes to the $H^0 \oplus H^4$ of the other $K3$ surface). Thus, the
smoothing of the singular point at the origin in $\overline{S}(c_1,c_2,c_3)$ corresponds, in a
slightly vague sense, to the resolution of the divisor at infinity of
$\overline{S}(b_1,b_2,b_3)$. From our point of view, since the rational curves at infinity in
$X(b_1,b_2,b_3)$ can be used to define a braid group action on its derived category, one would
like to have a similar configuration of Lagrangian two-spheres (vanishing cycles) in the
finite part of $X(c_1,c_2,c_3)$. On the level of homology, such a configuration exists of
course, but it is apparently unknown whether it can be realized geometrically (recall that
Lagrangian submanifolds can have many more non-removable intersection points than their
intersection number suggests).

\subsection{Singularities of threefolds\label{subsec:singularities}}
Throughout the following discussion, all varieties will be smooth projective threefolds which
are Calabi-Yau in the strict sense (some singular threefolds will also occur, but they will be
specifically designated as such). Let $X$ be such a variety.
\begin{examples}\label{ex:threefolds}
Any invertible sheaf on $X$ is a spherical object in $D^b(X)$. If $S$ is a smooth connected
surface in $X$ with $H^1(S,\O_S) = H^2(S,\O_S) = 0$ (e.g.\ a rational surface or Enriques
surface), the structure sheaf $\O_S$ is a spherical object, by Lemma \ref{th:subvarieties}.
Similarly, for $C$ a smooth rational curve in $X$ with normal bundle $\nu_C \iso \O_{{\mathbb
P}^1}(-1) \oplus \O_{{\mathbb P}^1}(-1)$ (usually referred to as a $(-1,-1)$-curve), $\O_C$ is
spherical. The ideal sheaf ${\mathcal J}_C$ of such a curve is also a spherical object; this
follows from ${\mathcal J}_C[1] \iso T_\O(\O_C)$.
\end{examples}
Supposing the ground field to be $k = \C$, we will now return to the conjectural duality
between smoothings and resolutions that already played a role in the previous section, and which has been considered by many physicists.
(Of course our interest in this is in trying to mirror Dehn twists on
smoothings, which arise as monodromy transformations around a degeneration
of the smoothing collapsing the appropriate spherical vanishing cycle, by
twists on the derived categories of the resolutions.)
To explain the approach of physicists (as described in
\cite{morrison99}, for instance), it is
better to adopt the traditional point of view in which mirror symmetry relates the combined
complex and (complexified) K{\"a}hler moduli spaces of two varieties, rather than Kontsevich's
conjecture which considers a fixed value of the moduli variables. Then the idea can be phrased
like this: moving towards the discriminant locus in the complex moduli space of a variety $X$,
which means degenerating it to a singular variety $Y$, should be mirror dual to going to a
`boundary wall' of the complexified K\"ahler cone of the mirror $\widehat X$ (the annihilator
of a face of the Mori cone), thus inducing an extremal contraction $\widehat X \rightarrow
\widehat{Y}$. A second application of the same idea, with the roles of the mirrors reversed,
shows that an arbitrary crepant resolution $Z \rightarrow Y$ should be mirror dual to a
smoothing $\widehat{Z}$ of $\widehat{Y}$. A case that is reasonably well-understood is that of
the ordinary double point (ODP: $x^2+y^2+z^2+t^2=0$ in local analytic coordinates)
singularity. ODPs should be self-dual, in the sense that if $Y$ has $d$ distinct ODPs then so
does $\widehat{Y}$ (this can be checked for Calabi-Yau hypersurfaces in toric varieties, for
instance). We now review briefly Clemens' work \cite{clemens83} on the homology of smoothings
and resolutions of such singularities.

A degeneration of $X$ to a variety $Y$ with $d$ ODPs determines $d$ vanishing cycles in $X$,
and hence a map $v: \Z^d \longrightarrow H_3(X)$. Let $v^\vee: H_3(X) \longrightarrow \Z^d$ be
the Poincar{\'e} dual of $v$. Suppose that $Y$ has a crepant resolution $Z$ which, in local
analytic coordinates near each ODP, looks like the standard small resolution. This means that
the exceptional set in $Z$ consists of $d$ disjoint $(-1,-1)$-curves. By \cite{clemens83}
\cite{gross95} one has
\[
H_3(Z) \iso \ker(v^\vee)/\im(v)
\]
and exact sequences
\begin{align*}
 & H_3(X) \stackrel{v^\vee}{\longrightarrow} \Z^d \longrightarrow H_2(Z) \longrightarrow
 H_2(X) \longrightarrow 0, \\
 & 0 \longrightarrow H_4(X) \longrightarrow H_4(Z) \longrightarrow \Z^d
 \stackrel{v}{\longrightarrow} H_3(X).
\end{align*}
Thus, if there are $r$ relations between the vanishing cycles (the image of $v$ is of rank
$d-r$), the Betti numbers are
\begin{equation} \label{eq:betti-numbers}
\begin{split}
b_2(Z) = b_2(X) + r, \quad b_3(Z) = b_3(X) - 2(d-r), \quad \\ b_4(Z) = b_4(X) + r.
\end{split}
\end{equation}
Topologically, $Z$ arises from $X$ through codimension three surgery along the vanishing
cycles, and the statements above can be proved e.g.\ by considering the standard cobordism
between them. More intuitively, one can explain matters as follows. Going from $X$ to $Y$
shrinks the vanishing cycles to points; at the same time, the relations between vanishing
cycles are given by four-dimensional chains which become cycles in the limit $Y$, because
their boundaries shrink to points. This means that we lose $d-r$ generators of $H_3$ and get
$r$ new generators of $H_4$. In $Z$, there are $d-r$ relations between the homology
classes of the exceptional ${\mathbb P}^1$s; these relations are pullbacks of closed
three-dimensional cycles on $Y$ which do not lift to cycles on $Z$, so that going from $Y$ to
$Z$ adds $r$ new generators to $H_2$ while removing another $d-r$ generators from $H_3$.
Finally $H_4(Z) = H_4(Y)$ for codimension reasons.

Mirror symmetry exchanges odd and even-dimensional homology, so if $X$ and $Z$ have mirrors
$\widehat{X}$ and $\widehat{Z}$ then
\[
\begin{split}
 b_2(\widehat{Z}) = b_2(\widehat{X}) - (d-r), \quad
 b_3(\widehat{Z}) = b_2(\widehat{X}) + 2r, \quad \\
 b_4(\widehat{Z}) = b_4(\widehat{X}) - (d-r).
\end{split}
\]
The suggested explanation, in the general framework explained above, is that $\widehat{Z}$
should contain $d$ vanishing cycles with $d-r$ relations between them, obtained from a
degeneration to a variety $\widehat{Y}$ with $d$ ODP, and that $\widehat{X}$ should be a
crepant resolution of $\widehat{Y}$. Thus, mirror symmetry exchanges ODPs with the opposite
number of relations between their vanishing cycles. Moreover, to the $d$ vanishing cycles in
the original variety $X$ correspond $d$ rational $(-1,-1)$-curves in its mirror $\widehat{X}$.
It seems plausible to think that the structure sheaves or ideal sheaves of these curves
(possibly twisted by some line bundle) should be mirror dual to the Lagrangian spheres
representing the vanishing cycles in $X$; however, as in the $K3$ case, such a statement is
not really meaningful unless one has chosen some specific equivalence $D^b(X) \iso
D^b\Fuk(\widehat{X})$.
\begin{remark} \label{remark:resolution}
When $r = 0$, $H_2(Y) \iso H_2(X)$ so that the exceptional cycles in $Y$ are
homologous to zero. This is not possible if the resolution is algebraic, so
we exclude this case, and also the case $d = r$ to avoid the same problem
on the mirror.
\end{remark}
Going a bit beyond this, we will now propose a possible mirror dual to the
$(A_{2d-1})$-singularity. Let $X$
be a variety which can be degenerated to some $Y$ with an
$(A_{2d-1})$-singularity, and let $v_1,\dots,v_{2d-1} \in H_n(X)$ be the corresponding vanishing
cycles. The signs will be fixed in such a way that $v_i \cdot v_{i+1} = 1$ for all $i$. We
impose two additional conditions. One is that $Y$ should have a partial smoothing $Y'$ (equivalently,
$X$ a partial degeneration) having $d$ ODPs, built according to the local model
\[
x^2+y^2+z^2+\prod_{i=1}^d(t-\epsilon_i)^2=0
\]
with the $\epsilon_i$s distinct and small. This means that in the $(A_{2d-1})$-configuration of
vanishing cycles in $X$, one can degenerate the 1st, 3rd, \dots, (2d-1)st to ODPs. The second
additional condition is that $Y'$ should admit a resolution $Z'$ of the standard kind
considered above. Then, according to Remark \ref{remark:resolution}, there is at least one
relation between $v_1,v_3,\dots,v_{2d-1}$. In fact, since the intersection matrix of all $v_i$
has only a one-dimensional nullspace, there must be precisely one relation.
\begin{remark}
This relation is in fact $v_1 + v_3 + \dots + v_{2d-1} = 0$. The corresponding situation on
$Z'$ is that all the exceptional $\mathbb P^1$s are homologous. This should not be too
surprising: they can be moved back together to give the $d$-times thickened $\mathbb P^1$ in
the resolution of the original $A_{2d-1}$-singularity that one gets by taking the $d$-fold
branch cover $t \mapsto t^d$ of the resolution of the ODP $x^2+y^2+z^2+t^2=0$. We note in
passing that out of the $2^d$ possible ways of resolving the ODPs in $Y'$ (differing by flops)
at most two can lead to an algebraic manifold, since an
exceptional $\mathbb P^1$ cannot be homologous to minus another one.
\end{remark}
In view of our previous discussion, we expect that the mirror $\widehat{X}$ of $X$ admits a
contraction $\widehat{X} \rightarrow \widehat{Y}'$ to a variety with $d$ ODPs; any smoothing
$\widehat{Z}'$ of $\widehat{Y}'$ should contain $d$ vanishing cycles with $(d-1)$ relations
between them. By \eqref{eq:betti-numbers} these give rise to a $(d-1)$-dimensional subspace of
$H_4(\widehat{X};\C) \iso H^{1,1}(\widehat{X})$. There is a natural basis for the relations
between the exceptional $\mathbb P^1$s in $Z'$, which comes from the even-numbered vanishing
cycles $v_{2i}$. The corresponding basis of the subspace of $H^{1,1}(\widehat{X})$ can be
represented by divisors $S_2,S_4,\dots,S_{2d-2}$ such that $S_{2i}$ intersects only the
$(i-1)$-th and $i$-th exceptional $\mathbb P^1$. Based on these considerations and others described below, we
make a concrete guess as to what $\widehat{X}$ looks like:
\begin{defn}
An $(\widehat{A}_{2d-1})$-configuration of subvarieties inside a smooth threefold consists of
embedded smooth surfaces $S_2$, $S_4$,\dots, $S_{2d-2}$ and curves $C_1$, $C_3$,\dots,
$C_{2d-1}$ such that
\begin{normallist}
\item the canonical sheaf of the threefold is trivial along each $S_{2i}$;
\item each $S_{2i}$ is isomorphic to ${\mathbb P^2}$ with two points blown up;
\item $S_{2i} \cap S_{2j} = \emptyset$ for $|i-j| \geq 2$;
\item $S_{2i-2}, S_{2i}$ are transverse and intersect in $C_{2i-1}$, which is a
rational curve and exceptional (i.e.\ has self-intersection $-1$) both in $S_{2i-2}$ and
$S_{2i}$.
\end{normallist}
\end{defn}
Note that the last condition implies that $C_{2i-1}$ is a $(-1,-1)$-curve in the threefold. What
we postulate is that the mirror $\widehat{X}$ contains such a configuration of subvarieties,
with the $C_{2i-1}$ being the exceptional set of the contraction $\widehat{X} \rightarrow
\widehat{Y}'$. Apart from being compatible with the informal discussion above, there are some
more feasibility arguments in favour of this proposal. Firstly, such configurations exist as
exceptional loci in crepant resolutions of singularities: Figure \ref{fig:toric} represents a
toric 3-fold with trivial canonical bundle containing such a configuration. The thick lines
represent the $C_{2i-1}$s joining consecutive surfaces $S_{2i-2},\,S_{2i}$, which are themselves
represented by the nodes. Removing these nodes and lines gives the singularity of which it is
a resolution by collapsing the whole chain of surfaces and lines; this singularity we think of
as the dual of the $(A_{2d-1})$-singularity.

\includefigure{toric}{toric.eps}{ht}

We could have deformed the $(A_{2d-1})$-singularity in $X$ differently, for instance by degenerating
an even-numbered vanishing cycle $v_{2i}$ to an ODP. This should correspond to contracting a
${\mathbb P}^1$ in $\widehat{X}$. Assuming that our guess is right, so that $\widehat{X}$
contains a $(\widehat{A}_{2d-1})$-configuration, this should be the
other exceptional curve in the $S_{2i}$ besides $C_{2i-1}$ and $C_{2i+1}$
(i.e. the line which we will call $C_{2i}\iso\mathbb P^1$
joining $C_{2i-1}$ and $C_{2i+1}$; in Figure \ref{fig:toric} these are represented by the vertical lines). Contracting these
curves while {\em not} contracting $C_{2i-1},\,C_{2i+1}$ turns $S_{2i}$ into a ${\mathbb P^1
\times \mathbb P^1}$. The whole 4-cycle $S_{2i}$ contracts to a lower dimensional cycle only when we contract another of the $\mathbb P^1$\,s in it,
leaving the final $\mathbb P^1$ (over which the surface fibres) still
uncontracted (on $X$, this corresponds to degenerating two consecutive vanishing cycles while
leaving the others finite). Thus, there are contractions of $\widehat{X}$ mirroring various
possible partial degenerations of $X$.

A final argument in favour of our proposal, and much of the motivation for it,
is that it leads to braid group actions on derived
categories of coherent sheaves. These are of interest in themselves,
independently of whether or not they can be considered to be mirror dual to
the braid groups of Dehn twist symplectomorphisms on smoothings of
$(A_{2d-1})$-singularities.
\begin{prop} \label{th:toric-chains}
Let $X$ be a smooth quasiprojective threefold, and $S_2$, $S_4$,\dots,$S_{2d-2}$, $C_1$,
$C_3$,\dots,$C_{2d-1}$ an $(\widehat{A}_{2d-1})$-configuration of subvarieties in $X$. Then
taking ${\mathcal E_i} = \O_{C_i}$ if $i$ is odd, or $\O_{S_i}$ if $i$ is even, gives an
$(A_{2d-1})$-configuration $(\mathcal E_1, \mathcal E_2,\dots, \mathcal E_{2d-1})$ of
spherical objects in $D^b(X)$. \qed
\end{prop}
The assumption that the $S_i$ are ${\mathbb P}^2$s with two points blown up can be weakened
considerably for this result to hold; any other rational surface will do. Proposition \ref{th:toric-chains} is a
three-dimensional analogue of Example \ref{ex:-2spheres} and hence, as a comparison with our
discussion of $K3$ surfaces shows, possibly too naive from the mirror symmetry point of view.
There is an alternative way of constructing spheri\-cal objects, closer to
Proposition \ref{th:k3}.
\begin{prop} \label{th:alternative-chains}
Let $X$ be a smooth projective threefold which is Calabi-Yau in the strict sense, containing
an $(\widehat{A}_{2d-1})$-configuration as in the previous Proposition. Take rational
curves $L_{2i}$ in $S_{2i}$ such that $L_{2i} \cap C_{2j+1} = \emptyset$ for all $i,j$
(the inverse image of the generic line in $\mathbb P^2$ is such a
rational curve in the blow-up of $\mathbb P^2$). Choose
\[
{\mathcal E}_i =
\begin{cases}
 \O_{C_i}(-1) \text{\em{} or } {\mathcal J}_{C_i} & \text{\em if $i$ is odd,} \\
 \O_{S_i}(-L_i) \text{\em{} or } \O_X(-S_i) & \text{\em if $i$ is even.}
\end{cases}
\]
Then the ${\mathcal E_i}$, $i = 1,2,\dots,2d-1$, form an $(A_{2d-1})$-configuration of
spherical objects in $D^b(X)$. \qed
\end{prop}
Here $\O_{C_i}(-1)$ is shorthand for $\iota_*(\O_{\mathbb P^1}(-1))$ where $\iota: {\mathbb
P}^1 \longrightarrow X$ is some embedding with image $C_i$, and $\O_{S_i}(-L_i)$ should be
interpreted in the same way. As in Lemma \ref{th:k3}, the choice
of ${\mathcal E}_i$ can be made independently for each $i$.

There are many other interesting configurations of spherical objects which arise in connection
with threefold singularities. Their mirror symmetry interpretations are mostly unclear. For
instance, a slight variation of the situation above yields braid group actions built only from
structure sheaves of surfaces:
\begin{prop}
Let $X$ be a smooth quasiprojective threefold, and $S_1$, $S_2$, \dots, $S_m$ a chain of
smooth embedded rational surfaces in $X$ with the following properties: $S_i \cap S_{i+1}$ is
transverse and consists of one rational curve, whose self-intersection in $S_i$ and $S_{i+1}$
is either zero or $-2$; $S_i \cap S_j = \emptyset$ for $|i-j| \geq 2$; and $\o_X|S_i$ is
trivial. Then $\mathcal E_i = \O_{S_i}$ is an $(A_m)$-configuration of spherical objects in
$D^b(X)$. \qed
\end{prop}
The conditions actually imply that every intersection $S_i \cap S_{i+1}$ is a rational curve
with normal bundle $\iso \O_{{\mathbb P}^1} \oplus \O_{{\mathbb P}^1}(-2)$ in $X$. Note also
that the presence of rational curves with self-intersection zero forces at least every second
of the surfaces $S_i$ to be fibred over ${\mathbb P^1}$. Configurations of this kind are the exceptional loci
of crepant resolutions of suitable toric singularities.

In a different direction, Nakamura's resolutions of abelian quotient singularities using Hilbert
schemes of clusters, with their toric representations as tessellations of hexagons
\cite{nakamura99} \cite{craw-reid99}, lead to situations similar to Proposition
\ref{th:toric-chains}. The nodes of the hexagons in Figure \ref{fig:hex} represent surfaces
that are the blow-ups of $\mathbb P^1\times\mathbb P^1$ in two distinct points; the six lines
emanating from a node represent the six exceptional $\mathbb P^1$s in the surface, in which it
intersects the other surfaces represented by the other nodes that the lines join.

\includefigure{hex}{hex.eps}{ht}

The structure sheaves of these curves and surfaces give rise to twists on the derived category
satisfying braid relations according to the Dynkin-type diagram obtained by adding a vertex in
the middle of each edge (these added vertices represent the structure sheaves of the curves -- see Figure \ref{fig:hex}).
The McKay correspondence (see section \ref{subsec:generalisations}) translates this into a
group of twists on the equivariant derived category of the threefold on which the finite group
acted.

\section{Faithfulness\label{sec:faithfulness}}

\subsection{Differential graded algebras and modules\label{subsec:dgas}}

The notions discussed in this section are, for the most part, familiar ones; we collect them
here to set up the terminology, and also for the reader's convenience. A detailed exposition
of the general theory of differential graded modules can be found in \cite[section
10]{bernstein-lunts94}.

Fix a field $k$ and an integer $m \geq 1$. Take the semisimple $k$-algebra
$R=k^m$ with generators $e_1,\dots,e_m$ and relations $e_i^2 = e_i$ for all
$i$, $e_ie_j = 0$ for $i \neq j$ (so $1 = e_1 + \dots + e_m$ is the unit
element). $R$ will play the role of ground ring in the following
considerations. In particular, by a graded algebra we will always mean a
$\Z$-graded unital associative $k$-algebra $A$, together with a homomorphism
(of algebras, and unital) $\iota_A: R \longrightarrow A^0$. This equips $A$
with the structure of a graded $R$-bimodule, and the multiplication becomes a
bimodule map. For the sake of brevity, we will denote the bimodule structure by
$e_ia$ and $ae_i$ ($a \in A$) instead of $\iota_A(e_i)\,a$ resp.\
$a\,\iota_A(e_i)$. All homomorphisms $A \longrightarrow B$ between graded
algebras will be required to commute with the maps $\iota_A,\iota_B$. A
differential graded algebra (dga) $\dga = (A,d_A)$ is a graded algebra $A$
together with a derivation $d_A$ of degree one, which satisfies $d_A^2 = 0$ and
$d_A \circ \iota_A = 0$. The cohomology $H(\dga)$ of a dga is a graded algebra.
A homomorphism of dgas is called a quasi-isomorphism if it induces an
isomorphism on cohomology. Two dgas $\dga,\mathcal{B}$ are called
quasi-isomorphic if there is a chain of dgas and quasi-isomorphisms $\dga
\leftarrow \mathcal{C}_1 \rightarrow \dots \leftarrow \mathcal{C}_k \rightarrow
\mathcal{B}$ connecting them (in fact it is sufficient to allow $k = 1$, since
the category of dgas admits a calculus of fractions \cite[Lemma
3.2]{kellercyclic}). A dga $\dga$ is called formal if it is quasi-isomorphic to
its own cohomology algebra $H(\dga)$, thought of as a dga with zero
differential.

By a graded module over a graded algebra $A$ we will always mean a graded {\em right}
$A$-module. Through the map $\iota_A$, any such module $M$ becomes a right $R$-module: again,
we will write $xe_i$ ($x \in M$) instead of $x\,\iota_A(e_i)$. A differential graded module
(dgm) over a dga $\dga = (A,d_A)$ is a pair $\dgm = (M,d_M)$ consisting of a graded $A$-module
$M$ and a $k$-linear map $d_M: M \longrightarrow M$ of degree one, such that $d_M^2 = 0$ and
$d_M(xa) = (d_Mx)a + (-1)^{\deg(x)} x (d_Aa)$ for $a \in A$. The cohomology $H(\dgm)$ is a
graded module over $H(\dga)$. For instance, $\dga$ is a dgm over itself, and as such it splits
into a direct sum of dgms
\begin{equation} \label{eq:projective}
{\mathcal P}_i = (e_iA,d_A|_{e_iA}), \quad 1 \leq i \leq m.
\end{equation}
By definition, a dgm homomorphism $M \longrightarrow N$ is a homomorphism of graded modules
which is at the same time a homomorphism of chain complexes. Dgms over $\dga$ and their
homomorphisms form an abelian category $Dgm(\dga)$. One can also define chain homotopies
between dgm homomorphisms. The category $K(\dga)$ with the same objects as $Dgm(\dga)$ and
with the homotopy classes of dgm homomorphisms as morphisms, is triangulated. The translation
functor in it takes $\dgm = (M,d_M)$ to $\dgm[1] = (M[1],-d_M)$, with no change of sign in the
module structure. Exact triangles are all those isomorphic to one of the standard triangles
involving a dgm homomorphism and its cone.

Having mentioned cones, we use the opportunity to introduce a slight generalisation, which
will be used later on. Assume that one has a chain complex in $Dgm(\dga)$, namely dgms
${\mathcal C}_i$, $i \in \Z$, and dgm homomorphisms $\delta_i: {\mathcal C}_i \longrightarrow
{\mathcal C}_{i+1}$ such that $\delta_{i+1}\delta_i = 0$. Then one can form a new dgm
${\mathcal C}$ by setting $C = \bigoplus_{i\in \Z} C_i[i]$ and
\[
d_C = \begin{pmatrix}
 \dots        &                &                        & \\
 \delta_{i-1} & (-1)^i d_{C_i} &                        & \\
              & \delta_i       & (-1)^{i+1} d_{C_{i+1}} & \\
              &                & \delta_{i+1}           & \dots
\end{pmatrix}.
\]
We refer to this as collapsing the chain complex (it can also be viewed as a
special case of a twisted complex, see e.g.\ \cite{bondalkapranov}), and write
${\mathcal C} = \{\dots {\mathcal C}_i \rightarrow {\mathcal C}_{i+1} \dots\}$;
for complexes of length two, it specializes to the cone of a dgm homomorphism.

Inverting the dgm quasi-isomorphisms in $K(\dga)$ yields another triangulated category
$D(\dga)$, in which any short exact sequence of dgms can be completed to an exact triangle. As
usual, $D(\dga)$ can also be defined by inverting the quasi-isomorphisms directly in
$Dgm(\dga)$, but then the triangulated structure is more difficult to see. We call $D(\dga)$
the derived category of dgms over $\dga$.

\begin{warning*}
Even though we use the same notation as in ordinary homological algebra, the expressions
$K(\dga)$ and $D(\dga)$ have a different meaning here. In particular $D(\dga)$ is not the
derived category, in the usual sense, of $Dgm(\dga)$.
\end{warning*}

For any dga homomorphism $f: \dga \longrightarrow {\mathcal B}$ there is a `restriction of
scalars' functor $Dgm({\mathcal B}) \longrightarrow Dgm(\dga)$. This preserves homotopy
classes of homomorphisms, takes cones to cones, and commutes with the shift functors. Hence it
descends to an exact functor $K({\mathcal B}) \longrightarrow K(\dga)$. Moreover, it obviously
preserves quasi-isomorphisms, so that it also descends to an exact functor $D({\mathcal B})
\longrightarrow D(\dga)$; we will denote any of these functors by $f^*$. The next result,
taken from \cite[Theorem 10.12.5.1]{bernstein-lunts94}, shows that two quasi-isomorphic dgas
have equivalent derived categories.

\begin{theorem} \label{th:bernstein-lunts}
If $f$ is a quasi-isomorphism, $f^*: D({\mathcal B}) \longrightarrow D({\mathcal A})$ is an
exact equivalence.
\end{theorem}

Let $\dga$ be a dga. The {\em standard twist functors} $t_1,\dots,t_m$ from $Dgm(\dga)$ to
itself are defined by
\[
t_i(\dgm) = \{\dgm e_i \otimes_k \mathcal{P}_i \longrightarrow \dgm\}.
\]
The tensor product of $\dgm e_i = (M e_i, d_M|_{M e_i})$ with the dgm ${\mathcal P}_i$ of
\eqref{eq:projective} is one of complexes of $k$-vector spaces; it becomes a dgm with the
module structure inherited from ${\mathcal P}_i$. The arrow is the multiplication map $Me_i
\otimes_k e_iA \longrightarrow M$, which is a homomorphism of dgms, and we are taking its
cone. $t_i$ descends to exact functors $K(\dga) \longrightarrow K(\dga)$ and $D(\dga)
\longrightarrow D(\dga)$, for which we will use the same notation. This is straightforward for
$K(\dga)$. As for $D(\dga)$, one needs to show that $t_i$ preserves quasi-isomorphisms; this
follows from looking at the long exact sequence
\[
\cdots \rightarrow H(\dgm)e_i \otimes_k e_i H(\dga) \rightarrow H(\dgm) \rightarrow
H(t_i(\dgm)) \rightarrow \cdots
\]

\begin{lemma} \label{th:twists-and-quasiisos}
Let $f: \dga \longrightarrow \mathcal{B}$ be a quasi-isomorphism of dgas. Then the following
diagram commutes up to isomorphism, for each $1 \leq i \leq m$:
\[
\xymatrix{
 {D(\mathcal{B})} \ar[r]^{t_i} \ar[d]_{f^*} & {D({\mathcal B})} \ar[d]_{f^*} \\
 {D(\dga)} \ar[r]^{t_i} & {D(\dga)}
}
\]
\end{lemma}

\proof Let $\dgm = (M,d_M)$ be a dgm over $\mathcal{B}$. Consider the commutative diagram of
dgms over $\dga$
\[
\xymatrix{
 {\mathcal M e_i \otimes_k e_i \mathcal A} \ar[r] \ar[d]_{\id \otimes (f|e_iA)} &
 {f^*\dgm} \ar[d]^{\id} \\
 {\mathcal M e_i \otimes_k f^*(e_i \mathcal B)} \ar[r] &
 {f^*\dgm}
}
\]
The upper horizontal arrow is $m \otimes a \mapsto m\,f(a)$, and the lower one is
multiplication. The cone of the upper row is $t_i(f^*(\dgm))$, while that of the lower one is
$f^*(t_i(\dgm))$. The two vertical arrows combine to give a quasi-isomorphism between these
cones. \qed

Now let $\ab' \subset \ab$ be as in section \ref{subsec:homological-algebra}, and ${\mathfrak
K}$ the category from Definition \ref{def:k}. Let $E_1,\dots,E_m$ be objects of ${\mathfrak
K}$, and $E$ their direct sum. The chain complex of endomorphisms
\[
end(E) := hom(E,E) = \bigoplus_{1 \leq i,j \leq m} hom(E_i,E_j)
\]
has a natural structure of a dga. Multiplication is given by composition of homomorphisms.
$\iota_{end(E)}$ maps $e_i \in R$ to $\id_{E_i} \in hom(E_i,E_i)$, so that left multiplication
with $e_i$ is the projection to $hom(E,E_i)$ while right multiplication is the projection to
$hom(E_i,E)$. In the same way, for any $F \in {\mathfrak K}$, the complex $hom(E,F)$ is a dgm
over $end(E)$. The functor $hom(E,-): {\mathfrak K} \longrightarrow K(end(E))$ defined in this
way is exact, because it carries cones to cones. The objects $E_i$ get mapped to the dgms
$hom(E,E_i) = e_i\, end(E)$, which are precisely the ${\mathcal P}_i$ from
\eqref{eq:projective}. We define a functor $\Psi_E$ to be the composition
\[
{\mathfrak K} \xrightarrow{hom(E,-)} K(end(E)) \xrightarrow{\text{quotient functor}}
D(end(E)).
\]

\begin{lemma} \label{th:twists-and-psis}
Assume that $E_1,\dots,E_m$ satisfy the conditions from Definition \ref{def:twists}, so that
the twist functors $T_{E_i}$ are defined. Then the following diagram is commutative up to
isomorphism, for each $1 \leq i \leq m$:
\[
\xymatrix{
 {{\mathfrak K}} \ar[r]^{T_{E_i}} \ar[d]_{\Psi_E} & {{\mathfrak K}} \ar[d]^{\Psi_E} \\
 {D(end(E))} \ar[r]^{t_i} & {D(end(E))}
}
\]
\end{lemma}

\proof For $F \in {\mathfrak K}$, consider the commutative diagram of dgms over $end(E)$
\[
\xymatrix{
 {hom(E_i,F) \otimes_k hom(E,E_i)} \ar[r] \ar[d] & {hom(E,F)} \ar[d]_{\id} \\
 {hom(E,hom(E_i,F) \otimes_k E_i)} \ar[r] & {hom(E,F)}
 }
\]
with the following maps: the horizontal arrow in the first row is the composition, that in the
second row is induced by the evaluation map $hom(E_i,F) \otimes_k E_i \longrightarrow F$. The
left hand vertical arrow is the first of the canonical maps from \eqref{eq:trivial-morphisms},
which is a quasi-isomorphism since $hom(E_i,F)$ has finite-dimensional cohomology. The cone of
the first row is $t_i(\Psi_E(F))$, while that of the second row is $\Psi_E(T_{E_i}(F))$. The
vertical arrows combine to give a natural quasi-isomorphism between these cones. \qed

Later on, in our application, the $E_i$ occur as resolutions of objects in $D^b(\ab')$. The
next two Lemmas address the question of how the choice of resolutions affects the
construction. This is not strictly necessary for our purpose, but it rounds off the picture.

\begin{lemma} \label{th:endomorphism-dga}
Let $E_i,E_i'$ ($1 \leq i \leq m$) be objects in ${\mathfrak K}$ such that $E_i \iso E_i'$ for
all $i$. Then the dgas $end(E)$ and $end(E')$ are quasi-isomorphic.
\end{lemma}

\proof Choose for each $i$ a map $g_i: E_i \longrightarrow E_i'$ which is a chain homotopy
equivalence. Set $C_i = \Cone(g_i)$, and let $C$ be the direct sum of these cones; this is the
same as the cone of $g = g_1 \oplus \dots \oplus g_m$. Let $end(C)$ be the endomorphism dga of
$C_1,\dots,C_m$. An element of $end(C)$ of degree $r$ is a matrix
\[
\phi =
\begin{pmatrix}
\phi_{11} & \phi_{12} \\ \phi_{21} & \phi_{22}
\end{pmatrix}
\]
with $\phi_{11} \in hom^r(E,E)$, $\phi_{21} \in hom^{r-1}(E,E')$, $\phi_{12} \in
hom^{r+1}(E',E)$, $\phi_{22} \in hom^r(E',E')$. The differential in $end(C)$ maps $\phi$ to
\[
\begin{pmatrix}
 -d_E\phi_{11} + (-1)^r \phi_{11} d_E - (-1)^r \phi_{12} g &
 -d_E \phi_{12} - (-1)^r \phi_{12} d_E \\
 g \phi_{11} - (-1)^r \phi_{22} g + d_{E'} \phi_{21} + (-1)^r \phi_{21} d_E &
 g \phi_{12} + d_{E'} \phi_{22} - (-1)^r \phi_{22} d_{E'}
\end{pmatrix}
\]
Let ${\mathcal C} \subset end(C)$ be the subalgebra of matrices which are lower-triangular
($\phi_{12} = 0$). The formula above shows that this is closed under the differential, and
hence again a dga. The projection $\pi_2: {\mathcal C} \longrightarrow end(E')$, $\pi_2(\phi)
= \phi_{22}$, is a homomorphism of dgas. Its kernel is isomorphic (as a complex of $k$-vector
spaces, and up to a shift) to the cone of the composition with $g$ map
\[
g \circ - : hom(E,E) \longrightarrow hom(E,E').
\]
Since $g$ is a homotopy equivalence this cone is acyclic, so that $\pi_2$ is a
quasi-isomorphism of dgas. A similar argument shows that the projection $\pi_1: {\mathcal C}
\longrightarrow end(E)$, $\pi_1(\phi) = (-1)^{\deg(\phi)}\phi_{11}$, is a quasi-isomorphism of
dgas. The two maps together prove that $end(E)$ and $end(E')$ are quasi-isomorphic. \qed

As a consequence of this and Theorem \ref{th:bernstein-lunts}, the categories $D(end(E))$ and
$D(end(E'))$ are equivalent. Actually, we have shown a more precise statement: any choice of
$g_i: E_i \longrightarrow E_i'$ yields, up to isomorphism of functors, an exact equivalence
$(\pi_2^*)^{-1}\pi_1^*: D(end(E)) \longrightarrow D(end(E'))$. We will now see that this
equivalence is compatible with the functors $\Psi_E,\Psi_{E'}$.

\begin{lemma} \label{th:endomorphism-functor}
In the situation of Lemma \ref{th:endomorphism-dga}, $(\pi_2^*)^{-1}\pi_1^* \circ \Psi_E \iso
\Psi_{E'}$.
\end{lemma}

\proof The obvious short exact sequence $0 \rightarrow E' \rightarrow C \rightarrow E[1]
\rightarrow 0$ induces, for any $F \in {\mathfrak K}$, a short exact sequence of dgms over
${\mathcal C}$
\[
0 \longrightarrow \pi_1^*hom(E,F)[-1] \longrightarrow hom(C,F) \longrightarrow
\pi_2^*hom(E',F) \longrightarrow 0.
\]
In the derived category $D({\mathcal C})$, this short exact sequence can be completed to an
exact triangle by a morphism
\begin{equation} \label{eq:connecting-map}
\pi_2^*hom(E',F) \longrightarrow \pi_1^*hom(E,F).
\end{equation}
One can define such a morphism explicitly by replacing the given sequence with a (canonically
constructed) quasi-isomorphic one, for which the corresponding morphism can be realized by an
actual homomorphism of dgms; compare \cite[Proposition III.3.5]{gelfand-manin}. The advantage
of this explicit construction is that \eqref{eq:connecting-map} is now natural in $F$. Since
$C$ is a contractible complex, $hom(C,F)$ is acyclic, which implies that
\eqref{eq:connecting-map} is an isomorphism in $D({\mathcal C})$ for any $F$. This shows that
the diagram
\[
\xymatrix{
 & {{\mathfrak K}} \ar[dl]_{\Psi_E} \ar[dr]^{\Psi_{E'}} & \\
 {D(end(E))} \ar[r]^-{\pi_1^*} & D({\mathcal C}) & {D(end(E'))} \ar[l]_-{\pi_2^*}
}
\]
commutes up to isomorphism, as desired. \qed

\subsection{Intrinsic formality\label{subsec:formality}}

Applications of dg methods to homological algebra often hinge on constructing a chain of
quasi-isomorphisms connecting two given dgas. For instance, in the situation explained in the
previous section, one can try to use the dga $end(E)$ to study the twists $T_{E_i}$ via the
functor $\Psi_E$. What really matters for this purpose is only the quasi-isomorphism type of
$end(E)$. In general, quasi-isomorphism type is a rather subtle invariant. However, there are
some cases where the cohomology already determines the quasi-isomorphism type.

\begin{defn} \label{def:formality}
A graded algebra $A$ is called intrinsically formal if any two dgas with cohomology $A$ are
quasi-isomorphic; or equivalently, if any dga ${\mathcal B}$ with $H({\mathcal B}) \iso A$ is
formal.
\end{defn}

For instance, one can show easily that any graded algebra $A$ concentrated in degree zero is
intrinsically formal (this particular example can be viewed as the starting point for
Rickard's theory of derived Morita equivalences \cite{rickard89}, as recast in dga language by
Keller \cite{keller93}). However, our intended application is to algebras of a rather
different kind.

An augmented graded algebra is a graded algebra $A$ together with a graded algebra
homomorphism $\epsilon_A: A \longrightarrow R$ which satisfies $\epsilon_A \circ \iota_A =
\id_R$. Its kernel is a two-sided ideal, called the augmentation ideal; we write it as $A^+
\subset A$. A special case is when $A$ is connected, which means $A^i = 0$ for $i<0$ and
$\iota_A: R \longrightarrow A^0$ is an isomorphism; then there is of course a unique
augmentation map, and $A^+$ is the subspace of elements of positive degree.

\begin{theorem} \label{th:formality}
Let $A$ be an augmented graded algebra. If $HH^q(A,A[2-q]) = 0$ for all $q>2$, then $A$ is
intrinsically formal.
\end{theorem}

We remind the reader that the Hochschild cohomology $HH^*(A,M)$ of a graded $A$-bimodule $M$
is the cohomology of the cochain complex
\begin{align*}
 & C^q(A,M) = \Hom_{R-R}(\overbrace{A^+ \otimes_R \dots \otimes_R A^+}^{q},M),\\
 & (\partial^q \phi)(a_1,\dots,a_{q+1}) = (-1)^{\epsilon} a_1 \phi(a_2,\dots,a_{q+1}) + \\
 & + \textstyle{\sum_{i=1}^q (-1)^{\epsilon_i} \phi(a_1,\dots,a_ia_{i+1},\dots,a_{q+1}) -
(-1)^{\epsilon_q} \phi(a_1,\dots,a_q) a_{q+1}},
\end{align*}
where $\Hom_{R-R}$ denotes homomorphisms of graded $R$-bimodules (by
definition, these are homomorphisms of degree zero). The signs are $\epsilon =
q\deg(a_1)$, $\epsilon_i = \deg(a_1) + \cdots + \deg(a_i) - i$. The bimodules
relevant for our application are $M = A[s]$ with the left multiplication
twisted by a sign: $a \cdot x \cdot a' = (-1)^{s \deg(a)} axa'$ for $a,a' \in
A$ and $x \in M$. Note that the chain complex $C^*(A,A[s])$ depends on $s$, so
that the cohomology groups which occur in the Theorem above belong to different
complexes.

We will give a proof of Theorem \ref{th:formality} for lack of an accessible
reference, and also because our framework (in which dgas may be nonzero in
positive and negative degrees) differs slightly from the usual one. However,
the result is by no means new. Originally, the phenomenon of intrinsic
formality was discovered by Halperin and Stasheff \cite{halperin-stasheff79} in
the framework of commutative dgas. They constructed a series of obstruction
groups, whose vanishing implies intrinsic formality. Later Tanr{\'e}
\cite{tanre86} identified these obstruction groups as Harrison cohomology
groups. To the best of our knowledge, the non-commutative version, in which
Hochschild cohomology replaces Harrison cohomology, is due to Kadeishvili
\cite{kadeishvili88}, who also realized the importance of $A_\infty$-algebras
in this context. A general survey of $A_\infty$-algebras and applications is
\cite{kellersurvey}. It is difficult to find a concrete counterexample, but
apparently Theorem \ref{th:formality} is not true without the augmentedness
assumption. This is related to a fundamental problem, which is that the notion
of $A_\infty$-algebra with unit is not homotopy invariant (there is no
`homological perturbation Lemma' for it).

Let $A$ be an augmented graded algebra and $\mathcal{B} = (B,d_B)$ a dga. An {\em
$A_\infty$-morphism} $\gamma: A \longrightarrow \mathcal{B}$ is a sequence of maps of graded
$R$-bimodules $\gamma_q \in \Hom_{R-R}((A^+)^{\otimes_R q},B[1-q])$, $q \geq 1$, satisfying
the equations
\begin{equation} \tag{$E_q$}
\begin{split}
d_B \gamma_q(a_1,\dots,a_q) = \textstyle{\sum_{i=1}^{q-1}} (-1)^{\epsilon_i} \big(
\gamma_{q-1}(a_1,\dots,a_ia_{i+1},\dots,a_q) - \\
\gamma_i(a_1,\dots,a_i)\gamma_{q-i}(a_{i+1},\dots,a_q) \big).
\end{split}
\end{equation}
The $\epsilon_i$ are as in the definition of $HH^*(A,M)$ above. The first two of these
equations are
\begin{align}
\tag{$E_1$} & d_B \gamma_1(a_1) = 0,\\ \tag{$E_2$} & d_B \gamma_2(a_1,a_2) =
(-1)^{\deg(a_1)-1}(\gamma_1(a_1a_2) - \gamma_1(a_1)\gamma_1(a_2)).
\end{align}
This means that $\gamma_1$, which needs not be a homomorphism of algebras, nevertheless
induces a multiplicative map $(\gamma_1)_*: A^+ \rightarrow H(\mathcal{B})$. In a sense, the
non-multiplicativity of $\gamma_1$ is corrected by the higher order maps $\gamma_q$, so that
$A_\infty$-morphisms are `approximately multiplicative maps'.

From a more classical point of view, one can see $A_\infty$-morphisms simply as a convenient
way of encoding dga homomorphisms from a certain large dga canonically associated to $A$, a
kind of `thickening of $A$'. Consider $V = A^+[1]$ as a graded $R$-bimodule, and let $T^+V =
\bigoplus_{q \geq 1} V^{\otimes_R q}$ be its tensor algebra, without unit. We will write
$\leftsc a_1,\dots,a_q \rightsc \in T^+V$ instead of $a_1 \otimes \dots \otimes a_q$. Now
consider $W = T^+V[-1]$ as a graded $R$-bimodule in its own right, and form its tensor algebra
with unit $TW = R \oplus \bigoplus_{r \geq 1} V^{\otimes_R r}$. The elements of $TW$ (apart
from $R \subset TW$) are linear combinations of expressions of the form
\[
x = \leftsc a_{11}, a_{12}, \dots, a_{1,q_1} \rightsc \otimes \dots \otimes \leftsc
a_{r1},\dots,a_{r,q_r} \rightsc
\]
with $r > 0$, $q_1,\dots,q_r > 0$, and $a_{ij} \in A^+$. The degree of such an expression is
$\deg_{TW}(x) = \sum_{ij} \deg_A(a_{ij}) - \sum_i q_i + r$. One defines a dga ${\mathcal X} =
(X,d_X)$ by taking $X = TW$ with the tensor multiplication, and $d_X$ to be the derivation
which acts on elements of $W$ as follows:
\begin{multline*}
 d_X \leftsc a_1, \dots, a_q \rightsc =
 \textstyle{\sum_{i=1}^{q-1}} (-1)^{\epsilon_i}
 \big(\leftsc a_1,\dots,a_ia_{i+1},\dots,a_q \rightsc - \\
 - \leftsc a_1,\dots,a_i \rightsc \otimes \leftsc a_{i+1},\dots,a_q \rightsc \big).
\end{multline*}
The passage from $A$ to $\mathcal{X}$ is usually written as composition of the bar and cobar
functors, which go from augmented dg algebras to dg coalgebras and back, see e.g.
\cite{moore70}. We can now make the above-mentioned connection with $A_\infty$-morphisms.

\begin{lemma} \label{th:ainfty-morphisms}
For any $A_\infty$-morphism $\gamma: A \longrightarrow {\mathcal B}$ one can define a dga
homomorphism $\Gamma: {\mathcal X} \longrightarrow {\mathcal B}$ by setting $\Gamma|R$ to be
the unit map $\iota_B$, and $\Gamma(\leftsc a_1,\dots,a_q \rightsc) = \gamma_q(a_1, \dots,
a_q)$. $\Gamma$ is a quasi-isomorphism iff $\iota_B \oplus \gamma_1$ induces an isomorphism
between $R \oplus A^+ \iso A$ and $H({\mathcal B})$.
\end{lemma}

\proof The first part follows immediately from comparing the equations $(E_q)$ with the
definition of the differential $d_X$. As for the second part, a classical computation due to
Moore \cite[Th{\'e}or{\`e}me 6.2]{moore59} \cite{moore70} shows that the inclusion $R \oplus
A^+ \hookrightarrow \ker\,d_X$ induces an isomorphism $R \oplus A^+ \iso H({\mathcal X})$.
This implies the desired result. \qed

As a trivial example, let $\dga = (A,0)$ be the dga given by $A$ with zero differential, and
take the $A_\infty$-morphism $\gamma: A \longrightarrow \dga$ given by $\gamma_1 = \id: A^+
\longrightarrow A$, $\gamma_q = 0$ for all $q \geq 2$. Then Lemma \ref{th:ainfty-morphisms}
shows that the corresponding map $\Gamma: {\mathcal X} \longrightarrow {\mathcal A}$ is a
quasi-isomorphism of dgas.

The next Lemma is an instance of `homological perturbation theory', see e.g.\
\cite{gugenheim-lambe-stasheff90}. Let $A$ be an augmented graded algebra, ${\mathcal B}$ be a
dga, and $\phi: A \longrightarrow H({\mathcal B})$ a homomorphism of graded algebras. This
makes the cohomology $H({\mathcal B})$ into a graded $A$-bimodule.

\begin{lemma} \label{th:perturbation-theory}
Assume that $HH^q(A,H({\mathcal B})[2-q]) = 0$ for all $q>2$. Then there is an
$A_\infty$-morphism $\gamma: A \longrightarrow {\mathcal B}$ such that the induced map
$(\gamma_1)_*: A^+ \longrightarrow H({\mathcal B})$ is equal to $\phi|A^+$.
\end{lemma}

\proof Choose a map of graded $R$-bimodules $\gamma_1: A^+ \longrightarrow \ker\, d_B \subset
B$ which induces $\phi|A^+$. Since $\gamma_1$ is multiplicative on cohomology, we can find a
map $\gamma_2$ such that $(E_2)$ is satisfied. From here onwards the construction is
inductive. Suppose that $\gamma_1, \dots ,\gamma_{q-1}$, for some $q \geq 3$, are maps such
that $(E_1), \dots, (E_{q-1})$ hold. Denote the right hand side of equation $(E_q)$ for these
maps by $\psi: (A^+)^{\otimes_R q} \longrightarrow B[2-q]$. One can compute directly that
\begin{equation} \label{eq:c-closed}
d_B \psi(a_1,\dots,a_q) = 0
\end{equation}
for all $a_1,\dots,a_q \in A^+$, and that
\begin{equation} \label{eq:hochschild-closed}
\begin{split}
&
\begin{gathered}
 \gamma_1(a_1)\psi(a_2,\dots,a_{q+1}) + \textstyle{\sum_{i=1}^q} (-1)^{\epsilon_i}
 \psi(a_1,\dots, a_i a_{i+1},\dots,a_{q+1}) - \\
 \hspace{5cm} - (-1)^{\epsilon_q} \psi(a_1,\dots,a_q) \gamma_1(a_{q+1}) =
\end{gathered}
\\
 & = d_B \big( \textstyle{\sum_{i=1}^q} (-1)^{\epsilon_i} \gamma_i(a_1,\dots,a_i)
 \gamma_{q+1-i}(a_{i+1},\dots,a_{q+1}) \big).
\end{split}
\end{equation}
By \eqref{eq:c-closed} $\psi$ induces a map $\bar{\psi}: (A^+)^{\otimes_R q} \longrightarrow
H({\mathcal B})[2-q]$, which is just an element of the Hochschild chain group
$C^q(A,H({\mathcal B})[2-q])$. Equation \eqref{eq:hochschild-closed} says that $\bar{\psi}$ is
a Hochschild cocycle. By assumption there is an $\bar{\eta} \in C^{q-1}(A,H({\mathcal
B})[2-q])$ such that $\partial^{q-1}\bar{\eta} = \bar{\psi}$. Choose any map of graded
$R$-bimodules $\eta: (A^+)^{\otimes_R q-1} \longrightarrow (\ker\, d_B)[1-q]$ which induces
$\bar{\eta}$, and set $\gamma_{q-1}^{\mathrm{new}} = \gamma_{q-1} - \eta$. The equations
$(E_1),\dots,(E_{q-1})$ will continue to hold if one replaces $\gamma_{q-1}$ by
$\gamma_{q-1}^{\mathrm{new}}$. Moreover, if $\psi^{\mathrm{new}}$ denotes the r.h.s. of
$(E_q)$ after this replacement, one computes that
\begin{equation} \label{eq:hochschild-boundary}
\begin{gathered}
(\psi - \psi^{\mathrm{new}})(a_1,\dots,a_q) = (-1)^{\deg(a_1)}
\gamma_1(a_1)\eta(a_2,\dots,a_{q-1}) + \\ + \textstyle{\sum_{i=1}^{q-1}} (-1)^{\epsilon_i}
\eta(a_1,\dots,a_ia_{i+1},\dots,a_q) - (-1)^{\epsilon_q} \eta(a_1,\dots,a_{q-1})\gamma_1(a_q).
\end{gathered}
\end{equation}
This means that $\bar{\psi}^{\mathrm{new}} = \bar{\psi} - \partial^{q-1} \bar{\eta} = 0$.
Clearly, the vanishing of $\bar{\psi}^{\mathrm{new}}$ ensures that one can extend the sequence
$\gamma_1, \dots,\gamma_{q-2},\gamma_{q-1}^{\mathrm{new}}$ by a map $\gamma_q$ such that
$(E_q)$ holds. This completes the induction step.

Note that in the $q$-th step, only the $(q-1)$-st of the given maps $\gamma_i$ is changed.
Therefore the sequence which we construct does indeed converge to an $A_\infty$-morphism
$\gamma$. \qed

\proof[Proof of Theorem \ref{th:formality}] Let $\mathcal{B}$ be a dga whose cohomology
algebra is isomorphic to $A$. Choose an isomorphism $\phi: A \longrightarrow H({\mathcal B})$.
By Lemma \ref{th:perturbation-theory} there is an $A_\infty$-morphism $\gamma: A
\longrightarrow {\mathcal B}$ such that $\gamma_1$ induces $\phi|A^+$. This obviously means
that $(\iota_B \oplus \gamma_1)_*: R \oplus A^+ \longrightarrow H({\mathcal B})$ is an
isomorphism. Hence, by Lemma \ref{th:ainfty-morphisms} the induced map $\Gamma: {\mathcal X}
\longrightarrow {\mathcal B}$ is a quasi-isomorphism of dgas. We have already seen that there
is a quasi-isomorphism ${\mathcal X} \longrightarrow \dga = (A,0)$. This shows that ${\mathcal
B}$ is quasi-isomorphic to $\dga$, hence formal. \qed

\subsection{The graded algebras $A_{m,n}$\label{subsec:a}}

We assume from now on that $m \geq 2$; this assumption will be retained throughout this
section and the following one. In addition, choose an $n \geq 1$.

Let $\Gamma$ be a quiver (an oriented graph) with vertices numbered $1,\dots,m$, and with a
`degree' (an integer label) attached to each edge. One can associate to it a graded algebra
$k[\Gamma]$, the path algebra, as follows. As a $k$-vector space $k[\Gamma]$ is freely
generated by the set of all paths (not necessarily closed, of arbitrary length $\geq 0$) in
$\Gamma$. The degree of a path is the sum of all `degrees' of the edges along which it runs.
The product of two paths is their composition if the endpoint of the first one coincides with
the starting point of the second one, and zero otherwise. The map $\iota_{k[\Gamma]}: R
\longrightarrow (k[\Gamma])^0$ maps $e_i$ to the path of length zero at the $i$-th vertex.
\includefigure{quiver}{quiver.eps}{hb}%

The example we are interested in is the quiver $\Gamma_{m,n}$ shown in Figure \ref{fig:quiver}.
Paths of length $l \geq 0$ in this quiver correspond to $(l+1)$-tuples $(i_0|\dots|i_l)$ with
$i_\nu \in \{1,\dots,m\}$ and $|i_{\nu+1}-i_{\nu}| = 1$. The product of two paths in
$k[\Gamma_{m,n}]$ is given by $(i_0|\dots|i_l)(i'_0|\dots|i'_{l'}) =
(i_0|\dots|i_l|i'_1|\dots|i'_{l'})$ if $i_l = i'_0$, or zero otherwise. The grading is $\deg\,
(i) = 0$, $\deg\, (i|i+1) = d_i$, $\deg\, (i+1|i) = n - d_i$, where we set
\begin{equation} \label{eq:d-i}
d_i =
 \begin{cases} \half n & \text{if $n$ is even,} \\
 \half(n+(-1)^i) & \text{if $n$ is odd.}
 \end{cases}
\end{equation}
We introduce a two-sided homogeneous ideal $J_{m,n} \subset k[\Gamma_{m,n}]$ as follows. If $m
\geq 3$ then $J_{m,n}$ is generated by $(i|i-1|i) - (i|i+1|i)$, $(i-1|i|i+1)$ and
$(i+1|i|i-1)$ for all $i = 2,\dots,m-1$; in the remaining case $m = 2$, $J_{m,n}$ is generated
by $(1|2|1|2)$ and $(2|1|2|1)$. Now define $A_{m,n} = k[\Gamma_{m,n}]/J_{m,n}$. This is again
a graded algebra. It is finite-dimensional over $k$; an explicit basis is given by the
$(4m-2)$ elements
\begin{equation} \label{eq:a-basis}
\left\{
\begin{split}
 & (1),\, \dots,\, (m),\\
 & (1|2),\, \dots,\, (m-1|m),\\
 & (2|1),\, \dots,\, (m|m-1),\\
 & (1|2|1),\, (2|3|2) = (2|1|2),\, \dots,\, (m-1|m|m-1) = \\
 &  \qquad \qquad = (m-1|m-2|m-1),\, (m|m-1|m).
\end{split}
\right.
\end{equation}
Here we have used the same notation for elements of $k[\Gamma_{m,n}]$ and their images in
$A_{m,n}$. We will continue to do so in the future, in particular $(i|i\pm 1|i)$ will be used
to denote the image of both $(i|i+1|i)$ and $(i|i-1|i)$ in $A_{m,n}$.

We will now explain why these algebras are relevant to our problem. Let ${\mathfrak K}$ be a
category as in Definition \ref{def:k} and $E_1,\dots,E_m \in {\mathfrak K}$ an
$(A_m)$-configuration of $n$-spherical objects.

\begin{lemma} \label{th:end-algebra}
Suppose that for each $i = 1,\dots,m-1$ the one-dimensional space $\Hom^*(E_{i+1},E_i)$ is
concentrated in degree $d_i$. Then the cohomology algebra of the dga $end(E)$ is isomorphic to
$A_{m,n}$.
\end{lemma}

We should say that the assumption on $\Hom^*(E_{i+1},E_i)$ is not really restrictive since,
given an arbitrary $(A_m)$-configuration, it can always be achieved by shifting each $E_i$
suitably.

\proof Since each $E_i$ is $n$-spherical, the pairings
\begin{equation} \label{eq:two-pairings}
\begin{aligned}
 \Hom^*(E_{i+1},E_i) \otimes \Hom^*(E_i,E_{i+1})
 & \longrightarrow \Hom^n(E_i,E_i) \iso k, \\
 \Hom^*(E_i,E_{i+1}) \otimes \Hom^*(E_{i+1},E_i)
 & \longrightarrow \Hom^n(E_{i+1},E_{i+1}) \iso k
\end{aligned}
\end{equation}
are nondegenerate for $i = 1,\dots,m-1$. Hence $\Hom^*(E_i,E_{i+1}) \iso k$ is concentrated in
degree $n-d_i$. Choose nonzero elements $\alpha_i \in \Hom^*(E_{i+1},E_i)$ and $\beta_i \in
\Hom^*(E_i,E_{i+1})$. Then, again because of the nondegeneracy of \eqref{eq:two-pairings}, one
has
\begin{equation} \label{eq:first-relation}
\alpha_i\beta_i = c_i \, (\beta_{i-1}\alpha_{i-1})
\end{equation}
in $\Hom^*(E,E)$ for some nonzero constants $c_2,\dots,c_{m-1} \in k$. Without changing
notation, we multiply each $\beta_i$ with $c_2c_3\dots c_i$; then the same equations
\eqref{eq:first-relation} hold with all $c_i$ equal to $1$. Since $\Hom^*(E_i,E_j) = 0$ for
all $|i-j| \geq 2$, we also have $\beta_i\beta_{i-1} = 0$, $\alpha_{i-1}\alpha_i = 0$ for all
$i = 2,\dots,m-1$. If $m \geq 3$ then this shows that there is a homomorphism of graded
algebras $A_{m,n} \longrightarrow \Hom^*(E,E)$ which maps $(i)$ to $\id_{E_i}$, $(i|i+1)$ to
$\alpha_i$, and $(i+1|i)$ to $\beta_i$. One sees easily that this is an isomorphism. In the
remaining case $m = 2$ one has to consider
\begin{equation} \label{eq:m2relation}
\beta_1\alpha_1\beta_1 \in \Hom^{2n-d_1}(E_1,E_2), \quad \alpha_1\beta_1\alpha_1 \in
\Hom^{n+d_1}(E_2,E_1).
\end{equation}
By assumption $\Hom^*(E_1,E_2)$ is concentrated in degree $n-d_1 < 2n-d_1$, and
$\Hom^*(E_2,E_1)$ is concentrated in degree $d_1 < n+d_1$. Hence both elements in
\eqref{eq:m2relation} are zero, which allows one to define $A_{m,n} \longrightarrow
\Hom^*(E,E)$ as before. The proof that this is an isomorphism is again straightforward. \qed

An inspection of the preceding proof shows that the result remains true for any other choice
of numbers $d_i$ in the definition of $A_{m,n}$. Our particular choice \eqref{eq:d-i} makes
the algebra as `highly connected' as possible: $A_{m,n}/R \cdot 1$ is concentrated in degrees
$\geq [n/2]$. This will be useful in the Hochschild cohomology computations of section
\ref{subsec:a-is-formal}.

Let $\dga_{m,n}$ be the dga given by $A_{m,n}$ with zero differential. We will now consider
the properties of the functors $t_i$ on the category $D(\dga_{m,n})$.

\begin{lemma} \label{th:a-invertibility}
The functors $t_i: D(\dga_{m,n}) \longrightarrow D(\dga_{m,n})$, $1 \leq i \leq m$, are exact
equivalences.
\end{lemma}

\proof This is closely related to the parallel statements in \cite{khovanov-seidel98} and in
our section \ref{subsec:inverse}. The strategy, as in Proposition \ref{th:invertibility}, is
to introduce a left adjoint $t_i'$ of $t_i$, and then to prove that the canonical natural
transformations $\Id \longrightarrow t_it_i'$, $t_i't_i \longrightarrow \Id$ are isomorphisms.

Set $\dga = \dga_{m,n}$ and ${\mathcal Q}_i = \mathcal{P}_i[n] \in Dgm(\dga)$. Define functors
$t_i'$ ($1 \leq i \leq m$) from $Dgm(\dga)$ to itself by
\[
t_i'(\dgm) = \{\dgm \stackrel{\eta_i}{\longrightarrow} \dgm e_i \otimes_k \mathcal{Q}_i\}
\]
where $\dgm$ is placed in degree zero, and $\eta_i(x) = x(i|i\pm 1|i) \otimes (i) + x (i+1|i)
\otimes (i|i+1) + x (i-1|i) \otimes (i|i-1) + x (i) \otimes (i|i \pm 1|i)$ (in this formula,
the second term should be omitted for $i = m$ and the third term for $i = 1$; the same
convention will be used again later on). To understand why $\eta_i$ is a module homomorphism,
it is sufficient to notice that the element
\begin{equation} \label{eq:central-element}
\begin{split}
 (i|i\!\pm\! 1|i) \otimes (i) + (i\!+\!1|i) \otimes (i|i\!+\!1) +
 (i\!-\!1|i) \otimes (i|i\!-\!1) + \\ +
 (i) \otimes (i|i \!\pm\! 1|i) \in Ae_i \otimes e_iA
\end{split}
\end{equation}
is central, in the sense that left and right multiplication (with respect to the obvious
$A$-bimodule structure of $Ae_i \otimes e_iA$) with any $a \in A$ have the same effect on it.
The same argument as for $t_i$ shows that $t_i'$ descends to exact functors on $K(\dga)$ and
$D(\dga)$. For any $\dgm \in Dgm(\dga)$ consider the complex of dgms
\[
 {\mathcal C}_{-1} = \dgm e_i \otimes {\mathcal P}_i
 \stackrel{\delta_{-1}}{\longrightarrow}
 {\mathcal C}_0 = \dgm \oplus (\dgm e_i \otimes e_iAe_i \otimes {\mathcal Q}_i)
 \stackrel{\delta_0}{\longrightarrow}
 {\mathcal C}_1 = \dgm e_i \otimes {\mathcal Q}_i,
\]
where $\delta_{-1}(x \otimes a) = (xa, x \otimes a(i|i \pm  1|i) \otimes (i) + x \otimes
a(i+1|i) \otimes (i|i+1) + x \otimes a(i-1|i) \otimes (i|i-1) + x \otimes a(i) \otimes (i|i\pm
1|i)) = (xa, x \otimes (i) \otimes (i|i\pm 1|i)a + x \otimes (i|i\pm 1|i) \otimes a)$ and
$\delta_0(x,y \otimes a \otimes b) = (\eta_i(x) - y a \otimes b)$. The reason why the second
expression for $\delta_{-1}$ is equal to the first one is again that the element
\eqref{eq:central-element} is central. A straightforward computation (including some tedious
sign checking) shows that the dgm ${\mathcal C}$ obtained by collapsing this complex is equal
to $t_i't_i(\dgm)$.

$e_iAe_i = k(i) \oplus k(i|i\pm 1|i)$ is simply a two-dimensional graded $k$-vector space,
nontrivial in degrees zero and $n$. Take the homomorphism of dgms
\begin{equation} \label{eq:nat-trans}
\begin{aligned}
 {\mathcal C}_0 =
 \dgm \oplus (\dgm e_i \otimes e_iAe_i \otimes {\mathcal Q}_i) & \longrightarrow \dgm, \\
 (x, y_1 \otimes (i) \otimes b_1 + y_2 \otimes (i|i\pm 1|i) \otimes b_2) & \longmapsto
 x - y_2b_2.
\end{aligned}
\end{equation}
Extending this by zero to ${\mathcal C}_{-1},{\mathcal C}_1$ yields a dgm homomorphism
$\psi_{\dgm}: {\mathcal C} = t_i't_i(\dgm) \longrightarrow \dgm$, because \eqref{eq:nat-trans}
vanishes on the image of $\delta_{-1}$. This homomorphism is surjective for any $\dgm$, and a
computation similar to that in Proposition \ref{th:invertibility} shows that the kernel is
always an acyclic dgm. Since $\psi_{\dgm}$ is natural in $\dgm$, we have indeed provided an
isomorphism $t_i't_i \iso \Id_{D(\dga)}$. The proof that $t_it_i' \iso \Id_{D(\dga)}$ is
parallel. \qed

\begin{lemma} \label{th:a-braiding}
The functors $t_i$ on $D(\dga_{m,n})$ satisfy the braid relations (up to graded natural
isomorphism):
\begin{align*}
 t_it_{i+1}t_i &\iso t_{i+1}t_it_{i+1} && \text{for } i = 1,\dots,m-1,\\
 t_it_j &\iso t_jt_i && \text{for } |i-j| \geq 2.
\end{align*}
\end{lemma}

\proof The second relation is easy (it follows immediately from the fact that $e_iA_{m,n}e_j =
0$ for $|i-j| \geq$), and we will therefore concentrate on the first one. Moreover, we will
only explain the salient points of the argument (a different version of it is described in
\cite{khovanov-seidel98} with full details). Note that the approach taken in Proposition
\ref{th:braid-two} cannot be adapted directly to the present case, since we have not developed
a general theory of twist functors on derived categories of dgms.

Set $\dga = {\dga}_{m,n}$ and ${\mathcal R}_i = {\mathcal P}_i[-n]$. For any $\dgm \in
Dgm(\dga)$ consider the complex of dgms
\begin{equation} \label{eq:boring-complex}
{\mathcal C}_{-3} \stackrel{\delta_{-3}}{\longrightarrow} {\mathcal C}_{-2}
\stackrel{\delta_{-2}}{\longrightarrow} {\mathcal C}_{-1}
 \stackrel{\delta_{-1}}{\longrightarrow} {\mathcal C}_0,
\end{equation}
where
\begin{align*}
 {\mathcal C}_{-3} &= \dgm e_i \otimes {\mathcal R}_i, \\
 {\mathcal C}_{-2} &= (\dgm e_i \otimes e_iAe_i \otimes {\mathcal P}_i) \oplus
 (\dgm e_i \otimes e_iAe_{i+1} \otimes {\mathcal P}_{i+1}) \oplus \\ & \hspace{5cm}
 \oplus (\dgm e_{i+1} \otimes e_{i+1}Ae_i \otimes {\mathcal P}_i), \\
 {\mathcal C}_{-1} &= (\dgm e_i \otimes {\mathcal P}_i)
 \oplus (\dgm e_{i+1} \otimes {\mathcal P}_{i+1}) \oplus
 (\dgm e_i \otimes {\mathcal P}_i), \\
 {\mathcal C}_0 & = \dgm \\
\intertext{and}
 \delta_{-3} & : \; (x \otimes a) \mapsto
 \begin{pmatrix}
  -x \otimes (i|i\!+\!1|i) \otimes a \\
  x \otimes (i|i\!+\!1) \otimes (i\!+\!1|i)a \\
  x(i|i+1) \otimes (i\!+\!1|i) \otimes a
 \end{pmatrix}, \\
 \delta_{-2} & : \;
 \begin{pmatrix}
  x_1 \otimes a_1 \otimes b_1 \\
  x_2 \otimes (i|i\!+\!1) \otimes b_2 \\
  x_3 \otimes (i\!+\!1|i) \otimes b_3
 \end{pmatrix}
 \mapsto
 \begin{pmatrix}
  x_1 \otimes a_1b_1 + x_2 \otimes (i|i\!+\!1)b_2 \\
  - x_2(i|i\!+\!1) \otimes b_2 + x_3 \otimes (i\!+\!1|i)b_3 \\
  - x_1a \otimes b_1 - x_3(i\!+\!1|i) \otimes b_3
 \end{pmatrix}, \\
 \delta_{-1} & : \;
 \begin{pmatrix}
  x_1 \otimes a_1 \\
  x_2 \otimes a_2 \\
  x_3 \otimes x_3 \\
 \end{pmatrix}
 \mapsto x_1a_1 + x_2a_2 + x_3a_3.
\end{align*}
As in the proof of the previous Lemma, one can contract this complex to a single dgm, which is
in fact canonically isomorphic to to $t_it_{i+1}t_i(\dgm)$. Now, one can map the whole complex
\eqref{eq:boring-complex} surjectively to an acyclic complex (concentrated in degrees $-3$ and
$-2$)
\[
\dgm e_i \otimes {\mathcal R}_i \stackrel{\id}{\longrightarrow} \dgm e_i \otimes {\mathcal
R}_i.
\]
This is done by taking the identity map on ${\mathcal C}_{-3}$ together with the homomorphism
${\mathcal C}_{-2} \supset \dgm e_i \otimes e_iAe_i \otimes {\mathcal P}_i \longrightarrow
\dgm e_i \otimes {\mathcal R}_i$, $m_1 \otimes (i) \otimes b_1 + m_2 \otimes (i|i\!+\!1|i)
\otimes b_2 \mapsto -m_2 \otimes b_2$, and extending this by zero to the other summands of
${\mathcal C}_{-2}$ and to ${\mathcal C}_{-1}$, ${\mathcal C_0}$. The kernel of the dgm
homomorphism defined in this way is a certain subcomplex of \eqref{eq:boring-complex}. When
writing this down explicitly (which we will not do here) one notices that it contains an
acyclic subcomplex isomorphic to
\[
\dgm e_i \otimes {\mathcal P}_i \stackrel{\id}{\longrightarrow} \dgm e_i \otimes {\mathcal
P}_i,
\]
located in degrees $-2$ and $-1$. If one divides out this acyclic subcomplex, what remains is
the complex
\begin{equation} \label{eq:remainder}
\begin{split}
 (\dgm e_i \otimes e_iAe_{i+1} \otimes {\mathcal P}_{i+1}) \oplus (\dgm e_{i+1} \otimes
 e_{i+1}Ae_i \otimes {\mathcal P}_i)
 \stackrel{\delta_{-1}'}{\longrightarrow} \\ \longrightarrow
 (\dgm e_i \otimes {\mathcal P}_i) \oplus (\dgm e_{i+1} \otimes {\mathcal P}_{i+1})
 \stackrel{\delta_0'}{\longrightarrow}
 \dgm
\end{split}
\end{equation}
with $\delta_{-1}'(x_1 \otimes (i|i\!+\!1) \otimes b_1,x_2 \otimes (i\!+\!1|i) \otimes b_2) =
(x_1 \otimes (i|i\!+\!1) b_1 - x_2(i\!+\!1|i) \otimes b_2, -x_1(i|i\!+\!1) \otimes b_1 + x_2
\otimes (i\!+\!1|i)b_2)$, $\delta_0'(x_1 \otimes a_1, x_2 \otimes a_2) = x_1a_1 + x_2a_2$. The
remarkable fact about \eqref{eq:remainder} is that it is symmetric with respect to exchanging
$i$ and $i+1$. Indeed, one can arrive at the same complex by starting with
$t_{i+1}t_it_{i+1}(\dgm)$ and removing acyclic parts. This shows that
$t_{i+1}t_it_{i+1}(\dgm)$ and $t_it_{i+1}t_i(\dgm)$ are quasi-isomorphic for all $\dgm$. We
leave it to the reader to verify that the argument provides a chain of exact functors and
graded natural isomorphisms between them, with $t_it_{i+1}t_i$ and $t_{i+1}t_it_{i+1}$ at the
two ends of the chain. \qed

\subsection{Geometric intersection numbers\label{subsec:gin}}

Consider the weak braid group action $\rho_{m,n}: B_{m+1} \longrightarrow
\Auteq(D(\dga_{m,n}))$ generated by $t_1,\dots,t_m$. The aim of this section is prove a strong
form of faithfulness for it:

\begin{thm} \label{th:a-faithful}
Let $R_{m,n}^g$ be a functor representing $\rho_{m,n}(g)$ for some $g \in B_{m+1}$. If
$R_{m,n}^g({\mathcal P}_j) \iso {\mathcal P}_j$ for all $j$, then $g$ must be the identity
element.
\end{thm}

We begin by looking at the center of $B_{m+1}$. It is infinite cyclic and generated by an
element which, in terms of the standard generators $g_1,\dots,g_m$, can be written as
$(g_1g_2\dots g_m)^{m+1}$.

\begin{lemma} \label{th:central-action}
For any $1 \leq j \leq m$, $(t_1 t_2 \dots t_m)^{m+1}({\mathcal P}_j)$ is isomorphic to
${\mathcal P}_j[2m - (m+1)n]$ in $D(\dga_{m,n})$.
\end{lemma}

\proof For each $1 \leq j \leq m$ there is a short exact sequence of dgms
\[
0 \longrightarrow {\mathcal P}_j[-n] \stackrel{\alpha}{\longrightarrow} e_jA_{m,n}e_j
\otimes_k {\mathcal P}_j \xrightarrow{\text{multiplication}} {\mathcal P}_j \longrightarrow 0,
\]
where $\alpha(x) = (j|j\pm 1|j) \otimes x - (j) \otimes (j|j\pm 1|j)x$. This implies that the
cone of the multiplication map, which is $t_j(\mathcal P_j)$, is isomorphic to ${\mathcal
P}_j[1-n]$ in $D(\dga_{m,n})$. Note also that $t_i(\mathcal P_j) \iso \mathcal P_j$ whenever
$|i-j| \geq 2$.

Consider the $m+1$ differential graded modules
\begin{align*}
 & \dgm_0 = \{{\mathcal P}_1[n\!-\!1] \rightarrow {\mathcal P}_2[2n\!-\!1\!-\!d_1]
 \rightarrow \dots
 \rightarrow {\mathcal P}_m[mn\!-\!1\!-\!d_1\!-\!\dots\!-\!d_{m-1}]\}, \\
 & \dgm_1 = {\mathcal P}_1, \\
 & \dgm_2 = {\mathcal P}_2[1\!-\!d_1], \\
 & \dgm_3 = {\mathcal P}_3[2\!-\!d_1\!-\!d_2], \\
 & \dots \\
 & \dgm_m = {\mathcal P}_m[m\!-\!1\!-\!d_1\!-\!\dots\!-\!d_{m-1}].
\end{align*}
The definition of $\dgm_0$ is by collapsing the complex of dgms in which ${\mathcal P}_1[n-1]$
is placed in degree zero, and where the maps are given by left multiplication with $(i+1|i)$.
We will prove that
\begin{equation} \label{eq:permutes}
\begin{cases}
 (t_1t_2\dots t_m)(\dgm_0) \iso \dgm_1, & \\
 (t_1t_2\dots t_m)(\dgm_i) \iso \dgm_{i+1} \text{ for $1 \leq i < m$,} & \\
 (t_1t_2\dots t_m)(\dgm_m) \iso \dgm_0[2m - (m+1)n], &
\end{cases}
\end{equation}
which clearly implies the desired result. By the definitions of $t_i$ and $t_i'$, the second
of which is given in the proof of Lemma \ref{th:a-invertibility}, one has
\[
t_{i+1}({\mathcal P_i}) = \{ {\mathcal P}_{i+1}[-d_i] \rightarrow {\mathcal P}_i \} \iso
t_i'(\mathcal P_{i+1})[1-d_i]
\]
where ${\mathcal P}_i$ is placed in degree zero and the arrow is left multiplication with
$(i|i+1)$. This shows that $t_it_{i+1}({\mathcal P}_i) \iso {\mathcal P}_{i+1}[1-d_i]$, and
since $t_i({\mathcal P}_j) \iso {\mathcal P}_j$ whenever $|i-j| \geq 2$, it proves the second
equation in \eqref{eq:permutes}. To verify the other two equations one computes
\begin{align*}
 & (t_1t_2\dots t_m)({\mathcal P}_m[n\!-\!1]) \iso \\
 & \iso (t_1t_2\dots t_{m-1})({\mathcal P}_m) \\
 & \iso (t_1t_2\dots t_{m-2})(\{{\mathcal P}_{m-1}[-n\!+\! d_{m-1}] \rightarrow
 {\mathcal P}_m\}) \\
 & \iso (t_1t_2\dots t_{m-3})(\{{\mathcal P}_{m-2}[-2n\!+\!d_{m-2}\!+\!d_{m-1}] \rightarrow
 {\mathcal P}_{m-1}[-n\!+\!d_{m-1}] \rightarrow {\mathcal P}_m\}) \\
 & \iso \dots \iso \dgm_0[m(1\!-\!n) \! +\! d_1\! +\! \dots\! +\! d_{m-1}] \\
\intertext{and}
 & (t_m'\dots t_2't_1')({\mathcal P}_1) \iso \\
 & \iso (t_m'\dots t_2')({\mathcal P}_1[n-1]) \\
 & \iso (t_m'\dots t_3')(\{{\mathcal P}_1[n-1] \rightarrow
 {\mathcal P}_2[2n\!-\!1\!-\!d_1]\}) \\
 & \iso (t_m'\dots t_4')(\{{\mathcal P}_1[n-1] \rightarrow
 {\mathcal P}_2[2n\!-\!1\!-\!d_1] \rightarrow
 {\mathcal P}_3[3n\!-\!1\!-\!d_1\!-\!d_2]\}) \\
 & \iso \dots \iso \dgm_0. \qed
\end{align*}

It seems likely that $(t_1 t_2 \dots t_m)^{m+1}$ is in fact isomorphic to the translation
functor $[2m - (m+1)n]$, but we have not checked this.

Before proceeding further, we need to recall some basic notions from the topology of curves on
surfaces. Let $D$ be a closed disc, and $\Delta \subset D \setminus \partial D$ a set of $m+1$
marked points. $\Diff(D,\partial D;\Delta)$ denotes the group of diffeomorphisms $f: D
\longrightarrow D$ which satisfy $f|\partial D = \id$ and $f(\Delta) = \Delta$. We write $f_0
\isotopic f_1$ for isotopy within this group. By a curve in $(D,\Delta)$ we mean a subset $c
\subset D \setminus \partial D$ which can be represented as the image of a smooth embedding
$\gamma: [0;1] \longrightarrow D$ such that $\gamma^{-1}(\Delta) = \{0;1\}$. In other words,
$c$ is an unoriented embedded path in $D \setminus \partial D$ whose endpoints lie in
$\Delta$, and which does not meet $\Delta$ anywhere else. There is an obvious notion of
isotopy for curves, denoted again by $c_0 \isotopic c_1$. For any two curves $c_0,c_1$ there
is a geometric intersection number $I(c_0,c_1) \geq 0$, which is defined by $I(c_0,c_1) =
|(c_0' \cap c_1) \setminus \Delta| + \half |(c_0' \cap c_1) \cap \Delta|$ for some $c_0'
\isotopic c_0$ which has minimal intersection with $c_1$ (this means, roughly speaking, that
$c_0'$ is obtained from $c_0$ by removing all unnecessary intersection points with $c_1$). We
refer to \cite[section 2a]{khovanov-seidel98} for the proof that this is well-defined. Once
one has shown this, the following properties are fairly obvious:
\begin{normallist}
\item[(I1)] $I(c_0,c_1)$ depends only on the isotopy classes of $c_0$ and $c_1$;
\item[(I2)] $I(c_0,c_1) = I(f(c_0),f(c_1))$ for all $f \in \Diff(D,\partial D;\Delta)$;
\item[(I3)] $I(c_0,c_1) = I(c_1,c_0)$.
\end{normallist}
Note that in general $I(c_0,c_1)$ is only a half-integer, because of the weight $1/2$ with
which the common endpoints of $c_0$ and $c_1$ contribute. The next Lemma, whose proof we omit,
is a modified version of \cite[Proposition III.16]{fathi-laudenbach-poenaru}.

\begin{lemma} \label{th:pairing}
Let $c_0,c_1$ be two curves in $(D,\Delta)$ such that $I(d,c_0) = I(d,c_1)$ for all $d$. Then
$c_0 \isotopic c_1$. \qed
\end{lemma}
\includefigure{basic-curves}{basic-curves.eps}{hb}%

From now on, fix a collection of curves $b_1,\dots,b_m$ as in Figure \ref{fig:basic-curves}, as
well as an orientation of $D$. Then one can identify $\pi_0(\Diff(D,\partial D;\Delta))$ with
the braid group by mapping the standard generators $g_1,\dots,g_m \in B_{m+1}$ to positive
half-twists along $b_1,\dots,b_m$.

\begin{lemma} \label{th:central-element}
Let $f \in \Diff(D,\partial D;\Delta)$ be a diffeomorphism which satisfies $f(b_j) \isotopic
b_j$ for all $1 \leq j \leq m$. The the corresponding element $g \in B_{m+1}$ must be of the
form $g = (g_1g_2\dots g_m)^{\nu(m+1)}$ for some $\nu \in \Z$.
\end{lemma}

\proof Since $f(b_j) \isotopic b_j$, $f$ commutes up to isotopy with the half-twist along
$b_j$, and hence with any element of $\Diff(D,\partial D;\Delta)$. This implies that $g$ is
central. \qed

The next Lemma, which is far more substantial than the previous ones, establishes a
relationship between the topology of curves in $(D,\Delta)$ and the algebraically defined
braid group action $\rho_{m,n}$.

\begin{lemma} \label{th:gin}
For $g \in B_{m+1}$, let $f \in \Diff(D,\partial D;\Delta)$ be a diffeomorphism in the isotopy
class corresponding to $g$, and $R_{m,n}^g$ a functor which represents $\rho_{m,n}(g)$. Then
\[
 \sum_{r \in \Z} \dim_k \Hom_{D(\dga_{m,n})}({\mathcal P}_i,R_{m,n}^g({\mathcal P}_j)[r]) =
 2\,I(b_i,f(b_j))
\]
for all $1 \leq i,j \leq m$.
\end{lemma}

A statement of the same kind, concerning a category and braid group action slightly different
from ours, has been proved in \cite[Theorem 1.1]{khovanov-seidel98}. In principle, the proof
given there can be adapted to our situation, but verifying all the details is a rather tedious
business. For this reason we take a slightly different approach, which is to derive the result
as stated here from its counterpart in \cite{khovanov-seidel98}. To do this, we first need to
recall the situation considered in that paper. In order to avoid confusion, objects which
belong to the setup of \cite{khovanov-seidel98} will be denoted by overlined symbols.
\includefigure{quiver0}{quiver0.eps}{hb}%

Consider the quiver $\overline{\Gamma}_m$ in Figure \ref{fig:quiver0} with vertices numbered
$0,\dots,m$ and whose edges are labelled with `degrees' zero or one. Paths of length $l$ in
$\overline{\Gamma}_m$ are described by $(l+1)$-tuples of numbers $i_0,\dots,i_l \in
\{0,\dots,m\}$; we will use the notation $\barpath{i_0|\dots|i_l}$ for them. The path algebra
$k[\overline{\Gamma}_m]$ is a graded algebra, whose ground ring is $\overline{R} = k^{m+1}$.
Let $\overline{J}_m$ be the homogeneous two-sided ideal in it generated by the elements
$\barpath{i-1|i|i+1}$, $\barpath{i+1|i|i-1}$, $\barpath{i|i+1|i} - \barpath{i|i-1|i}$ ($1 \leq
i \leq m-1$), and $\barpath{0|1|0}$. The quotient $\overline{A}_m =
k[\overline{\Gamma}_m]/\overline{J}_{m}$ is a finite-dimensional graded algebra; a concrete
basis is given by the $4m+1$ elements
\begin{equation} \label{eq:bar-basis}
\left\{
\begin{split}
 & \barpath{0},\, \dots,\, \barpath{m}, \,
 \barpath{0|1},\, \dots,\, \barpath{m\!-\!1|m} \quad \text{of degree zero, and} \\
 & \barpath{1|0},\, \dots,\, \barpath{m|m\!-\!1},\,
 \barpath{1|2|1} = \barpath{1|0|1}, \, \dots,\,
 \barpath{m\!-\! 1|m\!-\! 2|m\!-\! 1} = \\
 & \quad = \barpath{m\!-\!1|m|m\!-\!1}, \, \barpath{m|m\!-\! 1|m} \quad
 \text{of degree one.}
\end{split}
\right.
\end{equation}

$\overline{A}_m$ is evidently a close cousin of our algebras $A_{m,n}$. We will now make the
relationship precise on the level of categories. Let $\Amod$ be the abelian category of
finitely-generated graded right modules over $\overline{A}_m$, and $D^b(\Amod)$ its bounded
derived category (in contrast to the situation in section \ref{subsec:dgas}, this is the
derived category in the ordinary sense, not in the differential graded one). There is an
automorphism $\{1\}$ which shifts the grading of a module up by one. This descends to an
automorphism of $D^b(\Amod)$, which is not the same as the translation functor. In particular,
for any $X,Y \in D^b(\Amod)$ there is a bigraded vector space
\[
\bigoplus_{r_1,r_2} \Hom_{D^b(\Amod)}(X,Y\{r_1\}[r_2]).
\]
We denote by $\overline{P}_i \in \Amod$ the projective modules $\barpath{i} \overline{A}_m$,
for $0 \leq i \leq m$. Let ${\mathfrak P} \subset \Amod$ be the full subcategory whose objects
are direct sums of $\overline{P}_i\{r\}$ for $i = 1,\dots,m$ and $r \in \Z$; the important
thing is that $\overline{P}_0$ is not allowed. We write $K^b({\mathfrak P})$ for the full
subcategory of $K^b(\Amod)$ whose objects are finite complexes in ${\mathfrak P}$. This is an
abuse of notation since ${\mathfrak P}$ is not an abelian category; however, $K^b({\mathfrak
P})$ is still a triangulated category, because it contains the cone of any homomorphism.

\begin{lemma} \label{th:reduction}
There is an exact functor $\Pi: K^b({\mathfrak P}) \longrightarrow D(\dga_{m,n})$ with the
following properties:
\begin{normallist}
\item \label{item:projective-modules}
$\Pi(\overline{P}_i)$ is isomorphic to ${\mathcal P}_i$ up to some shift;
\item \label{item:shift}
There is a canonical isomorphism of functors $\Pi \circ \{1\} \iso [-n] \circ \Pi$;
\item \label{item:double-grading}
The natural map, which exists in view of property {\em \ref{item:shift}},
\[
\bigoplus_{r_2 = nr_1} \Hom_{K^b({\mathfrak P})}(X,Y\{r_1\}[r_2]) \longrightarrow
\Hom_{D(\dga_{m,n})}(\Pi(X), \Pi(Y)),
\]
is an isomorphism for all $X,Y \in K^b({\mathfrak P})$.
\end{normallist}
\end{lemma}

\proof As a first step, consider the functor $\Pi': {\mathfrak P} \longrightarrow
Dgm(\dga_{m,n})$ defined as follows. The object $\overline{P}_i\{r\}$ goes to the dgm
${\mathcal P}_i[\sigma_i-nr]$, where $\sigma_i = -d_1-d_2-\dots-d_{i-1}$, and this is extended
to direct sums in the obvious way. Let $\overline{A}_m^d$ be the space of elements of degree
$d$ in $\overline{A}_m$. Homomorphisms of graded modules ${\overline P}_i\{r\} \longrightarrow
{\overline P}_j\{s\}$ correspond in a natural way to elements of
$\barpath{j}\overline{A}_m^{r-s}\barpath{i}$. On the other hand, dgm homomorphisms between
${\mathcal P}_i[\sigma_i-nr]$ and ${\mathcal P}_j[\sigma_j-ns]$ correspond to elements of
degree $\sigma_j-\sigma_i - n(s-r)$ in $(j)A_{m,n}(i)$. There is an obvious isomorphism, for
any $1 \leq i,j \leq m$ and $d \in \Z$,
\begin{equation} \label{eq:change-grading}
\barpath{j}\overline{A}_m^d\barpath{i} \iso (j)A_{m,n}^{\sigma_j-\sigma_i+nd}(i)
\end{equation}
which sends any basis element in \eqref{eq:bar-basis} of the form $\barpath{i_0|\dots|i_\nu}$
to the corresponding element $(i_0|\dots|i_\nu) \in A_{m,n}$; one needs to check, case by
case, that the degrees turn out right. We use \eqref{eq:change-grading} to define $\Pi'$ on
morphisms; this is obviously compatible with composition, so that the outcome is indeed a
functor. Note that $\Pi' \circ \{1\} \iso [-n] \circ \Pi'$.

Now take a finite chain complex in ${\mathfrak P}$. Applying $\Pi'$ to each object in the
complex yields a chain complex in $Dgm(\dga_{m,n})$, which one can then collapse into a single
dgm. This procedure yields a functor $K^b({\mathfrak P}) \longrightarrow K(\dga_{m,n})$, which
is exact since it carries cones to cones. We define $\Pi$ to be the composition of this with
the quotient functor $K(\dga_{m,n}) \longrightarrow D(\dga_{m,n})$. Properties
\ref{item:projective-modules} and \ref{item:shift} are now obvious from the definition of
$\Pi'$. The remaining property \ref{item:double-grading} can be reduced, by repeated use of
the Five Lemma, to the case when $X = \overline{P}_i\{r\}$, $Y = \overline{P}_i\{s\}$; and
then it comes down to the fact that \eqref{eq:change-grading} is an isomorphism. \qed

Define exact functors $\bar{t}_1,\dots,\bar{t}_m$ from $D^b(\Amod)$ to itself by
\begin{equation} \label{eq:bar-t}
\bar{t}_i(X) = \{ X\barpath{i} \otimes_k \overline{P}_i \longrightarrow X \}.
\end{equation}
Here $X\barpath{i}$ is considered as a complex of graded $k$-vector spaces; tensoring with
$\overline{P}_i$ over $k$ makes this into a complex of graded $\overline{A}_m$-modules; and
the arrow is the multiplication map. We can now state the results of \cite{khovanov-seidel98}.

\begin{lemma} \label{th:ks-one}
$\bar{t}_1,\dots,\bar{t}_m$ are exact equivalences and generate a weak braid group action
$\bar{\rho}_m: B_{m+1} \longrightarrow \Auteq(D^b(\Amod))$.
\end{lemma}

\begin{lemma} \label{th:ks-two}
For $g \in B_{m+1}$, let $f \in \Diff(D,\partial D;\Delta)$ be a diffeomorphism in the isotopy
class corresponding to $g$, and $\overline{R}_m^g$ a functor which represents
$\bar{\rho}_m(g)$. Then
\[
 \sum_{r_1,r_2} \dim_k \Hom_{D^b(\Amod)}(\overline{P}_i,
 \overline{R}_m^g(\overline{P}_j)\{r_1\}[r_2]) =
 2\,I(b_i,f(b_j))
\]
for all $1 \leq i,j \leq m$.
\end{lemma}

Lemma \ref{th:ks-one} essentially summarizes the contents of \cite[section
3]{khovanov-seidel98}, and Lemma \ref{th:ks-two} is \cite[Theorem 1.1]{khovanov-seidel98}. The
notation here is slightly different (our $\overline{A}_m$, $\overline{P}_i$ and $\bar{t}_i$
are the $A_m$, $P_i$ and ${\mathcal R}_i$ of that paper). We have also modified the
definitions very slightly, namely, we use right modules instead of left modules as in
\cite{khovanov-seidel98}, and the coefficients are $k$ instead of $\Z$. These changes do not
affect the results at all (a very conscientious reader might want to check that inversion of
paths defines an isomorphism between $\overline{A}_m$ and its opposite, and that a result
similar to Lemma \ref{th:reduction} can be proved for an algebra $\overline{A}_m$ defined over
$\Z$).

\proof[Proof of Lemma \ref{th:gin}] Since the modules $\overline{P}_i$ are projective, the
obvious exact functor $K^b({\mathfrak P}) \longrightarrow D^b(\Amod)$ is full and faithful. To
save notation, we will consider $K^b({\mathfrak P})$ simply as a subcategory of $D^b(\Amod)$.
An inspection of \eqref{eq:bar-t} shows that the $\bar{t}_i$ preserve this subcategory, and
the same is true of their inverses, defined in \cite{khovanov-seidel98}. In other words, the
weak braid group action $\bar{\rho}_m$ restricts to one on $K^b({\mathfrak P})$. It follows
from the definition of $\Pi$ that $\Pi \circ \bar{t}_i|K^b({\mathfrak P}) \iso t_i \circ \Pi$.
Hence, if $\overline{R}^g_m$ and $R^g_{m,n}$ are functors representing $\bar{\rho}_m(g)$
respectively $\rho_{m,n}(g)$, the diagram
\[
\xymatrix{
 {K^b({\mathfrak P})} \ar[r]^-{\overline{R}^g_m} \ar[d]_{\Pi} &
 {K^b({\mathfrak P})} \ar[d]_{\Pi} \\
 {D(\dga_{m,n})} \ar[r]^-{R^g_{m,n}} &
 {D(\dga_{m,n})}
}
\]
commutes up to isomorphism. Using this, Lemma \ref{th:reduction}\ref{item:double-grading} and
Lemma \ref{th:ks-two}, one sees that
\begin{align*}
 & \textstyle{\sum_r} \dim_k \Hom_{D(\dga_{m,n})}({\mathcal P}_i, R_{m,n}^g({\mathcal P}_j)[r]) \\
 & = \textstyle{\sum_r} \dim_k \Hom_{D(\dga_{m,n})}(\Pi(\overline{P}_i), \Pi
 \overline{R}_m^g(\overline{P}_j)[r]) \\
 & = \textstyle{\sum_{r_1,r_2}} \dim_k \Hom_{D^b(\Amod)}(\overline{P}_i,
 \overline{R}_m^g(\overline{P}_j)\{r_1\}[r_2]) \\ & = 2\,I(b_i,f(b_j)). \qed
\end{align*}

\proof[Proof of Theorem \ref{th:a-faithful}] For $g \in B_{m+1}$, choose $f$ and $R^g_{m,n}$
as in Lemma \ref{th:gin}. Take also another element $g' \in B_{m+1}$ and correspondingly $f'$,
$R^{g'}_{m,n}$. Applying Lemma \ref{th:gin} to $(g')^{-1}g$ shows that
\begin{align*}
 & I(f'(b_i),f(b_j)) = I(b_i,(f')^{-1}f(b_j)) \\
 & = \half \textstyle{\sum_r} \dim_k \Hom({\mathcal P}_i,
 (R^{g'}_{m,n})^{-1}R^g_{m,n}({\mathcal P}_j)) \\
\intertext{and assuming that $R^g_{m,n}({\mathcal P}_j) \iso {\mathcal P}_j$ for all $j$,}
 & = \half \textstyle{\sum_r} \dim_k \Hom({\mathcal P}_i,
 (R^{g'}_{m,n})^{-1}({\mathcal P}_j)) \\
 & = I(b_i,(f')^{-1}(b_j)) = I(f'(b_i),b_j).
\end{align*}
Since $i$ and $f'$ can be chosen arbitrary, it follows from Lemma \ref{th:pairing} that
$f(b_j) \isotopic b_j$ for all $j$. Hence, by Lemma \ref{th:central-element}, $g =
(g_1g_2\dots g_m)^{\nu(m+1)}$ for some $\nu \in \Z$. But then $R^g_{m,n}({\mathcal P}_j) \iso
{\mathcal P_j}[\nu(2m - (m+1)n)]$ by Lemma \ref{th:central-action}. In view of the assumption
that $R^g_{m,n}({\mathcal P_j}) \iso {\mathcal P_j}$, this implies that $\nu = 0$, hence that
$g = 1$. \qed

\subsection{Conclusion\label{subsec:a-is-formal}}

The graded algebras $A_{m,n}$ are always augmented. For $n \geq 2$ they are even connected, so
that there is only one choice of augmentation map. This makes it possible to apply Theorem
\ref{th:formality}.

\begin{lemma} \label{th:a-is-formal}
$A_{m,n}$ is intrinsically formal for all $m,n \geq 2$.
\end{lemma}

The proof is by a straight computation of Hochschild cohomology (it would be
nice to have a more conceptual explanation of the result). Its difficulty
depends strongly on the parameter $n$. The easy case is when $n>2$, since then
already the relevant Hochschild cochain groups are zero; this is no longer true
for $n = 2$. At first sight the computation may appear to rely on our specific
choice \eqref{eq:d-i} of degrees $d_i$, but in fact this only serves to
simplify the bookkeeping: the Hochschild cohomology remains the same for any
other choice. Throughout, we will write $\Gamma,A$ instead of
$\Gamma_{m,n},A_{m,n}$.

\proof[Proof for $n > 2$] Note that the `degree' label on any edge of $\Gamma$ is $\geq
[n/2]$. Moreover, the labels on any two consecutive edges add up to $n$. These two facts imply
that the degree of any nonzero path $(i_0|\dots|i_l)$ of length $l$ in $k[\Gamma]$ is $\geq
[(nl)/2]$. Now, any element of $(A^+)^{\otimes_R q}$ can be written as a sum of expressions of
the form
\[
c = (i_{1,0}|\dots|i_{1,l_1}) \otimes (i_{2,0}|\dots|i_{2,l_2}) \otimes \dots \otimes
(i_{q,0}|\dots|i_{q,l_q}),
\]
with all $l_q > 0$. Because the tensor product is over $R$, such a $c$ can be nonzero only if
the paths $(i_{\nu,0}|\dots|i_{\nu,l_\nu})$ match up, in the sense that $i_{\nu,l_\nu} =
i_{\nu+1,0}$. Then, using the observation made above, one finds that
\[
\deg(c) = \deg (i_{1,0}|\dots|i_{1,l_1}|i_{2,1}|\dots|i_{2,l_2}|i_{3,1}|\dots|i_{q,l_q}) \geq
[n(l_1 + \dots + l_q)/2].
\]
Hence $(A^+)^{\otimes q}$ is concentrated in degrees $\geq [(nq)/2]$. On the other hand,
$A[2-q]$ is concentrated in degrees $\leq n + q - 2$, which implies that
\[
C^q(A,A[2-q]) = \Hom_{R-R}((A^+)^{\otimes_R q},A[2-q]) = 0 \quad \text{if $n \geq 4$ or $q
\geq 4$.}
\]
We will now focus on the remaining case $(n,q) = (3,3)$. Then $(A^+)^{\otimes_R 3}$ is
concentrated in degrees $\geq 4$ while $A[-1]$ is concentrated in degrees $\leq 4$. The degree
four part of $(A^+)^{\otimes_R 3}$ is spanned by elements $c = (i_0|i_1) \otimes (i_1|i_2)
\otimes (i_2|i_3)$, which obviously satisfy $i_3 \neq i_0$. It follows that as an
$R$-bimodule, the degree four part satisfies $e_i ((A^+)^{\otimes_R 3})^4 e_i = 0$. On the
other hand, the degree four part of $A[-1]$ is spanned by the elements $(i|i\pm 1|i)$, so it
satisfies $e_i A[-1]^4 e_j = 0$ for all $i \neq j$. This implies that there can be no nonzero
$R$-bimodule maps between $(A^+)^{\otimes_R 3}$ and $A[-1]$, and hence that $C^3(A,A[-1])$ is
after all trivial.

\proof[Proof for $n = 2$] Consider the relevant piece of the Hochschild complex,
\[
C^{q-1}(A,A[2-q]) \stackrel{\partial^{q-1}}{\longrightarrow} C^q(A,A[2-q])
\stackrel{\partial^q}{\longrightarrow} C^{q+1}(A,A[2-q]).
\]
$C^{q+1}(A,A[2-q])$ is zero for degree reasons. In fact, since all edges in $\Gamma$ have
`degree' labels one, paths are now graded by their length, so that $(A^+)^{\otimes_R q+1}$ is
concentrated in degrees $\geq q+1$, while $A[2-q]$ is concentrated in degrees $\leq q$. In
contrast $C^q(A,A[2-q])$ is nonzero for all even $q$. To give a more precise description of
this group we will use the basis of $A$ from \eqref{eq:a-basis}, and the basis of
$(A^+)^{\otimes_R q}$ derived from that. Let $(i_0|\dots|i_q)$, $i_q = i_0$, be a closed path
of length $q$ in $\Gamma$. Define $\phi_{i_0,\dots,i_m} \in C^q(A,A[2-q])$ by setting
\[
\phi_{i_0,\dots,i_q}(c) = \begin{cases}
 (i_0|i_0\!\pm\! 1|i_0) & \text{if } c = (i_0|i_1) \otimes \dots \otimes (i_{q-1}|i_q), \\
 0 & \text{on all other basis elements $c$.}
 \end{cases}
\]
We claim that the elements defined in this way, with $(i_0|\dots|i_q)$ ranging over all closed
paths, form a basis of $C^q(A,A[2-q])$. To prove this, note that there is only one degree,
which is $q$, where both $(A^+)^{\otimes q}$ and $A[2-q]$ are nonzero. The degree $q$ part of
$(A^+)^{\otimes q}$ is spanned by expressions $c = (i_0|i_1) \otimes \dots \otimes
(i_{q-1}|i_q)$, with $i_q$ not necessarily equal to $i_0$. The degree $q$ part of $A[2-q]$ is
spanned by elements $(i|i\pm 1|i)$. Hence, an argument using the $R$-bimodule structure shows
that if $i_q \neq i_0$, then $\phi(c) = 0$ for all $\phi \in C^q(A,A[2-q])$. This essentially
implies what we have claimed.

We now turn to $C^{q-1}(A,A[2-q])$; for this group we will not need a complete description,
but only some sample elements. Given a closed path $(i_0|\dots|i_q)$ as before in $\Gamma$, we
define $\phi' \in C^{q-1}(A,A[2-q])$ by setting $\phi'(c) = (i_0|i_{q-1})$ if $c = (i_0|i_1)
\otimes \dots (i_{q-2}|i_{q-1})$, and zero on all other basis elements $c$. A simple
computation shows that $\delta^{q-1}(\phi') = -\phi_{i_0,\dots,i_q} -
\phi_{i_{q-1},i_0,i_1,\dots,i_{q-1}}$. Also, for any closed path $(i_0|\dots|i_q)$ with $i_2 =
i_0$ and $i_1 = i_0+1$, define $\phi'' \in C^{q-1}(A,A[2-q])$ by setting $\phi''(c) = (i_0|i_0
\pm 1|i_0)$ for $c = (i_0|i_1|i_2) \otimes (i_2|i_3) \otimes \dots \otimes (i_{q-1}|i_q)$, and
again zero for all other basis elements $c$. Then $\delta^{q-1}(\phi'')$ is equal to
$-\phi_{i_0,i_1,\dots,i_q} - \phi_{i_0,i_1-2,i_2,\dots,i_q}$ for $i_0>1$, and to
$-\phi_{i_0,i_1,\dots,i_q}$ for $i_0 = 1$.

To summarize, we have now established that the following relations hold in $HH^q(A,A[2-q])$:
\begin{normallist}
\item \label{item:cyclic-relation}
$[\phi_{i_0,\dots,i_q}] = -[\phi_{i_{q-1},i_q,i_1,\dots,i_{q-1}}]$ for all closed paths
$(i_0|\dots|i_q)$ in the quiver $\Gamma$.
\item \label{item:kink-relation}
$[\phi_{i_0,\dots,i_q}] = -[\phi_{i_0,i_1-2,i_2,\dots,i_q}]$ whenever $i_0 = i_2 \geq 2$ and
$i_1 = i_0+1$.
\item \label{item:zero-relation}
$[\phi_{i_0,\dots,i_q}] = 0$ whenever $i_0 = i_2 = 1$ and $i_1 = 2$.
\end{normallist}
Take an arbitrary element $\phi_{i_0,\dots,i_q}$. By applying \ref{item:cyclic-relation}
repeatedly, one can find another element $\phi_{i_0',\dots,i_q'}$ which represents the same
Hochschild cohomology class, up to a sign, and such that $i_1'$ is maximal among all $i_\nu'$.
This implies that $i_0' = i_2' = i_1' - 1$. If $i_1' = 2$ then we can apply
\ref{item:zero-relation} to show that our Hochschild cohomology class is zero. Otherwise pass
to $\phi_{i_0',i_1'-2,\dots,i_q'}$, which represents the same Hochschild cohomology class up
to sign due to \ref{item:kink-relation}, and repeat the argument. The iteration terminates
after finitely many moves, because the sum of the $i_\nu$ decreases by two in each step. Hence
$HH^q(A,A[2-q])$ is zero for all $q \geq 1$. \qed

\proof[Proof of Theorem \ref{th:main}] We first need to dispose of the trivial case $m = 1$.
In that case, choose a resolution $F_1 \in {\mathfrak K}$ of $E_1$. Pick a nonzero morphism
$\phi: F_1 \longrightarrow F_1[n]$. This, together with $\id_{F_1}$, determines an isomorphism
of graded vector spaces $\Hom^*(F_1,F_1) \iso k \oplus k[-n]$, and hence an isomorphism in
${\mathfrak K}$ between $F_1 \oplus F_1[-n]$ and $\Hom^*(F_1,F_1) \otimes F_1$. Consider the
commutative diagram
\[
\xymatrix{
 {T_{F_1}(F_1)[-1]} \ar[r] &
 {\Hom^*(F_1,F_1) \otimes F_1} \ar[r]^-{ev} &
 {F_1} \\
 {F_1[-n]} \ar[r]^-{(-\phi,\id)} &
 {F_1 \oplus F_1[-n]} \ar[u]^-{\iso} \ar[r]^-{(\id,\phi)} &
 {F_1} \ar[u]^-{\id}
}
\]
The upper row is a piece of the exact triangle which comes from the definition
of $T_{F_1}$ as a cone, and the lower row is obviously also a piece of an exact
triangle. By the axioms of a triangulated category, the diagram can be filled
in with an isomorphism between $F_1[-n]$ and $T_{F_1}(F_1)[-1]$. Transporting
the result to $D^b(\ab')$ yields $T_{E_1}(E_1) \iso E_1[1-n]$. Since $n \geq 2$
by assumption, it follows that $T_{E_1}^r(E_1) \not\iso E_1$ unless $r = 0$.

From now on suppose that $m \geq 2$. After shifting each $E_i$ by some amount,
we may assume that $\Hom^*(E_{i+1},E_i)$ is concentrated in degree $d_i$ for $i
= 1,\dots,m-1$ (shifting will not affect the statement because $T_{E_i[j]}$ is
isomorphic to $T_{E_i}$ for any $j \in \Z$). Choose resolutions
$E_1',\dots,E_m' \in {\mathfrak K}$ for $E_1,\dots,E_m$. Lemma
\ref{th:end-algebra} shows that the endomorphism dga $end(E')$ has $H(end(E'))
\iso A_{m,n}$. By Lemma \ref{th:a-is-formal}, $end(E')$ must be
quasi-isomorphic to ${\mathcal A}_{m,n}$. Define an exact functor $\Psi$ to be
the composition
\[
D^b(\ab') \stackrel{\iso}{\longleftarrow} {\mathfrak K} \stackrel{\Psi_{E'}}{\longrightarrow}
D(end(E')) \stackrel{\iso}{\longrightarrow} D(\dga_{m,n}).
\]
The first arrow is the standard equivalence, and the last one is the equivalence induced by
some sequence of dgas and quasi-isomorphisms. By construction $\Psi(E_i) \iso {\mathcal P}_i$
for $i = 1,\dots,m$. In the diagram
\[
\xymatrix{
 {D^b(\ab')} \ar[d]_{T_{E_i}} &
 {{\mathfrak K}} \ar[r]^-{\Psi_{E'}} \ar[l]_-{\iso} \ar[d]_{T_{E_i'}} &
 {D(end(E'))} \ar[r]^-{\iso} \ar[d]_{t_i} &
 {D(\dga_{m,n})} \ar[d]_{t_i} \\
 {D^b(\ab')} &
 {{\mathfrak K}} \ar[r]^-{\Psi_{E'}} \ar[l]_-{\iso} &
 {D(end(E'))} \ar[r]^-{\iso} &
 {D(\dga_{m,n})}
}
\]
the first square commutes because that is the definition of $T_{E_i}$, the second square by
Lemma \ref{th:twists-and-psis}, and the third one by Lemma \ref{th:twists-and-quasiisos}. Now
let $g$ be an element of $B_{m+1}$, $R^g: D^b(\ab') \longrightarrow D^b(\ab')$ a functor which
represents $\rho(g)$, and $R_{m,n}^g: D(\dga_{m,n}) \longrightarrow D(\dga_{m,n})$ a functor
which represents $\rho_{m,n}(g)$. By applying the previous diagram several times one sees that
\[
R_{m,n}^g \circ \Psi \iso \Psi \circ R^g.
\]
Assume that $R_g(E_i) \iso E_i$ for all $i$; then also $R_{m,n}^g({\mathcal
P}_i) = R_{m,n}^g\Psi(E_i) \iso \Psi R^g(E_i) \iso \Psi(E_i) \iso {\mathcal
P}_i$. By Theorem \ref{th:a-faithful} it follows that $g$ must be the identity.
\qed

We have not tried to compute the Hochschild cohomology of $A_{m,n}$ for $n = 1$. However, an
indirect argument using the non-faithful $B_4$-action of section \ref{subsec:elliptic-curve}
shows that $A_{3,1}$ cannot be intrinsically formal. More explicitly, if one takes the sheaves
$\O_x, \O, \O_y$ used in that example, and chooses injective resolutions by quasi-coherent
sheaves for them, then the resulting dga $end(E')$ is not formal. One can give a more direct
proof of the same fact by using essentially the same Massey product computation as Polishchuk
in \cite[p.\ 3]{polishchuk98}.

\bibliographystyle{amsplain}

\providecommand{\bysame}{\leavevmode\hbox to3em{\hrulefill}\thinspace}

\end{document}